\documentclass[a4paper,11pt,reqno]{amsart}
\usepackage{graphics}
\newtheorem{theo}{Theorem}
\newtheorem{lem}{Lemma}
\newtheorem{cor}{Corollary}
\newtheorem{rem}{Remark}
\newtheorem{ex}{Example}

\newcommand{\mj}{\operatorname{STAND\_SPLIT}}
\newcommand{\jt}{\operatorname{JT}}
\newcommand{\cjt}{\operatorname{CJT}}
\newcommand{\rjt}{\operatorname{RJT}}
\newcommand{\crjt}{\operatorname{CRJT}}
\newcommand{\shift}{\operatorname{SHIFT}}
\newcommand{\rshift}{\operatorname{RSHIFT}}
\newcommand{\ssplit}{\operatorname{STAND\_SPLIT}}
\newcommand{\hsplit}{\operatorname{HOOK\_SPLIT}}
\newcommand{\smerge}{\operatorname{MERGE}}

\newcommand{\trans}{\operatorname{TRANS}}
\newcommand{\rtrans}{\operatorname{RTRANS}}
\newcommand{\js}{\operatorname{JS}}
\newcommand{\rjs}{\operatorname{RJS}}

\begin{document}

\title[The shifted hook-length formula]
{A bijective proof of the hook-length formula for shifted standard tableaux}

\author[Ilse Fischer]{\box\Adr}

\newbox\Adr
\setbox\Adr\vbox{
\centerline{ \large Ilse Fischer}
\vspace{0.3cm}
\centerline{Institut f\"ur Mathematik der Universit\"at Klagenfurt,}
\centerline{Universit\"atsstrasse 65-67, A-9020 Klagenfurt, Austria.}
\centerline{E-mail: {\tt Ilse.Fischer@uni-klu.ac.at}}
}

\date{}

\begin{abstract}
We present a bijective proof of the hook-length formula for shifted standard 
tableaux of a fixed shape based on a modified jeu de taquin and the ideas of the bijective proof of the hook-length
formula for ordinary standard tableaux by Novelli, Pak and Stoyanovskii
\cite{NPS}. In their proof Novelli, Pak and Stoyanovskii define a bijection 
between arbitrary fillings of the Ferrers diagram with the integers $1,2,\dots,n$ 
and pairs of standard tableaux and hook tabloids. In our shifted version of
their algorithm the map from the set of arbitrary
fillings of the shifted Ferrers diagram onto the set of shifted standard tableaux is 
analog to the construction of Novelli, Pak and Stoyanovskii, however, unlike to 
their algorithm, we are forced to use the 'rowwise' total order of the cells
in the shifted Ferrers diagram rather than the 'columnwise' total order as the underlying
order in the algorithm.
Unfortunately the
construction of the shifted hook tabloid is more complicated in the shifted case.
As a side-result we obtain a simple random algorithm for generating shifted
standard tableaux of a given shape, which produces every such tableau equally
likely.
\end{abstract}

\maketitle

\section{Introduction}
\begin{figure}
\setlength{\unitlength}{1cm}
\begin{picture}(14.5,4)
\put(1,1){\scalebox{0.35}{\includegraphics{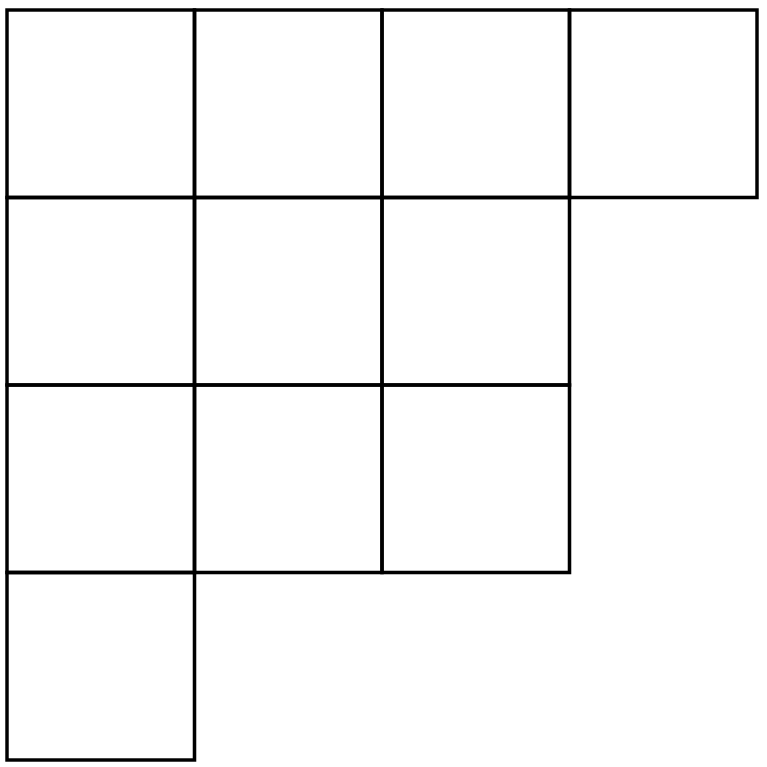}}}
\put(8,1){\scalebox{0.35}{\includegraphics{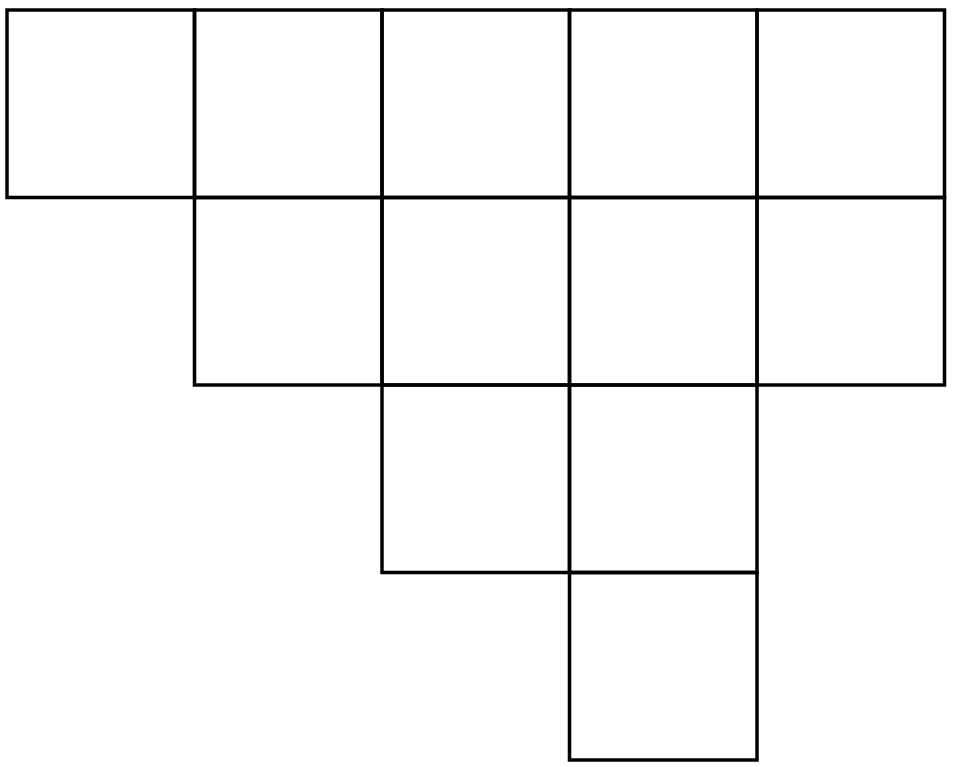}}}
\put(1,0.5){a. The Ferrers diagram}
\put(1,0){corresponding to $(4,3,3,1)$}
\put(8,0.5){b. The shifted Ferrers diagram}
\put(8,0){corresponding to $(5,4,2,1)$} 
\end{picture}
\caption{}
\label{ferrer}   
\end{figure}

A partition of a positive integer $n$ is a sequence of positive integers
$\lambda=(\lambda_1,\lambda_2,\dots,\lambda_r)$ with 
$\lambda_1+\lambda_2+\dots+\lambda_r=n$ and 
$\lambda_1 \ge \lambda_2 \ge \dots \ge \lambda_r$. 
The {\it (ordinary) Ferrers diagram} of shape $\lambda$ 
is an array of cells with $r$ left-justified rows
and $\lambda_i$ cells in row $i$. Figure~\ref{ferrer}.a
shows the Ferrers diagram corresponding to $(4,3,3,1)$.

If $\lambda$ is a partition with distinct components (strict partition)
then the {\it shifted Ferrers diagram} of shape $\lambda$ is 
an array of cells with $r$ rows, each row indented by one 
cell to the right with respect to the previous row and 
$\lambda_i$ cells in row $i$. Figure~\ref{ferrer}.b 
shows the shifted Ferrers diagram corresponding to 
$(5,4,2,1)$.

\begin{figure}
\setlength{\unitlength}{1cm}
\begin{picture}(14.5,4)
\put(1,1){\scalebox{0.35}{\includegraphics{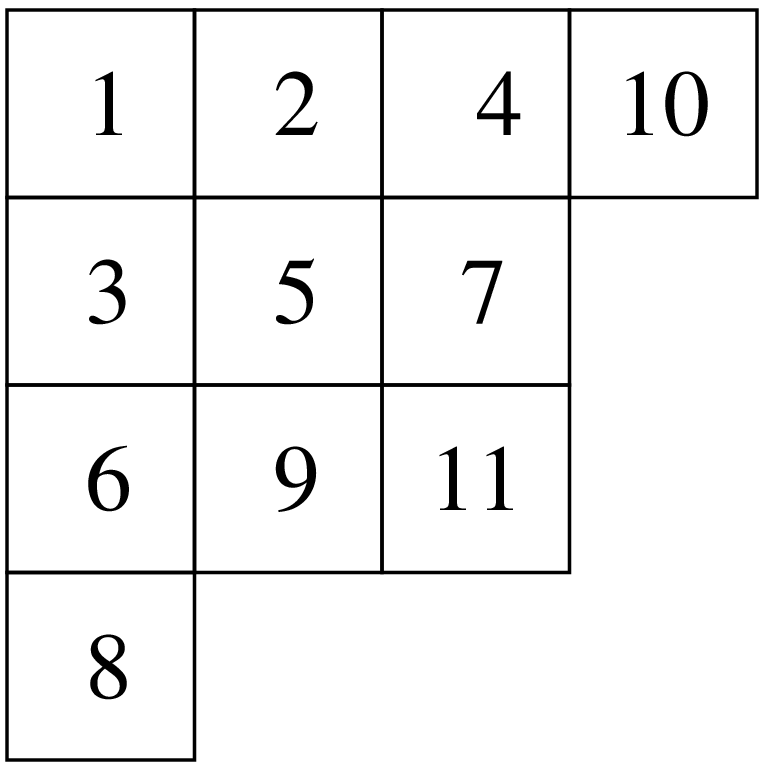}}}
\put(8,1){\scalebox{0.35}{\includegraphics{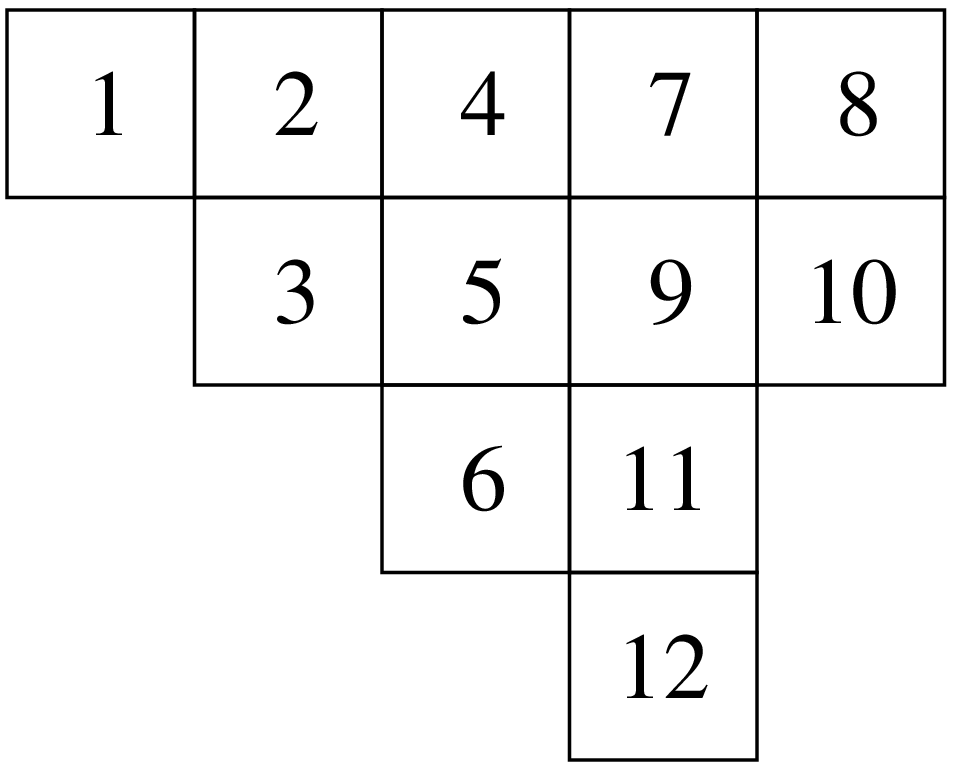}}}
\put(1,0.5){a. A standard tableau}
\put(8,0.5){b. A shifted standard}
\put(8,0){tableau} 
\end{picture}
\caption{}
\label{standard}   
\end{figure}

Given a partition $\lambda$ of $n$, respectively strict partition of $n$, a {\it standard
tableau}, respectively a {\it shifted standard tableau} of shape $\lambda$, 
is a filling of the cells of 
the ordinary Ferrers diagram, respectively shifted Ferrers diagram, of shape $\lambda$ with $1,2,\dots,n$,
such that the entries along rows and columns are increasing.
Figure~\ref{standard}.a displays an example of a standard tableau 
of shape $(4,3,3,1)$ and Figure~\ref{standard}.b displays an example 
of a shifted standard tableau of shape $(5,4,2,1)$,

\smallskip

Once we have accepted these definitions it is a natural question to 
ask for the number of standard tableaux, respectively shifted 
standard tableaux, of a given shape $\lambda$. Surprisingly there 
exists a simple product formula for these  numbers. 
It involves objects called {\it hooks}, which are defined in the 
following paragraph.

\begin{figure}
\setlength{\unitlength}{1cm}
\begin{picture}(14.5,4)
\put(1,1){\scalebox{0.35}{\includegraphics{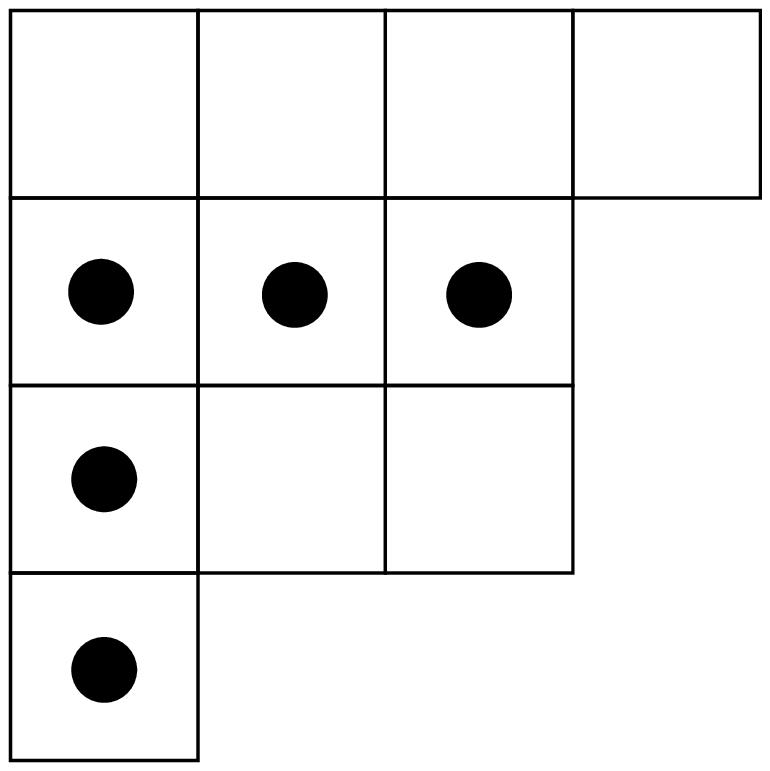}}}
\put(8,1){\scalebox{0.35}{\includegraphics{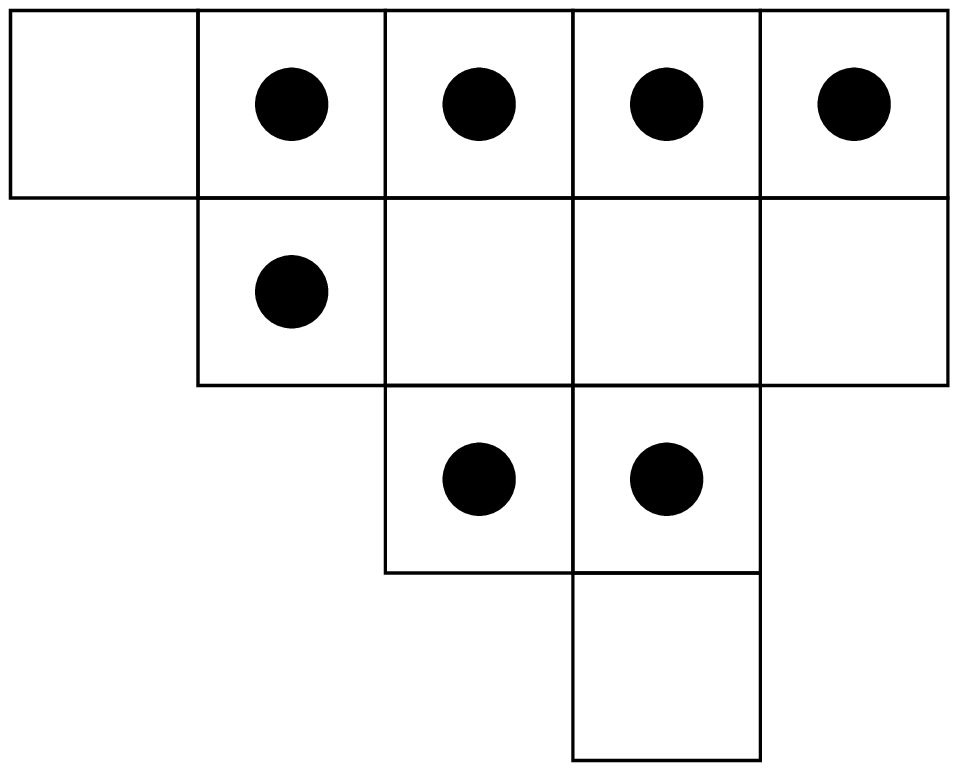}}}
\put(1,0.5){a. The hook of $(2,1)$}
\put(8,0.5){b. The shifted hook of $(1,2)$} 
\end{picture}
\caption{}
\label{hook}   
\end{figure}

We label the cell in the $i$-th row and $j$-th column
of the ordinary, respectively shifted, Ferrers diagram 
of shape $\lambda$ by the pair $(i,j)$. 
The {\it hook} of a cell $(i,j)$ in an {\it ordinary} Ferrers diagram 
is the set of cells that are either in the same row
as $(i,j)$ and to the right of $(i,j)$, or in the 
same column as $(i,j)$ and below $(i,j)$, $(i,j)$ included.
The dots in Figure~\ref{hook}.a indicate the hook 
of the cell $(2,1)$.
The {\it hook} of a cell $(i,j)$ in a {\it shifted} Ferrers diagram
again includes all cells that are either in the same 
row as $(i,j)$ and to the right of $(i,j)$, or in the 
same column as $(i,j)$ and below $(i,j)$, $(i,j)$ included,
but if this set contains the cell $(j,j)$ on the main 
diagonal, then also the cells of the 
$(j+1)$-st row belong to the hook of $(i,j)$. 
The dots in Figure~\ref{hook}.b indicate the hook 
of cell $(1,2)$. The {\it hook-length} $h_{i,j}$ 
of the cell $(i,j)$ is the number of cells in 
the hook of $(i,j)$. 

Now we are in the position to state the hook-length formula.

\begin{theo}[\cite{frt},\cite{gan}]
The number of standard tableaux, respectively shifted 
standard tableaux, of shape $\lambda$ is 
$$
\frac{n!}{\prod_{(i,j)} h_{i,j}},
$$ 
where the product in the denominator is taken 
over all cells in the Ferrers diagram, respectively shifted 
Ferrers diagram, of shape $\lambda$.
\end{theo}

Thus the number of standard tableaux of shape $(4,3,3,1)$ is 
$\frac{11!}{7 \cdot 5 \cdot 4 \cdot 5 \cdot 3 \cdot 2 \cdot 4 \cdot 2}=1188$
and the number of shifted standard tableaux of shape $(5,4,2,1)$ is 
$\frac{12!}{9 \cdot 7 \cdot 6 \cdot 5 \cdot 2 \cdot 6 \cdot 5 \cdot 4 \cdot 3
  \cdot 2}=176$.

For the history of the various proofs of the hook-length formula for ordinary standard tableaux
see \cite[page 400]{stan}. A recipe for an inductive proof of the hook-length 
formula for shifted standard tableaux can be found in \cite[page 266]{macd}. 
In \cite{sagan} it is shown that the nice probabilistic proof of the hook-length formula 
for ordinary standard tableaux in \cite{green} has an analog for shifted standard 
tableaux.    

\medskip

A majority among the combinatorialists considers a bijective proof as the 
most aesthetic type of proof for an enumeration result. In \cite{NPS} a 
bijective proof of the hook-length formula for 
ordinary standard tableaux based on a modified jeu de taquin is given. 
There exists a bijective proof of the hook-length formula 
for shifted standard tableaux as well \cite{kratt}, however, 
it makes use of the involution principle by Garsia and Milne.
The aim of this paper is to
present an involution principle-free bijective proof of the hook-length formula for 
shifted standard tableaux, which is in the spirit of 
the beautiful bijective proof in \cite{NPS}.

\medskip

If we discover two sets ${\mathcal S},{\mathcal T}$ of combinatorial objects with the same cardinality, 
we often believe that this fact is a projection of a canonical bijection between the two sets. 
Such a bijection is called a bijective proof of the equality 
$|{\mathcal S}|=|{\mathcal T}|$. Now suppose we are in the following more general situation: 
There exists an integer $h$ such that $h \cdot |{\mathcal S}| = |{\mathcal T}|$. Then this fact 
could be a projection of a canonical $h$ to $1$ surjection from ${\mathcal T}$ onto 
${\mathcal S}$, i.e. a map from ${\mathcal T}$ onto ${\mathcal S}$ where every 
element in ${\mathcal S}$ is assigned to exactly $h$ elements in ${\mathcal T}$.
In the following $S_{\lambda}$ denotes the set of shifted standard
tableaux of shape $\lambda$. A {\it shifted tabloid} of shape 
$\lambda$ is an (arbitrary) filling of the cells of the shifted Ferrers diagram 
of shape $\lambda$ with $1,2,\dots,n$. We denote the set of shifted tabloids of shape $\lambda$ by 
$T_\lambda$ and observe that its cardinality is $n!$.
In our main theorem (Theorem~\ref{main}) we  present a $\prod_{(i,j)} h_{i,j}$ to $1$ 
surjection 
from $T_{\lambda}$ onto $S_{\lambda}$, which is clearly a proof of the hook-length 
formula for shifted standard tableaux. We prove Theorem~\ref{main} by introducing the set of 
hook tabloids $H_\lambda$ with $|H_\lambda|=\prod_{(i,j)} h_{i,j}$ and extending the surjection to a bijection 
from $T_{\lambda}$ to $S_{\lambda} \times H_{\lambda}$. The corresponding surjection from the set of 
ordinary tabloids onto the set of standard tableaux is similar to the surjection in the shifted case, see
\cite{NPS}, the extension in the shifted case is however far more complicated
compared to the ordinary case.

\section{Modified jeu de taquin}
\label{jeu}

In this section we describe the ordering procedure which assigns to every 
`scrambled' shifted tabloid $T \in T_{\lambda}$ an `ordered'  shifted standard tableau $S
\in S_{\lambda}$. This map
has the property that the number of shifted tabloids which are mapped 
to a fixed shifted standard tableau is $\prod_{(i,j)} h_{i,j}$. The 
ordering procedure is based on a modified jeu de taquin, which we have to describe first. 

\medskip

{\bf Notation.}
Let $T$ be a shifted tabloid of shape $\lambda$. Then $T_{i,j}=T_{(i,j)}$ denotes  
the entry of cell $(i,j)$. If the cell $(i,j)$ does 
not exist in the shifted Ferrers diagram of shape $\lambda$, let 
$T_{i,j} = \infty$. For every $e$, $1 \le e \le n$, there exists a unique cell $(i,j)$ in 
the shifted Ferrers diagram of shape $\lambda$ with $T_{i,j}=e$. We define $c_T(e)=(i,j)$.
 
\medskip

\begin{figure}
\setlength{\unitlength}{1cm}
\begin{center}
\scalebox{0.35}{\includegraphics{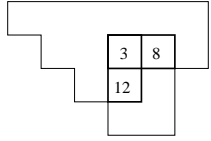}}
\end{center}
\caption{Entry $e=3$ is stable, for $3 \le 8, 12$.}
\label{stable}   
\end{figure}

\begin{figure}
\setlength{\unitlength}{1cm}
\begin{center}
\scalebox{0.35}{\includegraphics{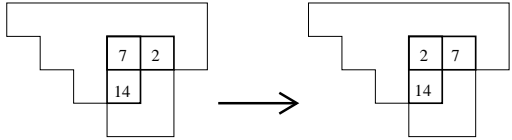}}
\end{center}
\caption{Entry $e=7$ is unstable in the left tabloid and thus $2=\min(14,2)$ and $7$ change place in the course of 
performing jeu de taquin with $7$. Thus $2$ is stable in the right tabloid. Next compare $7$ to 
its new neighbours to the right and below...until $7$ is either stable or the cell of $7$ 
is in $D$.}
\label{unstable}   
\end{figure}

\begin{figure}
\setlength{\unitlength}{1cm}
\begin{center}
\scalebox{0.30}{\includegraphics{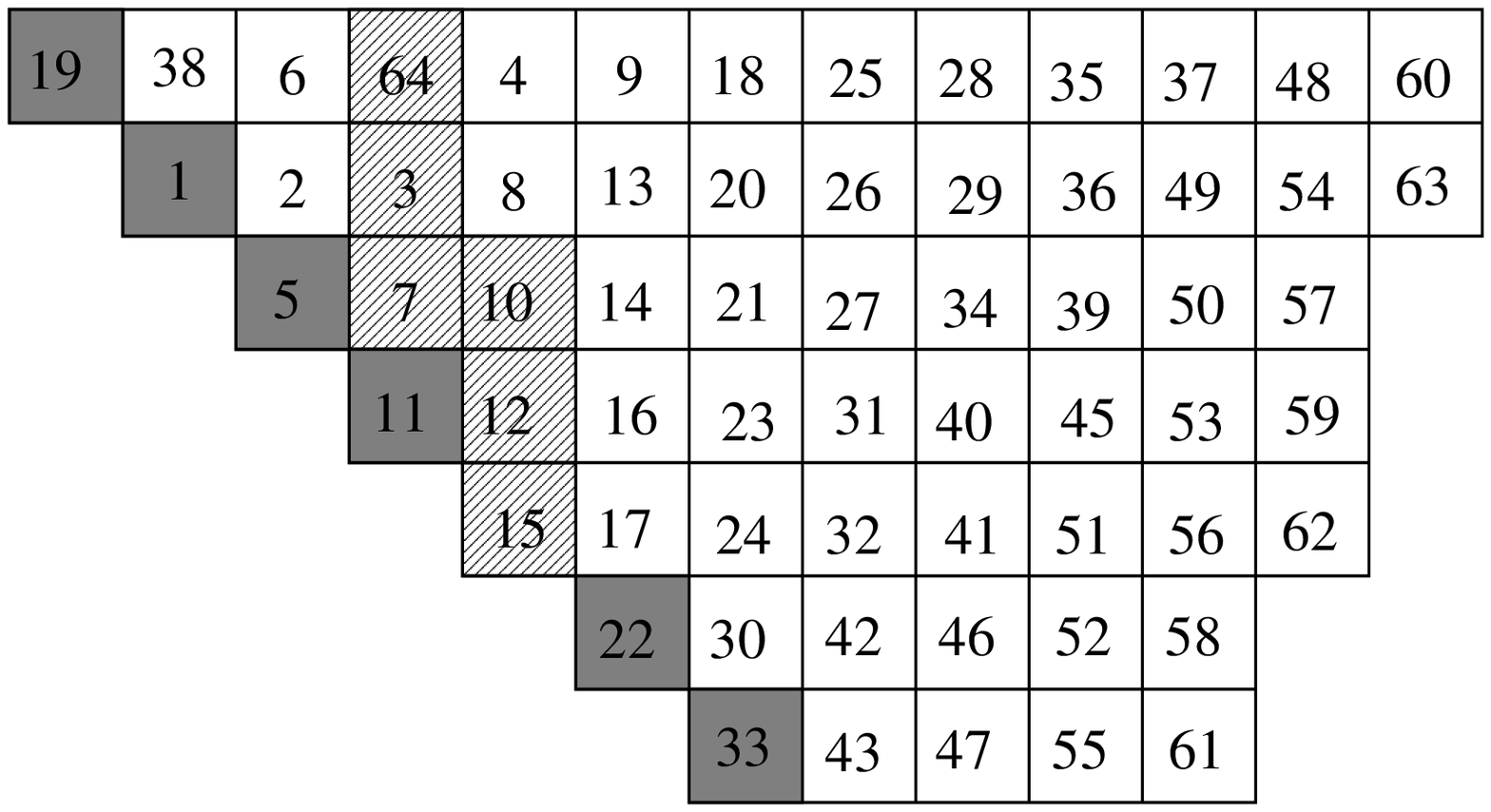}}
\hspace{1cm}
\scalebox{0.30}{\includegraphics{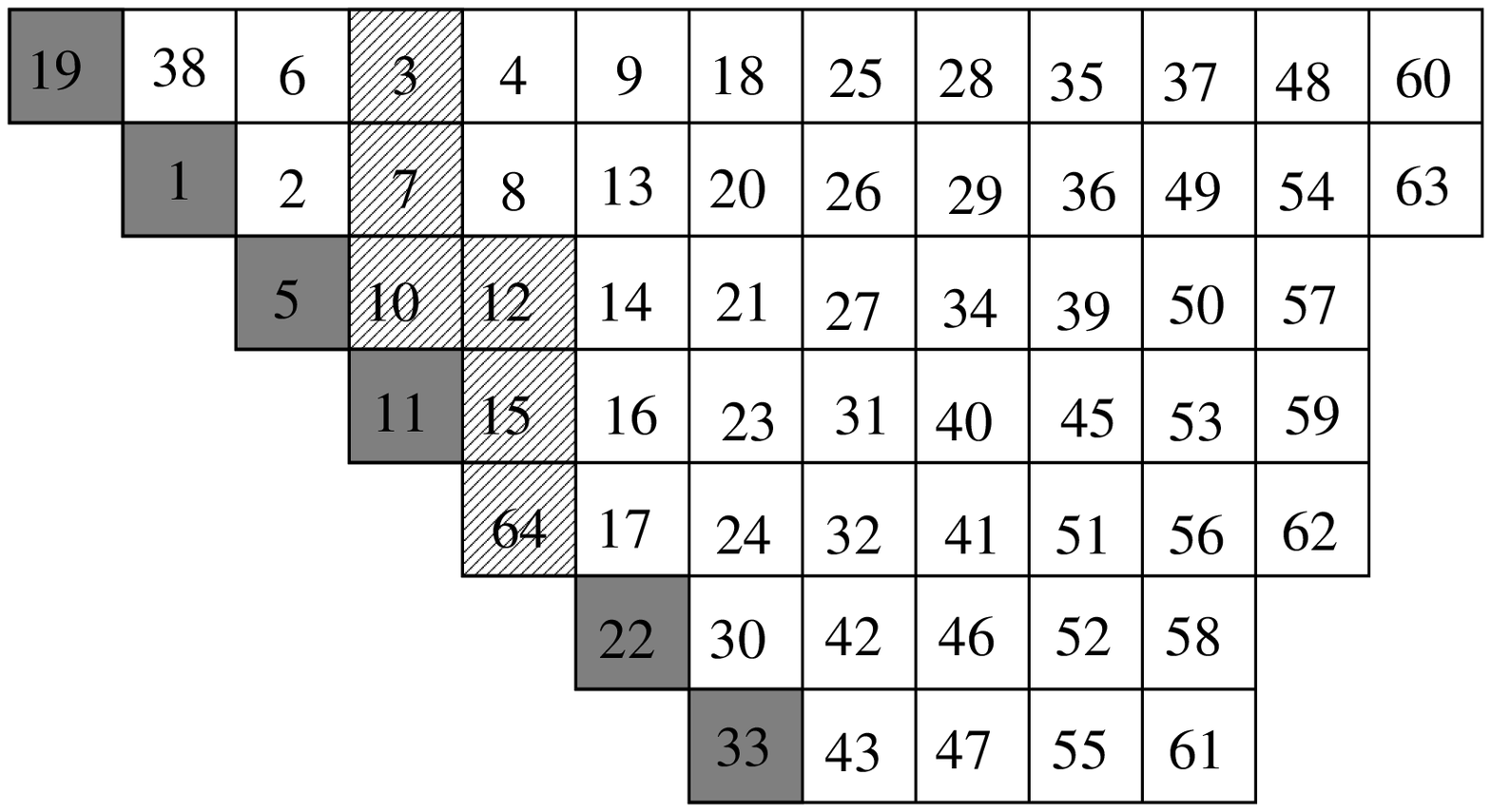}}
\end{center}
\caption{Performing jeu de taquin with $\rho=64$ and with respect to the cells on the main diagonal
in the tabloid on the left side
results in the tabloid on the right side.}
\label{jt}   
\end{figure}

{\bf Jeu de taquin in $T$ with entry $e$ and with respect to a set $D$.}
Let $1 \le e \le n$ and 
$D$ be a set of cells in the shifted Ferrers diagram of shape $\lambda$. We define the
routine {\it jeu de taquin in $T$ with entry $e$ and with 
respect to $D$} inductively. The output is another shifted 
tabloid $U$ of shape $\lambda$. Let $c_T(e)=(i,j)$. Set $U=T$ and stop if
either $(i,j) \in D$ or $T_{i,j} \le \min (T_{i+1,j}, T_{i,j+1})$. (In the 
latter case $e$ is said to be {\it stable}. See Figure~\ref{stable}.) 
Otherwise let $\rho \in \{ (i+1,j), (i, j+1) \}$ 
be such that $T_\rho = \min ( T_{i+1,j} , T_{i,j+1})$ and let $T'$ denote 
the tabloid we obtain by exchanging the entries $e$  and $T_\rho$ in $T$.
(See Figure~\ref{unstable}.) 
Next perform jeu de taquin in $T'$ with entry $e$ and with respect to $D$ 
in order to obtain the final tabloid $U$, i.e. we repeat this exchanging
procedure with $e$ and either its current neighbour to the right or below
until $e$ is either stable or in a cell of $D$.
(For an example see Figure~\ref{jt}.)

If $D$ is the empty set we omit `with respect to $D$'. The output 
tabloid $U$ is denoted by  $\jt_D(T,e)$ and the cell of $e$ in $U$ is denoted by 
$\cjt_D(T,e)$. If $D$ is the $k$-th row of the shifted Ferrers diagram we write 
$\jt_k(T,e)$ and $\cjt_k(T,e)$, respectively.

\smallskip

\begin{figure}
\setlength{\unitlength}{1cm}
\begin{center}
\scalebox{0.35}{\includegraphics{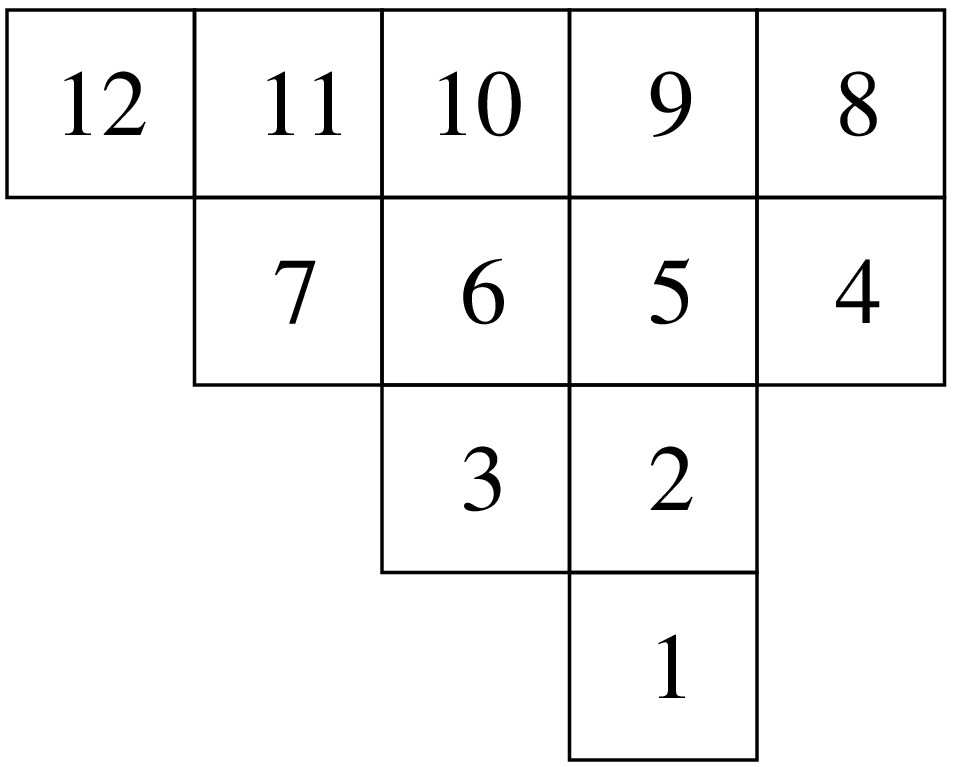}}
\end{center}
\caption{The total order.}
\label{ord}   
\end{figure}

{\bf The total order.}
We define a total order on the cells of a 
shifted Ferrers diagram of shape $\lambda$. This will be the order in which 
we perform the jeu de taquin just defined with the entries of a shifted tabloid.  A cell $\rho_1$ comes 
before cell $\rho_2$ if $\rho_1$ is in a lower row than 
$\rho_2$ or if both are in the same row but $\rho_1$ is to the 
right of $\rho_2$. Phrased differently, to obtain the total order 
one starts with the rightmost cell in the last row 
and reads each row from right to left, beginning with the
bottom row and continuing up to the first row.
Figure~\ref{ord} displays this total order for the shifted
Ferrers diagram of shape $(5,4,2,1)$. 

\medskip

{\bf The $\prod_{(i,j)} h_{i,j}$ to $1$ map from $T_\lambda$ onto $S_\lambda$.}
Let $T$ denote an arbitrary shifted tabloid of shape $\lambda$.
In order to construct the corresponding shifted standard 
tableau $S$ we perform step by step jeu de taquin with the entries of $T$
subject to the total order we have just defined and 
starting with the entry in the smallest cell.
To be more accurate: 

\medskip

\fbox{ \parbox{14cm}{ 

$\rho=$ smallest cell with respect to the total order

$S=T$

 Repeat 
 
$S=\jt(S,S_{\rho})$

$\rho$ = successor of $\rho$ in the total order

until $\rho=(1,1)$.

}}

\medskip

\begin{ex}
Consider the shifted tabloid 
$$
T= \begin{array}{ccccc} 11 & 4  &  9 & 8 & 1 \\
                           & 12 &  6 & 2 & 3 \\
                           &    & 10 & 5 &   \\
                           &    &    & 7 &   \\
    \end{array}
$$
of shape $(5,4,2,1)$. We construct the corresponding shifted standard 
tableau. According to the algorithm above we start with performing 
jeu de taquin with $7$, but $7$ is stable. The same is true for the entry 
$5$. Performing jeu de taquin with $10$ results in 
$
S =  \begin{array}{ccccc} 11 & 4  &  9 & 8 & 1 \\
                           & 12 &  6 & 2 & 3 \\
                           &    &  5 & 7 &   \\
                           &    &    & 10 &   \\
    \end{array}.
$
Since $3$ and $2$ are stable the next change happens to be when performing
jeu de taquin with entry $6$.
$
S =  \begin{array}{ccccc} 11 & 4  &  9 & 8 & 1 \\
                           & 12 &  2 & 3 & 6 \\
                           &    &  5 & 7 &   \\
                           &    &    & 10 &   \\
    \end{array}.
$
We perform jeu de taquin with $12$ --- 
$
S =  \begin{array}{ccccc} 11 & 4  &  9 & 8 & 1 \\
                           &   2 &  3 & 6 & 12 \\
                           &    &  5 & 7 &   \\
                           &    &    & 10 &   \\
    \end{array},
$
--- next with $8$ ---
$
S =  \begin{array}{ccccc} 11 & 4  &  9 & 1 & 8 \\
                           &   2 &  3 & 6 & 12 \\
                           &    &  5 & 7 &   \\
                           &    &    & 10 &   \\
    \end{array},
$
--- then with $9$ ---
$
S =  \begin{array}{ccccc} 11 & 4  &  1 & 6 & 8 \\
                           &   2 &  3 & 7 & 12 \\
                           &    &  5 & 9 &   \\
                           &    &    & 10 &   \\
    \end{array},
$
--- and with $4$ ---
$
S =  \begin{array}{ccccc} 11 & 1  &  3 & 6 & 8 \\
                           &   2 &  4 & 7 & 12 \\
                           &    &  5 & 9 &   \\
                           &    &    & 10 &   \\
    \end{array},
$
--- and finally with $11$ ---
$
S =  \begin{array}{ccccc} 1 & 2  &  3 & 6 & 8 \\
                           &   4 &  5 & 7 & 12 \\
                           &    &  9 & 10 &   \\
                           &    &    & 11 &   \\
    \end{array}.
$
\end{ex}

\bigskip

Observe that the output tabloid $S$ is a shifted standard tableau by construction. We denote it
by $\mj(T)$.  We are in the position to state our main theorem.
Note that the bijective proof in \cite{NPS} shows that
a similar theorem is true for ordinary standard tableaux.

\medskip

\begin{theo}
\label{main}
The map $T \to \mj(T)$ is a $\prod_{(i,j)} h_{i,j}$ to $1$ map from the 
set of shifted tabloids $T_\lambda$ onto the set of shifted standard tableaux $S_\lambda$.
\end{theo}

\medskip

As an interesting side-result we obtain a random algorithm which produces every shifted 
standard tableau of a given shape with the same probability.

\medskip

\begin{cor}
The following algorithm produces every 
shifted standard tableau of a given shape $\lambda$ with the same probability.
\begin{enumerate}
\item Generate a permutation $\pi$ of $\{1,2,\dots,n\}$ subject to uniform distribution.
\item Construct the corresponding shifted tabloid $T_\pi$ of shape $\lambda$ by filling the elements 
from $\pi$ into the shifted Ferrers diagram of shape $\lambda$ rowwise from top to bottom and 
in each row from left to right.
\item Apply $\ssplit(T_\pi)$ in order to obtain the shifted standard tableau.
\end{enumerate}
\end{cor}

\medskip

Note that Novelli, Pak and Stoyanovskii \cite{NPS}
use another total order of the cells in the Ferrers diagram, they perform jeu de taquin columnwise from right to left
and within a column from bottom to top. By 'transposing' their algorithm it is
clear that the order we defined in the shifted case would also induce a
$\prod_{i,j} h_{i,j}$ to $1$ map in 
the ordinary case. However, computer experiments with the strict partition  $(4,3,2,1)$
have shown that the order defined by Novelli, Pak and Stoyanovskii is not
admissible in the shifted case. Moreover it seems that the total order we
have defined is the only admissible in  the shifted case, whereas in
the ordinary case there exist many total orders with the property that they
induce a $\prod_{i,j} h_{i,j}$ to $1$ map. We plan to
discuss this phenomenon  in a forthcoming paper.

\medskip

\begin{figure}
\setlength{\unitlength}{1cm}
\begin{center}
\scalebox{0.50}{\includegraphics{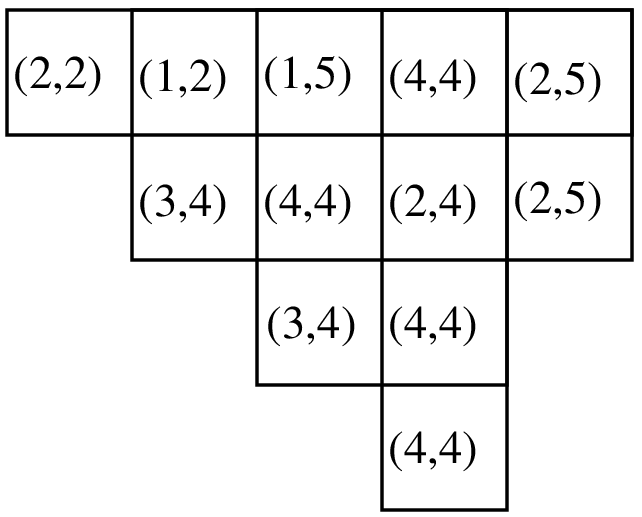}}
\end{center}
\caption{A shifted hook tabloid of shape $(5,4,2,1)$.}
\label{shbsp}   
\end{figure}

We prove Theorem~\ref{main} by giving a bijective proof of the 
hook-length formula. For that purpose we rewrite the hook-length formula as 
$$ n!= |S_{\lambda}| \cdot \prod_{(i,j)} h_{i,j}. $$
We define combinatorial objects that correspond to 
$\prod_{(i,j)} h_{i,j}$ in this formula in our bijective proof.
A {\it shifted hook tabloid} of shape $\lambda$ is a filling of the cells of the 
shifted Ferrers 
diagram of shape $\lambda$ with pairs of integers, such 
that the entry in a cell $\rho$ are the coordinates of a cell in 
the hook of $\rho$. (See Figure~\ref{shbsp}. This definition of 
a shifted hook tabloid has a natural analog for ordinary Ferrers diagram,
which is equivalent to the definition of a hook function in 
\cite{NPS}.) We denote the set of shifted hook tabloids by 
$H_\lambda$. Since $|H_\lambda| =  \prod_{(i,j)} h_{i,j}$
it suffices to find a bijection between the set of shifted tabloids $T_\lambda$ and 
the cartesian product of the set of  shifted standard tableaux $S_\lambda$
and the set of shifted hook tabloids $H_\lambda$ to prove the hook-length formula for 
shifted standard tableaux. 
$$
T_\lambda \stackrel{\text{bijection}}{\leftrightarrow} 
S_\lambda \times H_\lambda
$$
In order to prove Theorem~\ref{main} we construct such a bijection, where 
the shifted standard tableau is obtained from the shifted tabloid by the modified
jeu de taquin described above.
Unfortunately the second component of the bijection, i.e. the construction of the 
shifted hook tableau, is complicated.
Thus it would be nice to find an easier-to-describe 
bijection. Also a shorter proof of Theorem~\ref{main} would be of interest. 
Maybe the present paper serves as an inspiration in this task.  

\bigskip

\medskip

We give some further definitions we need for the rest of the paper.

\smallskip

{\bf Reverse jeu de taquin in $U$ with entry $e$ and with respect to a set $D$.}
We define an inverse to jeu de taquin, which we need to construct the inverse
of the map which assigns a pair of a shifted standard tableau and a
shifted hook tabloid to a given shifted tabloid. In all our algorithms where we use reverse jeu de taquin
we fix a row $i$. Reverse jeu de taquin depends on this row.
The output of reverse jeu de taquin in $U$ with entry $e$ and with 
respect to $D$ is again another shifted tabloid $T$ of shape $\lambda$. 
Let $c_{U}(e)=(i',j')$. Set $T=U$ and stop 
if either $(i',j')=(i,i)$ or $(i',j') \in D$. Otherwise let $\rho \in \{ (i'-1,j'), (i', j'-1) \}$ be such that
$U_\rho = \max ( U_{i'-1,j'} , U_{i',j'-1})$ if $i' \notin \{i,j'\}$, $\rho=(i',j'-1)$ if $i'=i$ and 
$\rho=(i'-1,j')$ if $i'=j'$ and let $U'$ denote 
the tabloid we obtain by exchanging the entries $e$ and $U_\rho$ in $U$. 
Next perform jeu de taquin in $U'$ with entry $e$ and with respect to $D$ 
in order to obtain the final tabloid $T$. (For an example read Figure~\ref{jt} from 
right to left. The left shifted tabloid can be obtained from the right by performing
reverse jeu de taquin with $64$ and with respect to $(1,4)$.)

If $D$ is the empty set we omit `with respect to $D$'. In the following the 
output tabloid $T$ is denoted by $\rjt_D(U,e)$ and the cell of $e$ in 
$T$ is denoted by $\crjt_D(U,e)$. If $D$ is the $k$-th row we write 
$\rjt_k(U,e)$ and $\crjt_k(U,e)$, respectively.



\medskip

Let $T$ be a shifted tabloid of shape $\lambda$ and $1 \le e \le n$.  
The {\it forward path of $e$ in $T$ with respect to $D$} is the set of cells $e$ comes 
across when 
performing jeu de taquin in $T$ with $e$ and with respect to $D$. 
For example the forward path of $64$ with respect to the cells on the 
main diagonal in the left tabloid of Figure~\ref{jt} is 
$\{(1,4),(2,4),(3,4),(3,5),(4,5),(5,5)\}$.
Similarly the {\it backward path 
in $T$ of $e$ with respect to $D$} is the set of cells $e$ comes across when performing
reverse jeu de taquin in $T$ with $e$ and with respect to $D$. 
Clearly the backward path of $64$ with respect to $(1,4)$ in the right shifted
tabloid in Figure~\ref{jt} coincides with the forward path in the left shifted 
tabloid.

\medskip

We need one more definition before we are able to relate jeu de taquin and 
reverse jeu de taquin.
Let $T$ be a shifted tabloid of shape $\lambda$ and $\rho$ a cell in the shifted Ferrers diagram of 
shape $\lambda$.
The shifted tabloid $T$ is said to be {\it ordered up to cell $\rho$} if rows and columns are 
increasing in the subtabloid consisting of the cells that are smaller or equal to $\rho$ with 
respect to the total order. (Note that in Figure~\ref{jt} the left tabloid is ordered up to 
$(1,5)$.)

\medskip

Observe that the routines {\it jeu de taquin} and {\it reverse jeu de taquin}
are inverse to each other in the following sense: 
Let $T$ be a shifted tabloid and $e$ an entry in $T$ with 
$c_T(e)=(i,j)$. Suppose 
$T$ is ordered up to the predecessor of $(i,j)$ in the total order. Perform jeu de taquin in 
$T$ with $e$ and with respect to $D$ and obtain the tabloid $U$.
If we perform reverse jeu de taquin in $U$ with $e$ and with respect to $(i,j)$ ($i$ being the fixed row), 
we reobtain $T$.
In symbols
$$
T=\rjt_{(i,j)} (\jt_D(T,e),e).
$$
Furthermore observe that $U$ is ordered up to $(i,j)$ if $D=\emptyset$.

Conversely: Let $e$ be an entry weakly below the $i$-th row of a shifted tabloid $T$ 
that is ordered up to $(i,j)$ and 
assume that the backward path of $e$ in $T$ contains $(i,j)$  
($i$ being the fixed row). Perform 
reverse jeu de taquin in $T$ with $e$ and respect to $(i,j)$ 
and denote the output tabloid by $U$. If we perform 
jeu de taquin in $U$ with $e$ and with respect to $(i',j')=c_T(e)$, we reobtain $T$.
In symbols
$$
T=\jt_{(i',j')} ( \rjt_{(i,j)} (T,e),e).
$$
Furthermore observe that $U$ is ordered up to the predecessor of $(i,j)$.
Therefore: {\it In order to be able to reconstruct $T$ from $\ssplit(T)$ we have to store
the endcells of the forward paths we obtain in the application of modified jeu
de taquin to $T$ in the shifted hook tabloid. Since we want to obtain a bijection this storage 
has to be organized most efficiently. The fact that the endcells are not independent from each
other (see Lemma~\ref{1}) makes this storage non-trivial.} 

\medskip

\begin{figure}
\setlength{\unitlength}{1cm}
\begin{center}
\scalebox{0.40}{\includegraphics{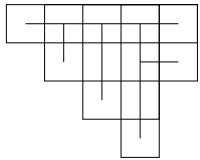}}
\end{center}
\caption{The backward paths of the entries in the shifted tabloid in Figure~\ref{standard}.b. Thus
$1 <_T 3 <_T 2 <_T 6 <_T 5 <_T 4 <_T 12 <_T 11 <_T 9 <_T 10 <_T 7 <_T 8$.}
\label{bpo}   
\end{figure}

{\bf Backward paths order.}
We define the {\it backward paths order} on the entries of a shifted tabloid $T$. 
Let $e_1$, $e_2$ 
be two entries in $T$ and $P_1$,
$P_2$ their backward paths ($1$ being the fixed row) in $T$. 
Furthermore let $(i,j)$ denote the smallest cell in 
$P_1 \cap P_2$ with respect to the total order. If either 
$(i+1,j) \in P_1$ or $(i,j+1) \in P_2$ we define 
$e_1 <_T e_2$ and say that $e_2$ is greater than $e_1$ with respect to 
the backward paths order in $T$. For an example see Figure~\ref{bpo}.

\medskip

\medskip

We introduce a manner-of-speaking: Let $T$ 
be a shifted tabloid and $H$ an accompanying partial shifted hook tabloid, 
a partial shifted hook tabloid being
a filling of some cells of the shifted Ferrer diagram with  entries satisfying 
the requirement for the entries in a shifted hook tabloid. Let 
$\rho$ be an entry of a cell $\sigma$ in the fixed $i$-th row of $H$,
i.e. $H_{\sigma}=\rho$. We say that 
$T_\rho$ is a {\it horizontal candidate} 
in $T$ with respect to $H$, if $\rho$ is in the same column as 
$\sigma$ in the shifted Ferrers diagram. For example 
consider the shifted hook tabloid in Figure~\ref{shbsp} and an 
arbitrary shifted tabloid $T$  of shape $(5,4,2,1)$. If $i=1$ then  $T_{1,2}$, $T_{4,4}$
and $T_{2,5}$ are the horizontal candidates.
If $\rho$ is neither in the same row nor in the same column as $\sigma$, 
$T_\rho$ is a {\it vertical candidate} with respect to $H$.
In our example in Figure~\ref{shbsp}  $T_{2,2}$ is the only vertical candidate if $i=1$.

\medskip

The rest of the paper is organized as follows.  
In Section~\ref{alg} we describe the 
Algorithm SPLIT that converts a shifted tabloid of shape $\lambda$ into a pair 
of a shifted standard tableau of shape $\lambda$ and a shifted hook tabloid of shape 
$\lambda$. In Section~\ref{examples} we give some examples of the application of SPLIT. 
In Section~\ref{inv} we describe the Algorithm MERGE that 
`merges' a pair of a shifted standard tableau and a shifted hook tabloid to 
a shifted tabloid. In Section~\ref{proof} we prove that the Algorithm MERGE 
is the inverse of the Algorithm SPLIT.

\section{The Algorithm SPLIT}
\label{alg}

In this section we describe the Algorithm SPLIT
that transforms a shifted tabloid into a pair of a shifted standard tableau 
and a shifted hook tabloid. 

\medskip

The construction of the shifted hook tabloid is more involved compared to the construction of the 
shifted standard tabloid, which we have already described. It
depends on the endcells of the forward paths we obtain in the course of performing jeu de taquin 
in the shifted tabloid and 
on the intersection of two such paths near the main diagonal. 
We have to modify the order of some 
steps in the construction of the shifted standard tableau 
so that the shifted hook tabloid can be  built up simultaneously. But since these 
steps in question commute the  modified algorithm for building the shifted standard tableau is equivalent to the 
origin algorithm.

\medskip

Before we are in the position to describe the algorithm, we 
introduce two routines on a partial shifted hook tabloid $H$ of shape $\lambda$.

\smallskip

\begin{figure}
\setlength{\unitlength}{1cm}
\begin{center}
\scalebox{0.35}{\includegraphics{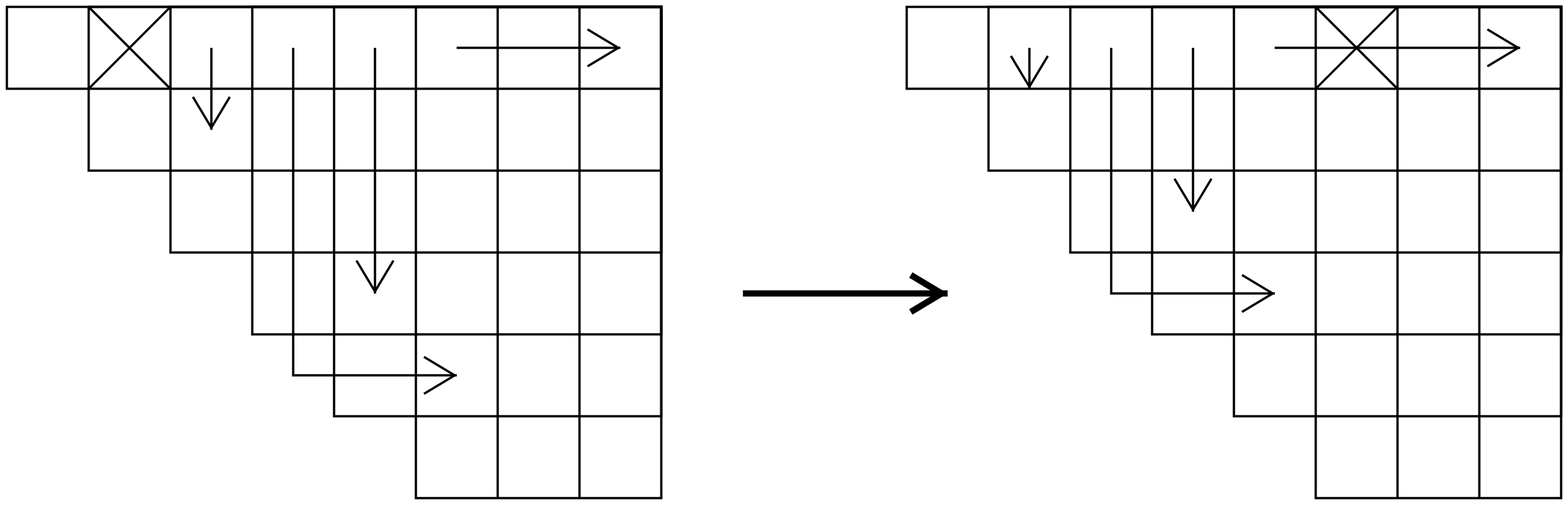}}
\end{center}
\caption{A shift from $(1,2)$ to $(1,6)$: The ends of the arrows indicate the entries of the cells at the origins of the arrows.
The x denotes an empty cell.}
\label{shift}   
\end{figure}

{\bf A shift from cell $(i,j)$ to cell $(i,j')$ in a partial shifted hook tabloid.} 
Let $(i,j)$, $(i,j')$ 
be two cells in $H$, $j \le j'$. We define the term {\it shift from 
$(i,j)$ to $(i,j')$ in $H$}. The output of this operation is another
partial shifted hook tabloid $H'$ which coincides with $H$ except for the 
cells $(i,k)$, $j \le k \le j'$. Let $j \le k < j'$. If 
$H_{i,k+1}=(i',k+1)$, $i < i'$,  or $H_{i,k+1}=(k+2,j')$ 
then set $H'_{i,k}=(i'-1,k)$  or $H'_{i,k}=(k+1,j'-1)$, respectively. 
If $H_{i,k+1}=(i,l)$ then
set $H'_{i,k}=(i,l)$. The cell $(i,j')$ in $H'$ remains 
empty. We denote $H'$ by $\shift(H,i,j,j')$. 
For an example see Figure~\ref{shift}.

\medskip

Whenever we perform jeu de taquin with an entry $e$ and with 
respect to a set $D$ in the shifted tabloid $T$ in the algorithms below, 
we store the end of the forward path in the accompanying partial shifted hook tabloid $H$. If we are 
in the course of performing jeu de taquin with the entries in the $i$-th row and 
$(i',j')$ is the end of the forward path, this is either done by setting $H_{i,j'}=(i',j')$ 
or by setting $H_{i,i'-1}=(i',j')$. 
The latter possibility is used if $e$ was previously a vertical candidate.
By reverse jeu de taquin the knowledge 
of the end of the forward path is enough to undo the performance of jeu de taquin.
However, when $H_{i,j'}$, respectively  $H_{i,i'-1}$, is occupied we have to perform a shift from the 
cell in $H$ that has previously pointed to $e$ to $(i,j')$, respectively  $(i,i'-1)$,
in order to empty either $(i,j')$ or $(i,i'-1)$. 
More accurate: If $e$ is either a horizontal candidate or no candidate in $(T,H)$ we denote by $\js_D(T,H,e)$  
the pair of a shifted tabloid $T'$ and a shifted hook tabloid $H'$, which is 
obtained in the following way: $T'=\jt_D(T,e)$ , $H'=\shift(H,i,q,j')$ and $H'_{i,j'}=(i',j')$, 
where 
$q$ is the column of $e$ in $T$ and $c_{T'}(e)=(i',j')$. If $e$ is a vertical candidate we denote by  $\js_D(T,H,e)$ 
the pair of a shifted tabloid 
$T'$ and a shifted hook tabloid $H'$, which is obtained in the following way: $T'=\jt_D(T,e)$, $H'=\shift(H,i,p-1,i'-1)$
and $H'_{i,i'-1}=(i',j')$, where $p$ is the row of $e$ in $T$ and $c_{T'}(e)=(i',j')$. 
If $D$ is empty we write $\js(T,H,e)$ and if $D$ is the $k$-th row we write $\js_k(T,H,e)$.

\medskip

\begin{figure}
\setlength{\unitlength}{1cm}
\begin{center}
\scalebox{0.35}{\includegraphics{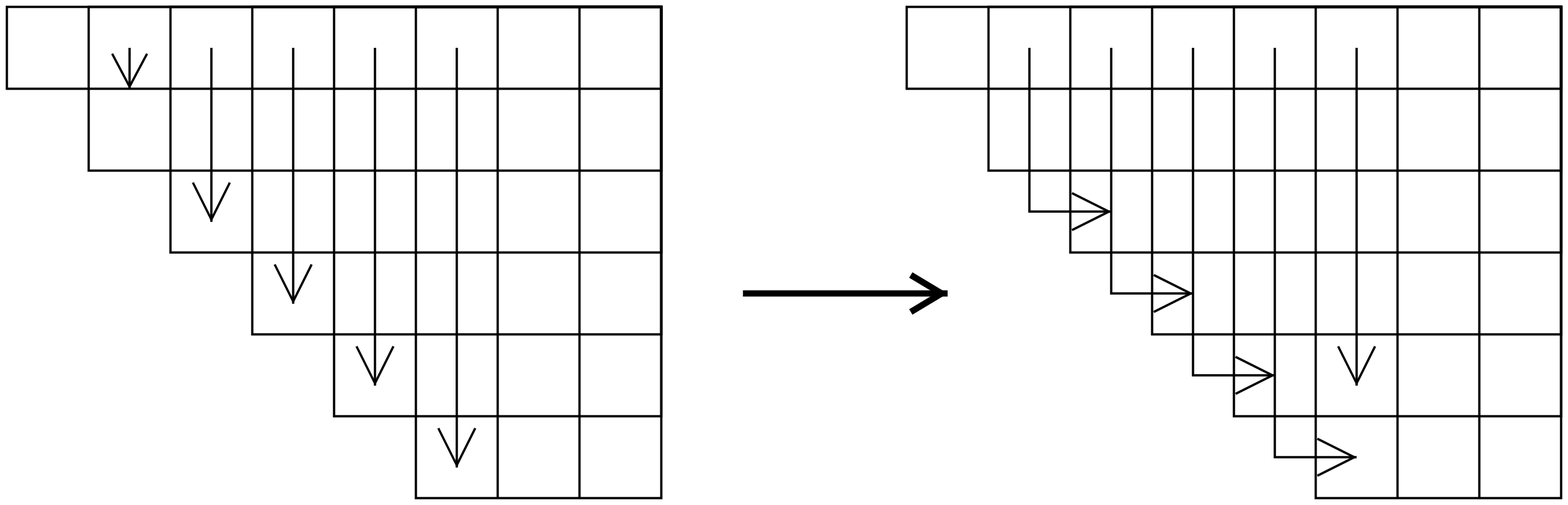}}
\end{center}
\caption{A transfer from $(1,2)$ to $(1,6)$: The ends of the arrows indicate the entries of the cells at the origins of the arrows.}
\label{trans} 
\end{figure}

{\bf A transfer from cell $(i,j)$ to cell $(i,k)$ in a partial
shifted hook tabloid.} Let $j \le k \le r$ and $H_{i,j}=(i',j')$ with 
either $i'=i$ or $j'=j$. We define the term {\it 
transfer from cell $(i,j)$ to cell $(i,k)$} in $H$. The output of this operation is 
again another partial shifted hook tabloid $H'$ which coincides
with $H$ except for the cells $(i,l)$, $j \le l \le k$. If $i'=i$ and $k \le j'$ let 
$H'_{i,k}=(i',j')$,  otherwise let $H'_{i,k}=(i'+k-j',k)$. For $j \le l < k$ let 
$H'_{i,l}=(l+1,l+1)$. We denote $H'$ by $\trans(H,i,j,k)$.
For an example see Figure~\ref{trans}.
Whenever we apply a transfer from $(i,j)$ to $(i,k)$ in $H$ we  
have $H_{i,l}=(l,l)$ for $j < l \le k$ as in the example.
This operation we need in the algorithms for 
converting horizontal candidates into vertical candidates.

\medskip

\medskip

\medskip

The Algorithm SPLIT is divided into $r$ steps, where in the $i$-th step 
we perform jeu de taquin with the entries in the $(r-i+1)$-st row. Within a row $i$, SPLIT is divided into 
3 steps, SPLIT~1, SPLIT~2 and SPLIT~3. 
Assume we just start performing jeu de taquin with the 
entries in the $i$-th row. Let $T$ denote the shifted tabloid we have constructed so far 
and $H$ the partial shifted hook tabloid. At this point $T$ is ordered up to 
$(i+1,i+1)$ (as it is in the origin algorithm for constructing the shifted standard tableau 
in Section~\ref{jeu}), the first $i$ rows of $H$ are 
empty and the last $r-i$ rows of $H$ form a shifted hook tabloid. 
In SPLIT~1 we perform jeu de taquin with the entries in the $i$-th row from right to left
and with respect to the set 
$MD=\{ (1,1), (2,2), \dots, (r,r) \}$ of cells on the main diagonal.

\medskip

\fbox{ \parbox{14cm}{

{\bf SPLIT 1.} 
Repeat for $j=\lambda_i + i -1$ down to $j=i$: Set $(T,H)=\js_{MD} (T,H,T_{i,j})$. 
}}

\smallskip

\medskip

After SPLIT~1 for every entry $e$ whose forward path
terminates in SPLIT~1 in a cell $(k,k)$ on the main diagonal, we have 
$T_{k,k}=e$ and $H_{i,k}=(k,k)$ and therefore all unstable entries in $T$
 are horizontal candidates with respect to $H$  after the application 
of SPLIT~1. This follows from the fact that
whenever a forward path of an entry $T_{i,j}$ ends in a cell 
$(k,k)$ on the main diagonal, then the forward paths of the following 
entries $T_{i,g}$, $g < j$, end strictly left of the $k$-th column.
(See Lemma~\ref{6}.)

\medskip

\begin{figure}
\setlength{\unitlength}{1cm}
\begin{center}
\scalebox{0.35}{\includegraphics{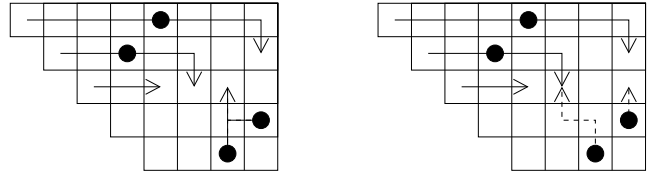}}
\end{center}
\caption{In the left tabloid we are in Case~1 of SPLIT~2 and in the right tabloid we are in Case~2 of SPLIT~2: 
In both diagrams $i=1$, $i'=h=3$ and the full lines indicate the forward paths in the course of 
constructing $U$. The circles indicate
the horizontal candidates and if they are strictly below the $3$-rd row the 
dashed lines indicate their backward paths 
with respect to the $3$-rd row.}   
\label{split2}   
\end{figure}

\fbox{ \parbox{14cm}{

{\bf SPLIT 2.} 
Choose $i'$, $i \le i' \le r$, maximal such that $H_{i,k}=(k,k)$ and 
$T_{k-1,k-1}$ is unstable (i.e. $T_{k-1,k-1} > T_{k-1,k}$) for $i <  k \le i'$.   
If $i'=i$ jump to Case~2.

\smallskip

Set $U=T$. Repeat for $g=i'$ down to $g=i$:  Set $U=\jt_{g+1}(U,T_{g,g})$.

\medskip

We distinguish between two cases. 
Let $h$ be minimal, $i \le h \le i'$, such that $T_{h,h}$
is not in the $(h+1)$-st row of $U$. We continue with Case 2  below if 
$h$ does not exist. We also continue with Case 2  if $h=i'$ and there exists a horizontal 
candidate with respect to $H$ 
strictly below the $i'$-th row of $U$ which is smaller than $T_{i'-1,i'-1}$ with 
respect to the backward paths order in $U$.
In all other cases we continue with Case~1. (See Figure~\ref{split2}.)

\medskip

Reject the tabloid $U$ we constructed so far in SPLIT~2.

\medskip

{\it Case 1.} 
Repeat for $g=i'$ down to $g=h+1$:
[ Let $k$ be such that either $(g,k)$ or $(g+1,k)$ is the endcell of the 
forward path of $T_{g,g}$ in the procedure for constructing $U$.
Set $T=\jt_{(g,k)}(T,T_{g,g})$ and $H_{i,g-1}=(g,k)$.   
{ \small (Note that the forward path of $T_{g,g}$ ends in $(g,k)$.)}  ]

Let $(h,k)=\cjt(T,T_{h,h})$ and set $T=\jt(T,T_{h,h})$.
If $h - k  \le i - i'$ let $H_{i,i'}=(i,i - h + k)$ otherwise let $H_{i,i'}=(i'+h-k,i')$.

Repeat for $g=h-1$ down to $g=i$: [ Let $H_{i,g}=\cjt_{g+1}(T,T_{g,g})$ and 
$T=\jt_{g+1}(T,T_{g,g})$. ]

\smallskip

{\it Case 2.}
Set $(T,H)=\js_{i'+1} (T,H,T_{i',i'})$. 

Repeat for $g=i'-1$ down to $g=i$: [ Let $H_{i,g}=\cjt_{g+1} (T,T_{g,g})$ and $T=\jt_{g+1}(T,T_{g,g})$. ]
}}

\medskip

\medskip

Observe the following (See also Figure~\ref{split2}.): Consider the backward paths with respect to the $i'$-th row 
of the horizontal candidates strictly below the $i'$-th row in the output tabloid of SPLIT~2. 
Then either one of these paths ends weakly to the left of the vertical candidate 
in the $i'$-th row or there exists no vertical candidate if  and only if we were in Case~2 of SPLIT~2. (If 
we were in Case~1 of SPLIT~2 and the backward path of 
a horizontal candidate strictly below the $i'$-th row would end weakly to the left of 
the vertical candidate in the $i'$-th row after the application of Case~1 SPLIT~2 
then $h \not= i'$ by the cases distinction in SPLIT~2.
Thus the  backward path of the same horizontal candidate includes the cell $(i'+1,i'+1)$ on the main diagonal
before the application of Case~1 in SPLIT~2 (by the argument in the proof of Lemma~\ref{1})
which implies $H_{i,i'+1}=(i'+1,i'+1)$ (by property AFTER\_SPLIT'~1 in the
proof of Claim~1 in Section~\ref{proof}) and this is 
not possible by the choice of $i'$ for our assumptions imply that $T_{i',i'}$ is unstable before the
application of Case~1 in SPLIT~2. Furthermore there always exists a vertical candidate 
after the performance of Case~1 of SPLIT~2.) 
In other words: We were in Case~2 of SPLIT~2 if and only if
either the smallest
horizontal candidate strictly below the row of the smallest vertical candidate is smaller than this 
smallest vertical candidate or there exists no vertical candidate  in the output pair.

\medskip

Besides the pair $(T,H)$ we need two sets of entries $C'$ and $C$ as an input for SPLIT~3. 
The set $C'$ contains the vertical candidates 
together with the smallest candidate with respect to the backward paths order
if it is unstable.
The set $C$ contains the horizontal 
candidates on the main diagonal of $T$ which are not in $C'$.
Observe that the entries in $C'$ are strictly above of the entries
in $C$, i.e. the maximal row of an entry in $C'$ is smaller than 
the minimal row of an entry in $C$. This will also be true in the course of performing SPLIT~3.
Similarly all vertical candidates are strictly above of the 
entries in $C$ before and while performing SPLIT~3.

\medskip

\fbox{ \parbox{14cm}{

{\bf SPLIT 3.}
Repeat the following until $C \cup C' = \emptyset$.

\smallskip

[If $C' = \emptyset$ choose 
$\overline{e} \in C$ such that the row of $\overline{e}$ is minimal, set $C' = C' \cup \{ \overline{e}\}$ 
and  $C = C \setminus \{ \overline{e} \}$. Choose $e \in C'$ with maximal row.

If $C \not= \emptyset$,  choose $e'=T_{h,h} \in C$ such that $h$ is minimal. 

Consider the forward path $P$ of $e$ in $T$. We distinguish between three cases.

\smallskip

{\it Case 1.} $P$ ends weakly left of $(h-1)$-st column  or $e'$ does not exist.

Set $(T,H)=\js(T,H,e)$ and  $C'= C' \setminus \{e\}$.

{\it Case 2.} 
$(h-1,h-1) \in P$ and  $(h-1,h) \in P$. 

Set $(T,H)=\js_{(h-1,h-1)}(T,H,e)$.

Let $h$, $h \le h'$, be maximal such that $T_{j-1,j-1}$ is unstable and 
$T_{j,j} \in C$ for $h < j \le h'$.

If $e$ is a horizontal candidate set $H=\trans(H,i,h-2,h'-1)$, otherwise $H=\trans(H,i,h-1,h')$. 

Furthermore let $C'= C' \cup \{T_{h,h}, T_{h+1,h+1}, \dots, T_{h',h'}\}$ and 
$C = C \setminus \{ T_{h,h}, T_{h+1,h+1}, \dots, T_{h',h'} \}$. 

\smallskip

{\it Case 3.} $P$ ends weakly right of column $h$, but  does not contain $(h-1,h-1)$. 
Set $C'=C' \cup \{e'\}$ and $C = C \setminus \{e'\}$. ]
}}

\medskip

\begin{figure}
\setlength{\unitlength}{1cm}
\begin{center}
\scalebox{0.40}{\includegraphics{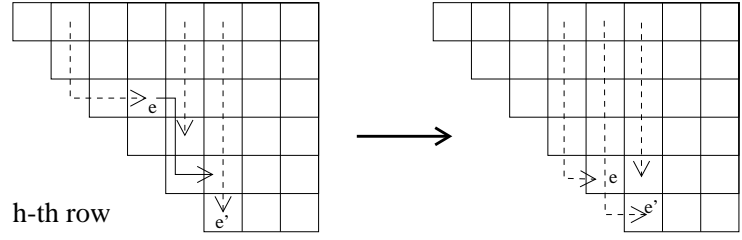}}
\end{center}
\caption{Situation before (left) and after (right) Case~2 in SPLIT~3. The dashed lines 
indicate relevant candidates with respect to the accompanying shifted hook
tabloid (see also Figure~\ref{shift}).
The full line indicates the forward path $P$ of $e$.}
\label{split3}   
\end{figure}

In Figure~\ref{split3} an example of Case~2 in SPLIT~3 is displayed.

\medskip

Let $T$ be a shifted tabloid ordered up to $(i+1,i+1)$ 
and $H$ a partial shifted hook tabloid. Let $(T',H')$ denote the 
pair we obtain after the application of SPLIT~1, SPLIT~2 and SPLIT~3 
to the $i$-th row of $(T,H)$. We denote 
$T'=\ssplit(T,H,i)=\ssplit(T,i)$ and $H'=\hsplit(T,H,i)$.  We are 
finally in the position to formulate the Algorithm SPLIT. The input is 
a shifted tabloid $T$ of shape $\lambda$ and an empty shifted tabloid $H'$ of shape $\lambda$.

\medskip

\fbox{ \parbox{14cm}{
{\bf SPLIT.} $T'=T$. Repeat for $i=r$ down to $i=1$: 
[ $H'=\hsplit(T',H',i)$ and $T'=\ssplit(T',H',i)$.]
}}

\medskip

The output $T'$ and $H'$ of the algorithm is denoted by $\ssplit(T,H)=\ssplit(T)$ and 
$\hsplit(T,H)$.

\medskip

\section{Examples}
\label{examples}

In this section we give five examples for the application of SPLIT.
All examples are of shape $(11,10,9,8,7,6,5,4,3,2)$ and the input shifted tabloid is 
ordered up to $(2,2)$. 
We perform jeu de taquin with the entries in the 
first row of the shifted tabloid and simultaneously built up the first row of the 
shifted hook tabloid. Our examples cover all possible cases in SPLIT~2:
In Example~\ref{ex1} we are in Case~2 of SPLIT~2 for $h$ does not exist, 
in Example~\ref{best} we are in Case~2 of SPLIT~2 for $i=i'$, 
in Example~\ref{ex3} we are in Case~1 of SPLIT~2 for $h=1$, 
in Example~\ref{ex4} we are in Case~1 of SPLIT~2 with $h=i'$ and 
in Example~\ref{ex5} we are in Case~2 of SPLIT~2 with $h=i'$.

\medskip

\begin{ex} 
\label{ex1}
We consider the shifted tabloid 
$$  T=
\begin{array}{cccccccccccc}
63 & 64 & 61 &  41 &  34 &  20 & 62 & 65   & 56 &  15 &  54 \\
    &   1 &   2 &   3 &    4 &   7 &   8 &  13  &  16 &  18 &  24 \\
    &     &   5 &   6 &   9 &  10 &  17 &  19  &  30 &  31 &  45 \\
    &     &     &  11 &  12 &  21 &  23 &  26  &  32 &  37 &  47 \\
    &     &     &     &  14 &  22 &  25 &  27  &  36 &  39 &  51 \\
    &     &     &     &     &  28 &  29 &  35  &  40 &  44 &  52 \\
    &     &     &     &     &     &  33 &  38  &  43 &  48 &  53 \\
    &     &     &     &     &     &     &  42  &  46 &  49 &  55 \\
    &     &     &     &     &     &     &      &  50 & 57 & 58 \\
    &     &     &     &     &     &     &      &     & 59 & 60
\end{array}.  
$$
First we have to apply SPLIT~1 to the first row. If we perform jeu de taquin 
with $54$ and with respect to the cells on the main diagonal, the entry 
ends in cell $(7,11)$. 
$$
T=
\begin{array}{cccccccccccc}
63 & 64 & 61 &  41 &  34 &  20 & 62 & 65   & 56 &  15 &  24 \\
    &   1 &   2 &   3 &    4 &   7 &   8 &  13  &  16 &  18 &  45 \\
    &     &   5 &   6 &   9 &  10 &  17 &  19  &  30 &  31 &  47 \\
    &     &     &  11 &  12 &  21 &  23 &  26  &  32 &  37 &  51 \\
    &     &     &     &  14 &  22 &  25 &  27  &  36 &  39 &  52 \\
    &     &     &     &     &  28 &  29 &  35  &  40 &  44 &  53 \\
    &     &     &     &     &     &  33 &  38  &  43 &  48 &  54 \\
    &     &     &     &     &     &     &  42  &  46 &  49 &  55 \\
    &     &     &     &     &     &     &      &  50 & 57 & 58 \\
    &     &     &     &     &     &     &      &     & 59 & 60
\end{array}  
$$
Thus we set $H_{1,11}=(7,11)$ and therefore the first row of the 
shifted hook tabloid is
$$
H_1= \left(
(-,-) , (-,-) , (-,-) , (-,-) , (-,-) , (-,-) , (-,-) , (-,-) , (-,-) , (-,-) , (7,11) 
 \right) .
$$
The entry $15$ is stable, thus the shifted tabloid does not change in the next step and 
we set $H_{1,10}=(1,10)$, the cell of $15$.
$$ 
H_1= \left(
(-,-) , (-,-) , (-,-) , (-,-) , (-,-) , (-,-) , (-,-) , (-,-) , (-,-) , (1,10) , (7,11) 
\right)
$$
Next we perform jeu de taquin with $56$. 
$$
T=
\begin{array}{cccccccccccc}
63 & 64 & 61 &  41 &  34 &  20 & 62 & 65   &  15 &  18 &  24 \\
    &   1 &   2 &   3 &    4 &   7 &   8 &  13  &  16 &  31 &  45 \\
    &     &   5 &   6 &   9 &  10 &  17 &  19  &  30 &  37 &  47 \\
    &     &     &  11 &  12 &  21 &  23 &  26  &  32 &  39 &  51 \\
    &     &     &     &  14 &  22 &  25 &  27  &  36 &  44 &  52 \\
    &     &     &     &     &  28 &  29 &  35  &  40 &  48 &  53 \\
    &     &     &     &     &     &  33 &  38  &  43 &  49 &  54 \\
    &     &     &     &     &     &     &  42  &  46 &  55 & 56 \\
    &     &     &     &     &     &     &      &  50 & 57 & 58 \\
    &     &     &     &     &     &     &      &     & 59 & 60
\end{array}  
$$
The forward path ends in cell $(8,11)$ 
and therefore we want to set $H_{1,11}=(8,11)$. However, $H_{1,11}$ is 
occupied, thus we perform a shift from $(1,9)$ to $(1,11)$.
$$ H_1 = \left(
(-,-) , (-,-) , (-,-) , (-,-) , (-,-) , (-,-) , (-,-) , (-,-) , (1,10) , (6,10) , (8,11) 
\right)
$$
If we perform jeu de taquin with $65$ and with respect to the cells on the 
main diagonal, the entry gets stuck in $(9,9)$ and remains unstable there.
$$ 
T=
\begin{array}{cccccccccccc}
63 & 64 & 61 &  41 &  34 &  20 & 62 &  13  &  15 &  18 &  24 \\
    &   1 &   2 &   3 &    4 &   7 &   8 &  16  &  30 &  31 &  45 \\
    &     &   5 &   6 &   9 &  10 &  17 &  19  &  32 &  37 &  47 \\
    &     &     &  11 &  12 &  21 &  23 &  26  &  36 &  39 &  51 \\
    &     &     &     &  14 &  22 &  25 &  27  &  40 &  44 &  52 \\
    &     &     &     &     &  28 &  29 &  35  &  43 &  48 &  53 \\
    &     &     &     &     &     &  33 &  38  &  46 &  49 &  54 \\
    &     &     &     &     &     &     &  42  &  50 &  55 & 56 \\
    &     &     &     &     &     &     &      & 65  & 57 & 58 \\
    &     &     &     &     &     &     &      &     & 59 & 60
\end{array}  
$$
Again we want to set $H_{1,9} = (9,9)$, but $H_{1,9}$ is occupied 
and therefore we have to perform a shift.
$$
H_1 = \left(
(-,-) , (-,-) , (-,-) , (-,-) , (-,-) , (-,-) , (-,-) , (1,10) , (9,9) , (6,10) , (8,11) 
\right)
$$
Next we perform jeu de taquin 
with $62$ and this entry again gets stuck in a cell on the main 
diagonal. Then we perform jeu de taquin with $20$, $34$, $41$, $61$, $64$ and finally with 
$63$. These entries either end in a stable position or on the main diagonal.
Simultaneously we built up the first row of the shifted hook tabloid.
There we put down the endcells of the paths in the appropriate columns and whenever we want to use 
a cell of the shifted hook tabloid which is already occupied we perform the appropriate  shift.
This procedure results in 
$$ T = 
\begin{array}{cccccccccccc}
63 &   1 &   2 &   3 &    4 &   7 &   8 &  13  &  15 &  18 &  24 \\
    & 64 &   5 &   6 &   9 &  10 &  16 &  19  &  30 &  31 &  45 \\
    &     & 61 &  11 &  12 &  17 &  20 &  26  &  32 &  37 &  47 \\
    &     &     &  41 &  14 &  21 &  23 &  27  &  36 &  39 &  51 \\
    &     &     &     &  34 &  22 &  25 &  35  &  40 &  44 &  52 \\
    &     &     &     &     &  28 &  29 &  38  &  43 &  48 &  53 \\
    &     &     &     &     &     &  33 &  42  &  46 &  49 &  54 \\
    &     &     &     &     &     &     & 62  &  50 &  55 & 56 \\
    &     &     &     &     &     &     &      & 65  & 57 & 58 \\
    &     &     &     &     &     &     &      &     & 59 & 60
\end{array} 
$$
and 
$$ 
H_1 = \left( 
(1,1) , (2,2) , (3,3) , (4,4) , (5,5) , (1,10) , (3,7) , (8,8) , (9,9) , (6,10) , (8,11) 
\right).
$$

We apply SPLIT~2 to this pair. Observe that $i'=5$ since $H_{1,j}=(j,j)$ for $1 \le j \le 5$
and $63$, $64$, $61$ and $41$ are unstable. We perform jeu de taquin with 
$34$, $41$, $61$, $64$ and $63$ in that order and with repsect to the 
relative next row. 
$$ \small
U=
\begin{array}{cccccccccccc}
  1 &   2 &   3 &    4 &   7 &   8 &  13 &  15  &  18 &  24 &  45 \\
    &   5 &   6 &   9 &  10 &  16 &  19 &  30  &  31 &  37 & 63 \\
    &     &  11 &  12 &  17 &  20 &  26 &  32  &  36 & 64 &  47 \\
    &     &     &  14 &  21 &  23 &  27 &  35  & 61 &  39 &  51 \\
    &     &     &     &  22 &  25 &  29 &  41  &  40 &  44 &  52 \\
    &     &     &     &     &  28 &  34 &  38  &  43 &  48 &  53 \\
    &     &     &     &     &     &  33 &  42  &  46 &  49 &  54 \\
    &     &     &     &     &     &     & 62  &  50 &  55 & 56 \\
    &     &     &     &     &     &     &      & 65  & 57 & 58 \\
    &     &     &     &     &     &     &      &     & 59 & 60
\end{array}
$$
We observe that every entry changes row and therefore $h$ does not 
exist. Consequently we are in Case~2 of SPLIT~2. We set $T=U$. Since 
$34$ is in cell $(6,7)$, we set $H_{1,7}=(6,7)$ after 
performing a shift, $41$ is in cell $(5,8)$ and therefore
$H_{1,4}=(5,8)$, $61$ is in cell $(4,9)$ and therefore $H_{1,3}=(4,9)$, 
$64$ is in cell $(3,10)$ and therefore $H_{1,2}=(3,10)$ and $63$ is in 
cell $(2,11)$ and therefore $H_{1,1}=(2,11)$. Thus 
$$
H_1= \left(
(2,11) , (3,10) , (4,9) , (5,8) , (1,10) , (2,6) , (6,7) , (8,8) , (9,9) , (6,10) , (8,11) 
\right).
$$

Finally we apply SPLIT~3. Observe that $C'=\{34,41,61,64,63\}$ and $C=\{62,65\}$.
In the first step of SPLIT~3 $e=34$ and $e'=62$. The forward path of $e$ ends 
left of the column of $e'$ and therefore we are in Case~1. We obtain 
$$ 
T=
\begin{array}{cccccccccccc}
  1 &   2 &   3 &    4 &   7 &   8 &  13 &  15  &  18 &  24 &  45 \\
    &   5 &   6 &   9 &  10 &  16 &  19 &  30  &  31 &  37 & 63 \\
    &     &  11 &  12 &  17 &  20 &  26 &  32  &  36 & 64 &  47 \\
    &     &     &  14 &  21 &  23 &  27 &  35  & 61 &  39 &  51 \\
    &     &     &     &  22 &  25 &  29 &  41  &  40 &  44 &  52 \\
    &     &     &     &     &  28 &  33 &  38  &  43 &  48 &  53 \\
    &     &     &     &     &     &  34 &  42  &  46 &  49 &  54 \\
    &     &     &     &     &     &     & 62  &  50 &  55 & 56 \\
    &     &     &     &     &     &     &      & 65  & 57 & 58 \\
    &     &     &     &     &     &     &      &     & 59 & 60
\end{array}  
$$
and 
$$ 
H_1= \left(
(2,11) , (3,10) , (4,9) , (5,8) , (1,10) , (2,6) , (7,7) , (8,8) , (9,9) , (6,10) , (8,11) 
\right).
$$
We delete $34$ from $C'$: $C'=\{41,61,64,63\}$ and $C=\{62,65\}$.

In the next step of SPLIT~3 we have $e=41$ and $e'=62$. Since the forward path of $41$ 
ends weakly right of the column of $62$ but does not contain $(7,7)$, we 
are in Case~3. Here we do not change $T$ or $H$, we only delete 
$62$ from $C$ and put it into $C'$: $C'=\{62,41,61,64,63\}$ and $C=\{65\}$.

Now $e=62$ and $e'=65$ and we are in Case~2 for the forward path of $62$ contains 
$(8,8)$ and $(8,9)$. Therefore we want to change $62$ into a vertical
candidate by setting  $H_{1,7}=(8,8)$, but since 
$H_{1,7}$ is occupied we perform a transfer from $(1,7)$ to $(1,8)$. 
$$ H_1 = \left( 
(2,11) , (3,10) , (4,9) , (5,8) , (1,10) , (2,6) , (8,8) , (8,8) , (9,9) , (6,10) , (8,11) 
\right)
$$
The shifted tabloid $T$ does not change. 
Furthermore $C'=\{65,62,41,61,64,63\}$, $C=\emptyset$ and we are always in Case~1
for the rest of the application of SPLIT~3.

Next $e=65$ and its forward path ends in $(10,11)$. Therefore we want to set 
$H_{1,11}=(10,11)$ and since $H_{1,11}$ is occupied we perform a shift 
from the previous pointer $(1,9)$ of $65$ to $(1,11)$.
We obtain 
$$ T= 
\begin{array}{cccccccccccc}
  1 &   2 &   3 &    4 &   7 &   8 &  13 &  15  &  18 &  24 &  45 \\
    &   5 &   6 &   9 &  10 &  16 &  19 &  30  &  31 &  37 & 63 \\
    &     &  11 &  12 &  17 &  20 &  26 &  32  &  36 & 64 &  47 \\
    &     &     &  14 &  21 &  23 &  27 &  35  & 61 &  39 &  51 \\
    &     &     &     &  22 &  25 &  29 &  41  &  40 &  44 &  52 \\
    &     &     &     &     &  28 &  33 &  38  &  43 &  48 &  53 \\
    &     &     &     &     &     &  34 &  42  &  46 &  49 &  54 \\
    &     &     &     &     &     &     & 62  &  50 &  55 & 56 \\
    &     &     &     &     &     &     &      & 57 & 58 & 60 \\
    &     &     &     &     &     &     &      &     & 59 & 65 
\end{array}
$$ 
$$
H_1 = \left(
(2,11) , (3,10) , (4,9) , (5,8) , (1,10) , (2,6) , (8,8) , (8,8) , (5,9) , (7,10) , (10,11) 
 \right)
$$
and $C'=\{62,41,61,64,63\}$.

The forward path of $e=62$ ends in $(9,11)$. Therefore we want to set 
$H_{1,8}=(9,11)$ and since $H_{1,8}$ is occupied we perform a shift from 
$(1,7)$ to $(1,8)$.  This yields 
$$ T= 
\begin{array}{cccccccccccc}
  1 &   2 &   3 &    4 &   7 &   8 &  13 &  15  &  18 &  24 &  45 \\
    &   5 &   6 &   9 &  10 &  16 &  19 &  30  &  31 &  37 & 63 \\
    &     &  11 &  12 &  17 &  20 &  26 &  32  &  36 & 64 &  47 \\
    &     &     &  14 &  21 &  23 &  27 &  35  & 61 &  39 &  51 \\
    &     &     &     &  22 &  25 &  29 &  41  &  40 &  44 &  52 \\
    &     &     &     &     &  28 &  33 &  38  &  43 &  48 &  53 \\
    &     &     &     &     &     &  34 &  42  &  46 &  49 &  54 \\
    &     &     &     &     &     &     &  50  &  55 & 56 & 60 \\
    &     &     &     &     &     &     &      & 57 & 58 & 62 \\
    &     &     &     &     &     &     &      &     & 59 & 65 
\end{array}.
$$
Moreover
$$ 
H_1 = \left(
(2,11) , (3,10) , (4,9) , (5,8) , (1,10) , (2,6) , (7,7) , (9,11) , (5,9) , (7,10) , (10,11) 
\right)
$$
and $C'=\{41,61,64,63\}$.

In the next step $e=41$ etc. $\dots$

Finally we obtain the following shifted standard tableau
$$
T=
\begin{array}{cccccccccccc}
  1 &   2 &   3 &    4 &   7 &   8 &  13 &  15  &  18 &  24 &  45 \\
    &   5 &   6 &   9 &  10 &  16 &  19 &  30  &  31 &  37 &  47 \\
    &     &  11 &  12 &  17 &  20 &  26 &  32  &  36 &  44 &  51 \\
    &     &     &  14 &  21 &  23 &  27 &  35  &  39 &  48 &  52 \\
    &     &     &     &  22 &  25 &  29 &  38  &  40 &  49 & 60 \\
    &     &     &     &     &  28 &  33 &  41  &  43 &  53 & 61 \\
    &     &     &     &     &     &  34 &  42  &  46 &  54 & 62 \\
    &     &     &     &     &     &     &  50  &  55 & 56 & 63 \\
    &     &     &     &     &     &     &      & 57 & 58 & 64 \\
    &     &     &     &     &     &     &      &     & 59 & 65 
\end{array}  
$$
and the following first row of $H$
$$ 
H_1 = \left(
(1,10) , (3,5) , (1,5) , (4,4) , (6,9) , (7,9) , (8,11) , (9,11) , (5,9) , (7,10) , (10,11) 
\right).
$$

\end{ex}

\medskip

\begin{ex}
\label{best}
We consider
$$ T=
\begin{array}{cccccccccccc}
 64 &  35 &   8 &  12 &  36 &  49 &  51 &   1  &  34 &  63 &  54 \\
    &   2 &   3 &   5 &   6 &  10 &  13 &  14  &  15 &  20 &  27 \\
    &     &   4 &   7 &   9 &  11 &  16 &  21  &  25 &  30 &  31 \\
    &     &     &  17 &  18 &  19 &  24 &  28  &  29 &  37 &  41 \\
    &     &     &     &  22 &  23 &  26 &  33  &  42 &  43 &  45 \\
    &     &     &     &     &  32 &  38 &  40  &  46 &  47 &  53 \\
    &     &     &     &     &     &  39 &  44  &  50 &  52 &  58 \\
    &     &     &     &     &     &     &  48  &  55 &  56 &  61 \\
    &     &     &     &     &     &     &      &  57 &  59 &  62 \\
    &     &     &     &     &     &     &      &     &  60 &  65
\end{array}.  
$$
After the application of SPLIT~1 we obtain the pair 
$$
T=
\begin{array}{cccccccccccc}
 64 &   1 &   3 &   5 &   6 &  10 &  13 &  14  &  15 &  20 &  27 \\
    &   2 &   4 &   7 &   9 &  11 &  16 &  21  &  25 &  30 &  31 \\
    &     &  35 &   8 &  12 &  19 &  24 &  28  &  29 &  37 &  41 \\
    &     &     &  17 &  18 &  23 &  26 &  33  &  34 &  43 &  45 \\
    &     &     &     &  22 &  32 &  38 &  40  &  42 &  47 &  53 \\
    &     &     &     &     &  36 &  39 &  44  &  46 &  52 &  54 \\
    &     &     &     &     &     &  49 &  48  &  50 &  56 &  58 \\
    &     &     &     &     &     &     &  51  &  55 &  59 &  61 \\
    &     &     &     &     &     &     &      &  57 &  60 &  62 \\
    &     &     &     &     &     &     &      &     &  63 &  65
\end{array},
$$
$$ 
H_1= \left(
(1,1), (1,8), (3,3), (3,4), (3,5) , (6,6), (7,7),  (8,8), (4,9) , (10,10) , (6,11)
\right).
$$

Since $i'=i$ we are in Case~2 of SPLIT~2. We move $64$ to the second row 
$$
T=
\begin{array}{cccccccccccc}
  1 &   2 &   3 &   5 &   6 &  10 &  13 &  14  &  15 &  20 &  27 \\
    &  64 &   4 &   7 &   9 &  11 &  16 &  21  &  25 &  30 &  31 \\
    &     &  35 &   8 &  12 &  19 &  24 &  28  &  29 &  37 &  41 \\
    &     &     &  17 &  18 &  23 &  26 &  33  &  34 &  43 &  45 \\
    &     &     &     &  22 &  32 &  38 &  40  &  42 &  47 &  53 \\
    &     &     &     &     &  36 &  39 &  44  &  46 &  52 &  54 \\
    &     &     &     &     &     &  49 &  48  &  50 &  56 &  58 \\
    &     &     &     &     &     &     &  51  &  55 &  59 &  61 \\
    &     &     &     &     &     &     &      &  57 &  60 &  62 \\
    &     &     &     &     &     &     &      &     &  63 &  65
\end{array}  
$$
and set $H_{1,2}=(2,2)$ after performing the appropriate shift
$$ 
H_1= \left(
(1,8), (2,2), (3,3), (3,4), (3,5) , (6,6), (7,7),  (8,8), (4,9) , (10,10) , (6,11)
\right).
$$

We set 
$C'=\{64\}$ and $C=\{35,36,49,51,63\}$. In the first step of SPLIT~3
we have $e=64$ and $e'=35$. Since the forward path of $64$ includes 
$(2,2)$ and $(2,3)$ we are in Case~2 of SPLIT~3. The shifted tabloid 
$T$ does not change, but we want to set $H_{1,1}=(2,2)$ and therefore 
perform a transfer from $(1,1)$ to $(1,2)$ in $H$. 
$$ 
H_1= \left(
(2,2), (1,8), (3,3), (3,4), (3,5) , (6,6), (7,7),  (8,8), (4,9) , (10,10) , (6,11)
\right)
$$ 
Furthermore $C'=\{35,64\}$ and $C=\{36,49,51,63\}$.

Next $e=35$ and $e'=36$. Again we are in Case~2.
$$
T=
\begin{array}{cccccccccccc}
  1 &   2 &   3 &   5 &   6 &  10 &  13 &  14  &  15 &  20 &  27 \\
    &  64 &   4 &   7 &   9 &  11 &  16 &  21  &  25 &  30 &  31 \\
    &     &   8 &  12 &  18 &  19 &  24 &  28  &  29 &  37 &  41 \\
    &     &     &  17 &  22 &  23 &  26 &  33  &  34 &  43 &  45 \\
    &     &     &     &  35 &  32 &  38 &  40  &  42 &  47 &  53 \\
    &     &     &     &     &  36 &  39 &  44  &  46 &  52 &  54 \\
    &     &     &     &     &     &  49 &  48  &  50 &  56 &  58 \\
    &     &     &     &     &     &     &  51  &  55 &  59 &  61 \\
    &     &     &     &     &     &     &      &  57 &  60 &  62 \\
    &     &     &     &     &     &     &      &     &  63 &  65
\end{array}  
$$
We set $H_{1,4}=(5,5)$ after performing an appropriate shift and a transfer from 
$(1,4)$ to $(1,5)$.
$$ 
H_1= \left(
(2,2), (1,8), (2,3), (5,5), (3,5) , (6,6), (7,7),  (8,8), (4,9) , (10,10) , (6,11)
\right)
$$ 
We update $C'=\{36,35,64\}$ and $C=\{49,51,63\}$.

Next $e=36$ and $e'=49$ and since $36$ is stable we are in Case~1. 
Thus $T$ and $H$ do not change, but $C'=\{35,64\}$ and $C=\{49,51,63\}$.

Therefore $e=35$ and $e'=49$ and we are in Case~1.
$$
T=
\begin{array}{cccccccccccc}
  1 &   2 &   3 &   5 &   6 &  10 &  13 &  14  &  15 &  20 &  27 \\
    &  64 &   4 &   7 &   9 &  11 &  16 &  21  &  25 &  30 &  31 \\
    &     &   8 &  12 &  18 &  19 &  24 &  28  &  29 &  37 &  41 \\
    &     &     &  17 &  22 &  23 &  26 &  33  &  34 &  43 &  45 \\
    &     &     &     &  32 &  35 &  38 &  40  &  42 &  47 &  53 \\
    &     &     &     &     &  36 &  39 &  44  &  46 &  52 &  54 \\
    &     &     &     &     &     &  49 &  48  &  50 &  56 &  58 \\
    &     &     &     &     &     &     &  51  &  55 &  59 &  61 \\
    &     &     &     &     &     &     &      &  57 &  60 &  62 \\
    &     &     &     &     &     &     &      &     &  63 &  65
\end{array}  
$$
$$ 
H_1= \left(
(2,2), (1,8), (2,3), (5,6), (3,5) , (6,6), (7,7),  (8,8), (4,9) , (10,10) , (6,11)
\right)
$$ 
Moreover 
$C'=\{64\}$ and $C=\{49,51,63\}$.

In the next step $e=64$ and $e'=49$. Because of the run of the forward path of 
$64$ we are in Case~3. Thus $T$ and $H$ do not change, only 
$C'=\{49,64\}$ and $C=\{51,63\}$.

Therefore $e=49$ and $e'=51$. We are in Case~2 and $T$ does not change.
We want to set $H_{1,6}=(7,7)$.
Since $H_{1,6}$ is occupied by a pointer which does not point to $49$, we 
perform a transfer from $(1,6)$ to $(1,7)$.
$$ 
H_1= \left(
(2,2), (1,8), (2,3), (5,6), (3,5) , (7,7), (7,7),  (8,8), (4,9) , (10,10) , (6,11)
\right)
$$ 
Moreover 
$C'=\{51,49,64\}$ and $C=\{63\}$.

Next $e=51$ and $e'=63$. We are in Case~1 for $51$ is stable and thus 
$C'=\{49,64\}$ and $C=\{63\}$. 

Therefore $e=49$ and $e'=63$. Again we are in Case~1.
$$
T=
\begin{array}{cccccccccccc}
  1 &   2 &   3 &   5 &   6 &  10 &  13 &  14  &  15 &  20 &  27 \\
    &  64 &   4 &   7 &   9 &  11 &  16 &  21  &  25 &  30 &  31 \\
    &     &   8 &  12 &  18 &  19 &  24 &  28  &  29 &  37 &  41 \\
    &     &     &  17 &  22 &  23 &  26 &  33  &  34 &  43 &  45 \\
    &     &     &     &  32 &  35 &  38 &  40  &  42 &  47 &  53 \\
    &     &     &     &     &  36 &  39 &  44  &  46 &  52 &  54 \\
    &     &     &     &     &     &  48 &  49  &  50 &  56 &  58 \\
    &     &     &     &     &     &     &  51  &  55 &  59 &  61 \\
    &     &     &     &     &     &     &      &  57 &  60 &  62 \\
    &     &     &     &     &     &     &      &     &  63 &  65
\end{array}  
$$
$$ 
H_1= \left(
(2,2), (1,8), (2,3), (5,6), (3,5) , (7,8), (7,7),  (8,8), (4,9) , (10,10) , (6,11)
\right)
$$
Moreover $C'=\{64\}$ and $C=\{63\}$.

Now we have $e=64$ and $e=63$. We are in Case~2. 
$$
T=
\begin{array}{cccccccccccc}
  1 &   2 &   3 &   5 &   6 &  10 &  13 &  14  &  15 &  20 &  27 \\
    &   4 &   7 &   9 &  11 &  16 &  21 &  25  &  29 &  30 &  31 \\
    &     &   8 &  12 &  18 &  19 &  24 &  28  &  34 &  37 &  41 \\
    &     &     &  17 &  22 &  23 &  26 &  33  &  42 &  43 &  45 \\
    &     &     &     &  32 &  35 &  38 &  40  &  46 &  47 &  53 \\
    &     &     &     &     &  36 &  39 &  44  &  50 &  52 &  54 \\
    &     &     &     &     &     &  48 &  49  &  55 &  56 &  58 \\
    &     &     &     &     &     &     &  51  &  57 &  59 &  61 \\
    &     &     &     &     &     &     &      &  64 &  60 &  62 \\
    &     &     &     &     &     &     &      &     &  63 &  65
\end{array}  
$$
Since $64$ was already a vertical candidate we want to set $H_{1,9}=(10,10)$. For that
purpose we perform a transfer from $(1,9)$ to $(1,10)$.
$$ 
H_1= \left(
(1,8), (1,2), (4,5), (2,4), (6,7) , (6,8), (7,7),  (9,9), (10,10) , (5,10) , (6,11)
\right)
$$
We update $C'=\{63,64\}$ and $C=\emptyset$. Now $63$ is stable and  
$C'=\{64\}$ after the next step. We terminate with the pair
$$
T=
\begin{array}{cccccccccccc}
  1 &   2 &   3 &   5 &   6 &  10 &  13 &  14  &  15 &  20 &  27 \\
    &   4 &   7 &   9 &  11 &  16 &  21 &  25  &  29 &  30 &  31 \\
    &     &   8 &  12 &  18 &  19 &  24 &  28  &  34 &  37 &  41 \\
    &     &     &  17 &  22 &  23 &  26 &  33  &  42 &  43 &  45 \\
    &     &     &     &  32 &  35 &  38 &  40  &  46 &  47 &  53 \\
    &     &     &     &     &  36 &  39 &  44  &  50 &  52 &  54 \\
    &     &     &     &     &     &  48 &  49  &  55 &  56 &  58 \\
    &     &     &     &     &     &     &  51  &  57 &  59 &  61 \\
    &     &     &     &     &     &     &      &  60 &  62 &  64 \\
    &     &     &     &     &     &     &      &     &  63 &  65
\end{array},
$$ 
$$ 
H_1= \left(
(1,8), (1,2), (4,5), (2,4), (6,7) , (6,8), (7,7),  (9,11), (10,10) , (5,10) , (6,11)
\right).
$$
\end{ex}

\medskip

\begin{ex}
\label{ex3}
Let 
$$ 
T=
\begin{array}{cccccccccccc}
 22 &  15 &  37 &  21 &  44 &  11 &  62 &  14  &  57 &  13 &  34 \\
    &   1 &   2 &   4 &   5 &   7 &  10 &  17  &  19 &  28 &  32 \\
    &     &   3 &   6 &   8 &   9 &  23 &  27  &  31 &  38  &  41 \\
    &     &     &  12 &  16 &  18 &  24 &  29  &  36 &  40 &  53 \\
    &     &     &     &  20 &  25 &  26 &  30  &  46 &  47 &  54 \\
    &     &     &     &     &  33 &  35 &  42  &  48 &  50 &  55 \\
    &     &     &     &     &     &  39 &  43  &  49 &  51 &  60 \\
    &     &     &     &     &     &     &  45  &  52 &  58 &  63 \\
    &     &     &     &     &     &     &      &  56 &  59 &  64 \\
    &     &     &     &     &     &     &      &     &  61 &  65
\end{array}.
$$
After the application of SPLIT~1 we obtain 
$$ T= 
\begin{array}{cccccccccccc}
 22 &   1 &   2 &   4 &   5 &   7 &  10 &  13  &  14 &  28 &  32 \\
    &  15 &   3 &   6 &   8 &   9 &  17 &  19  &  31 &  34 &  41 \\
    &     &  37 &  11 &  16 &  18 &  23 &  27  &  36 &  38  &  53 \\
    &     &     &  12 &  20 &  24 &  26 &  29  &  40 &  47 &  54 \\
    &     &     &     &  21 &  25 &  30 &  42  &  46 &  50 &  55 \\
    &     &     &     &     &  33 &  35 &  43  &  48 &  51 &  57 \\
    &     &     &     &     &     &  39 &  44  &  49 &  58 &  60 \\
    &     &     &     &     &     &     &  45  &  52 &  59 &  63 \\
    &     &     &     &     &     &     &      &  56 &  61 &  64 \\
    &     &     &     &     &     &     &      &     &  62 &  65
\end{array}  
$$
and the first row of the shifted hook tabloid is 
$$ 
H_1 = \left(
(1,1) , (2,2) , (3,3) , (1,4) , (5,5) , (1,10) , (1,9) , (7,8) , (1,10) , (10,10) , (6,11) 
\right).
$$
We apply SPLIT~2 to the pair. Observe that $i'=3$ since $H_{1,j}=(j,j)$ for $1 \le j \le 3$ 
and the entries $22$, $15$ are unstable.

If we perform jeu de taquin with $37$, $15$ and $22$ and with respect to the relative 
next row we observe that $22$ does not change row and thus $h=1$. Therefore we are 
in Case~1 of SPLIT~2. We perform jeu de taquin with $37$, $15$ and $22$ in that order
and with respect to the last cell in the current row of the forward path in the course of constructing 
$U$. This results in 
$$ T =
\begin{array}{cccccccccccc}
  1 &   2 &   4 &   5 &   7 &  10 &  13 &  14  &  22 &  28 &  32 \\
    &   3 &   6 &   8 &   9 &  15 &  17 &  19  &  31 &  34 &  41 \\
    &     &  11 &  37 &  16 &  18 &  23 &  27  &  36 &  38  &  53 \\
    &     &     &  12 &  20 &  24 &  26 &  29  &  40 &  47 &  54 \\
    &     &     &     &  21 &  25 &  30 &  42  &  46 &  50 &  55 \\
    &     &     &     &     &  33 &  35 &  43  &  48 &  51 &  57 \\
    &     &     &     &     &     &  39 &  44  &  49 &  58 &  60 \\
    &     &     &     &     &     &     &  45  &  52 &  59 &  63 \\
    &     &     &     &     &     &     &      &  56 &  61 &  64 \\
    &     &     &     &     &     &     &      &     &  62 &  65
\end{array}  .
$$
Since $37$ moved to cell $(3,4)$ we set $H_{1,2}=(3,4)$, since $15$ moved to 
$(2,6)$ we set $H_{1,1}=(2,6)$ and since $22$ moved to $(1,9)$ we set 
$H_{1,3}=(1,9)$.
$$
H_1 = \left( 
(2,6) , (3,4) , (1,9) , (1,4) , (5,5) , (1,10) , (1,9) , (7,8) , (1,10) , (10,10) , (6,11)
\right)
$$
After the application of SPLIT~3 we obtain 
$$ 
T=
\begin{array}{cccccccccccc}
  1 &   2 &   4 &   5 &   7 &  10 &  13 &  14  &  22 &  28 &  32 \\
    &   3 &   6 &   8 &   9 &  15 &  17 &  19  &  31 &  34 &  41 \\
    &     &  11 &  12 &  16 &  18 &  23 &  27  &  36 &  38  &  53 \\
    &     &     &  20 &  21 &  24 &  26 &  29  &  40 &  47 &  54 \\
    &     &     &     &  25 &  30 &  35 &  42  &  46 &  50 &  55 \\
    &     &     &     &     &  33 &  37 &  43  &  48 &  51 &  57 \\
    &     &     &     &     &     &  39 &  44  &  49 &  58 &  60 \\
    &     &     &     &     &     &     &  45  &  52 &  59 &  63 \\
    &     &     &     &     &     &     &      &  56 &  61 &  64 \\
    &     &     &     &     &     &     &      &     &  62 &  65
\end{array}  
$$
and 
$$
H_1= \left( 
(2,6) , (1,9) , (4,4) , (1,4) , (6,7) , (1,10) , (1,9) , (7,8) , (1,10) , (10,10) , (6,11) 
\right).
$$
\end{ex}

\medskip

\begin{ex}
\label{ex4}
Let 
$$
T=
\begin{array}{cccccccccccc}
 20 &   9 &  64 &   7 &  33 &  54 &  65 &  13  &  61 &  50 &  46 \\
    &   1 &   2 &   4  &   8 &  11 &  14 &  16  &  22 &  24 &  37 \\
    &     &   3 &   5 &  10 &  12 &  18 &  23  &  26 &  29 &  41 \\
    &     &     &   6 &  15 &  17 &  21 &  25  &  30 &  36 &  43 \\
    &     &     &     &  19 &  27 &  28 &  34  &  35 &  39 &  53 \\
    &     &     &     &     &  31 &  32 &  38  &  42 &  44 &  55 \\
    &     &     &     &     &     &  40 &  45  &  48 &  51 &  58 \\
    &     &     &     &     &     &     &  47  &  49 &  52 &  59 \\
    &     &     &     &     &     &     &      &  56 &  57 &  62 \\
    &     &     &     &     &     &     &      &     &  60 &  63
\end{array}.
$$
After the application of SPLIT~1 we obtain the pair 
$$
T=
\begin{array}{cccccccccccc}
 20 &   1 &   2 &   4  &   8 &  11 &  13 &  16  &  22 &  24 &  37 \\
    &   9 &   3 &   5 &  10 &  12 &  14 &  23  &  26 &  29 &  41 \\
    &     &  64 &   6 &  15 &  17 &  18 &  25  &  30 &  36 &  43 \\
    &     &     &   7 &  19 &  21 &  28 &  34  &  35 &  39 &  46 \\
    &     &     &     &  33 &  27 &  32 &  38  &  42 &  44 &  53 \\
    &     &     &     &     &  31 &  40 &  45  &  48 &  50 &  55 \\
    &     &     &     &     &     &  54 &  47  &  49 &  51 &  58 \\
    &     &     &     &     &     &     &  65  &  52 &  57 &  59 \\
    &     &     &     &     &     &     &      &  56 &  60 &  62 \\
    &     &     &     &     &     &     &      &     &  61 &  63
\end{array}  
$$
and 
$$
H_1= \left(
(1,1) , (2,2) , (3,3) , (4,4) , (5,5) , (1,8) , (7,7) , (8,8) , (5,9) , (10,10) , (4,11)
\right).
$$

Next we apply SPLIT~2. Observe that $i'=4$ (since $7$ is stable) and that 
$$
U=
\begin{array}{cccccccccccc}
  1 &   2 &   4  &   7 &   8 &  11 &  13 &  16  &  22 &  24 &  37 \\
    &   3 &   5 &  20 &  10 &  12 &  14 &  23  &  26 &  29 &  41 \\
    &     &   6 &   9 &  15 &  17 &  18 &  25  &  30 &  36 &  43 \\
    &     &     &  64 &  19 &  21 &  28 &  34  &  35 &  39 &  46 \\
    &     &     &     &  33 &  27 &  32 &  38  &  42 &  44 &  53 \\
    &     &     &     &     &  31 &  40 &  45  &  48 &  50 &  55 \\
    &     &     &     &     &     &  54 &  47  &  49 &  51 &  58 \\
    &     &     &     &     &     &     &  65  &  52 &  57 &  59 \\
    &     &     &     &     &     &     &      &  56 &  60 &  62 \\
    &     &     &     &     &     &     &      &     &  61 &  63
\end{array}.  
$$
Thus $h=4$ and no backward paths with respect to 
the $4$-th row of a horizontal candidate strictly below the $4$-th row 
ends weakly left of $64$. Therefore we are in Case~1 of SPLIT~2. 
We set $U=T$. Since $7$ original ended in $(4,4)$ we set $H_{1,4}=(4,4)$, 
since $64$ is in $(4,4)$ we set $H_{1,3}=(4,4)$, 
since $9$ is in $(3,4)$ we set $H_{1,2}=(3,4)$ 
and since $20$ is in $(2,4)$ we set $H_{1,1}=(2,4)$.
$$
H_1 = \left(
(2,4) , (3,4) , (4,4) , (4,4) , (5,5) , (1,8) , (7,7) , (8,8) , (5,9) , (10,10) , (4,11)
\right)
$$
After the performance of SPLIT~3 we obtain the pair 
$$
T=
\begin{array}{cccccccccccc}
  1 &   2 &   4  &   7 &   8 &  11 &  13 &  16  &  22 &  24 &  37 \\
    &   3 &   5 &   9 &  10 &  12 &  14 &  23  &  26 &  29 &  41 \\
    &     &   6 &  15 &  17 &  18 &  20 &  25  &  30 &  36 &  43 \\
    &     &     &  19 &  21 &  28 &  32 &  34  &  35 &  39 &  46 \\
    &     &     &     &  27 &  31 &  38 &  42  &  44 &  50 &  53 \\
    &     &     &     &     &  33 &  40 &  45  &  48 &  54 &  55 \\
    &     &     &     &     &     &  47 &  49  &  51 &  57 &  58 \\
    &     &     &     &     &     &     &  52  &  56 &  59 &  62 \\
    &     &     &     &     &     &     &      &  60 &  61 &  64 \\
    &     &     &     &     &     &     &      &     &  63 &  65
\end{array}
$$
and 
$$
H_1 = \left(
(2,3) , (3,7) , (3,3) , (5,5) , (6,9) , (1,8) , (3,7) , (9,11) , (10,11) , (10,10) , (4,11)
\right).
$$
\end{ex}

\medskip

\begin{ex}
\label{ex5}
We consider 
$$
T=
\begin{array}{cccccccccccc}
 26 &   4 &  13 &  31 &  47 &  24 &  58 &  65  &  53 &  25 &  60 \\
    &   1 &   2 &   3 &   6 &   8 &   9 &  11  &  14 &  22 &  29 \\
    &     &   5 &   7 &  10 &  16 &  17 &  20  &  23 &  30 &  36 \\
    &     &     &  12 &  15 &  18 &  19 &  21  &  33 &  41 &  43 \\
    &     &     &     &  27 &  28 &  32 &  37  &  38 &  42 &  49 \\
    &     &     &     &     &  34 &  35 &  39  &  40 &  46 &  51 \\
    &     &     &     &     &     &  44 &  45  &  48 &  52 &  54 \\
    &     &     &     &     &     &     &  50  &  55 &  57  &  61 \\
    &     &     &     &     &     &     &      &  56 &  59 &  62 \\
    &     &     &     &     &     &     &      &     &  63 &  64
\end{array}.
$$
After the application of SPLIT~1 we obtain 
$$ T=
\begin{array}{cccccccccccc}
 26 &   1 &   2 &   3 &   6 &   8 &   9 &  11  &  14 &  22 &  29 \\
    &   4 &   5 &   7 &  10 &  16 &  17 &  20  &  23 &  25 &  36 \\
    &     &  13 &  12 &  15 &  18 &  19 &  21  &  30 &  41 &  43 \\
    &     &     &  31 &  24 &  28 &  32 &  33  &  38 &  42 &  49 \\
    &     &     &     &  25 &  34 &  35 &  37  &  40 &  46 &  51 \\
    &     &     &     &     &  47 &  39 &  45  &  48 &  52 &  54 \\
    &     &     &     &     &     &  44 &  50  &  53 &  57  &  60 \\
    &     &     &     &     &     &     &  58  &  55 &  59 &  61 \\
    &     &     &     &     &     &     &      &  56 &  62 &  64 \\
    &     &     &     &     &     &     &      &     &  63 &  65
\end{array}  
$$
and 
$$
H_1= \left( 
(1,1) , (2,2) , (3,3) , (4,4) , (3,5) , (6,6) , (1,9) , (8,8) , (6,9) , (6,10) , (10,11)
\right).
$$

Next we apply SPLIT~2. Observe that $i'=2$ and that 
$$
U=
\begin{array}{cccccccccccc}
  1 &   2 &   3 &   6 &   8 &   9 &  11 &  14  &  22 &  25 &  29 \\
    &   4 &   5 &   7 &  10 &  16 &  17 &  20  &  23 &  26 &  36 \\
    &     &  13 &  12 &  15 &  18 &  19 &  21  &  30 &  41 &  43 \\
    &     &     &  31 &  24 &  28 &  32 &  33  &  38 &  42 &  49 \\
    &     &     &     &  27 &  34 &  35 &  37  &  40 &  46 &  51 \\
    &     &     &     &     &  47 &  39 &  45  &  48 &  52 &  54 \\
    &     &     &     &     &     &  44 &  50  &  53 &  57  &  60 \\
    &     &     &     &     &     &     &  58  &  55 &  59 &  61 \\
    &     &     &     &     &     &     &      &  56 &  62 &  64 \\
    &     &     &     &     &     &     &      &     &  63 &  65
\end{array}.
$$
Thus $h=2$ and the backward path with respect to the $2$-nd row of the 
horizontal candidate $13$ ends weakly left of $26$. We are in Case~2 of
SPLIT~2 and  set $T=U$ and 
$$
H_1 = \left(
(2,10) , (2,2) , (3,3) , (4,4) , (3,5) , (6,6) , (1,9) , (8,8) , (6,9) , (6,10) , (10,11)
\right).
$$
After the application of SPLIT~3 we finally obtain 
$$
T=
\begin{array}{cccccccccccc}
  1 &   2 &   3 &   6 &   8 &   9 &  11 &  14  &  22 &  25 &  29 \\
    &   4 &   5 &   7 &  10 &  16 &  17 &  20  &  23 &  26 &  36 \\
    &     &  12 &  13 &  15 &  18 &  19 &  21  &  30 &  41 &  43 \\
    &     &     &  24 &  27 &  28 &  32 &  33  &  38 &  42 &  49 \\
    &     &     &     &  31 &  34 &  35 &  37  &  40 &  46 &  51 \\
    &     &     &     &     &  39 &  44 &  45  &  48 &  52 &  54 \\
    &     &     &     &     &     &  47 &  50  &  53 &  57  &  60 \\
    &     &     &     &     &     &     &  55  &  56 &  59 &  61 \\
    &     &     &     &     &     &     &      &  58 &  62 &  64 \\
    &     &     &     &     &     &     &      &     &  63 &  65
\end{array}  
$$
and 
$$
H_1 = \left(
(2,10) , (3,4) , (3,3) , (2,4) , (5,5) , (1,9) , (7,7) , (5,8) , (9,9) , (6,10) , (10,11)
\right).
$$
\end{ex}

\section{The Algorithm MERGE}
\label{inv}

Now we describe the Algorithm MERGE that merges a pair of a shifted standard tableau and 
a shifted hook tabloid to a shifted tabloid. This will turn out to be the inverse of the Algorithm SPLIT. 
Again we have to introduce two routines on a partial shifted hook tabloid $H$ of shape $\lambda$.

\medskip

{\bf A reshift from cell $(i,j')$ to cell $(i,j)$ in a partial shifted hook tabloid.} 
Let $(i,j)$, $(i,j')$ 
be two cells in $H$, $j \le j'$. We define the term {\it reshift from 
$(i,j')$ to $(i,j)$ in $H$}. The output of this operation is another
shifted tabloid $H'$ which coincides with $H$ except for the 
cells $(i,k)$, $j \le k \le j'$. Let $j < k \le j'$. If 
$H_{i,k-1}=(i',k-1)$ or  
$H_{i,k-1}=(k,j')$ then set $H'_{i,k}=(i'+1,k)$  or 
$H'_{i,k}=(k+1,j'+1)$, respectively. 
If $H_{i,k-1}=(i,l)$, $l \ge k$, then
set $H'_{i,k}=(i,l)$. The cell $(i,j)$ in $H'$ remains 
empty. We denote $H'$ by $\rshift(H,i,j',j)$. 
For an example for a reshift from $(1,6)$ to $(1,2)$ read Figure~\ref{shift} from right to left.

\smallskip

Observe that $\shift$ and $\rshift$ are inverse to each other in the following sense:
Let $i \le j \le j' \le \lambda_i + i -1$ and $H$ be a partial shifted hook tabloid
in which only cell $(i,j)$ is empty. Then 
$$
H = \rshift ( \shift(H,i,j,j'), i , j',j).
$$
If $H'$ is a partial shifted hook tabloid where only cell $(i,j')$ is empty then 
$$
H' = \shift ( \rshift (H, i , j' , j), i, j , j').
$$

\smallskip

In general the output of a reshift in a partial shifted hook tabloid need not to be a 
partial shifted hook tabloid for $(i'+1,k)$, respectively $(k+1,j'+1)$, 
need not to be a cell in the shifted Ferrers diagram of shape $\lambda$ if 
$(i',k-1)$, respectively $(k,j')$, is a cell in this diagram. However, 
in the following two situations the reshift is indeed a partial shifted hook tabloid 
and we apply the reshift 
only if one of these situations encounters. 

\begin{enumerate}

\item 
Let $T$ be a shifted tabloid and $H$ an accompanying partial shifted hook tabloid. Let 
$e$ be a vertical candidate 
in $T$ with respect to $H$ and the $i$-th row (meaning that the pointer to $e$
is in the $i$-th row of $H$)and perform reverse jeu de
taquin ($i$ being the fixed row)
with $e$ in $T$ and with 
respect to a set $D$. Let $(g,j)$ denote the cell of $e$ in the output tabloid and 
$(g',j')$ denote the cell of $e$ in $T$. Suppose that $e$ is the maximal vertical candidate 
in $T$ under the vertical candidates in the
rows $h$, $g \le h \le g'$, with respect to the backward paths order 
of $T$. Then the reshift from cell $(i,g'-1)$ to cell $(i,g-1)$ in $H$
produces another partial shifted hook tabloid.

\medskip

\item 
Again let $T$ be a shifted tabloid and $H$ an accompanying partial shifted hook tabloid.
Furthermore we assume in this situation that there exists no vertical candidate 
in $T$ with respect to $H$. Let $e$ be a horizontal candidate 
in $T$ with respect to $H$ and the $i$-th row and perform reverse jeu de
taquin ($i$ being the fixed row) with $e$ in $T$ and with 
respect to a set $D$. Let $(g,j)$ denote the cell of $e$ in the output tabloid and $(g',j')$ 
denote the cell of $e$ in $T$. Suppose that $e$ is the smallest horizontal candidate  in $T$ 
under the horizontal candidates in the columns
$k$, $j \le k \le j'$, with respect to the backward paths order of 
$T$. Then the reshift from cell $(i,j')$ to cell $(i,j)$ in $H$ produces another partial shifted 
hook tabloid. 

\end{enumerate}

\smallskip

(In order to see that use Lemma~\ref{0}.)

\smallskip

As jeu de taquin in SPLIT is mostly followed by a shift in the accompanying shifted hook tabloid, 
reverse jeu de taquin in MERGE is mostly followed by a reshift. Let $T$ be a shifted tabloid, 
$H$ a 
shifted hook tabloid, $e$ a candidate and $D$ a set of cells. If $e$ is a horizontal 
candidate then $(T',H')=\rjs_D(T,H,e)$ is obtained as follows: $T'=\rjt_D(T,e)$, 
$H'=\rshift(H,i,q,j')$ and $H'_{i,j'}=(i',j')$, where $q$ is the column of $e$ in $T$
and $c_{T'}(e)=(i',j')$. If $e$ is a vertical candidate then $(T',H')=\rjs_D(T,H,e)$ is 
obtained as follows: $T'=\rjt_D(T,e)$,
$H'=\rshift(H,i,p-1,i'-1)$ and $H'_{i,i'-1}=(i',j')$, where $p$ is the row of $e$ in $T$
and $c_{T'}(e)=(i',j')$. Observe that $\js$ and $\rjs$ are each other's respective inverses.

\medskip

{\bf A retransfer from cell $(i,k)$ to cell $(i,j)$ in a partial shifted hook tabloid.}
Let $j \le k \le r$ and $H_{i,k}=(i',k')$ with either $i'=i$ or $k'=k$. We define 
the term {\it retransfer from cell $(i,k)$ to cell $(i,j)$ in $H$.} The output of this operation 
is again another partial shifted hook tabloid $H'$ which coincides with $H$ except for the 
cells $(i,l)$, $j \le l \le k$. If $k'=k$ and $i'-i \ge k - j$ let $H'_{i,j}=(i'+j-k,j)$, 
otherwise $H'_{i,j}=(i,k'+i-i')$. For $j < l \le k$ set $H'_{i,l}=(l,l)$. We 
denote $H'$ by $\rtrans(H,i,k,j)$. For an example of a retransfer from $(1,6)$ to $(1,2)$ read Figure~\ref{trans}
from right to left. Again, whenever we will apply a retransfer from $(i,k)$ to
$(i,j)$ in $H$, we have $H_{i,l}=(l+1,l+1)$ for $j \le l < k$ as in the example.

\smallskip

Observe that $\trans$ and $\rtrans$ are inverse to each other in the following sense.
Let $i \le j \le k \le \lambda_i + i -1$. If $H$ is a shifted hook tabloid such that 
$H_{i,j}$ is in the same row or column as $(i,j)$ and with $H_{i,l}=(l,l)$ for 
$j < l \le k$ then 
$$
H = \rtrans ( \trans (H, i, j , k), i, k ,j).
$$
If $H'$ is a shifted hook tabloid such that $H'_{i,k}$ is in the same row or column
as $(i,k)$ and with $H'_{i,l}=(l+1,l+1)$ for $j \le l < k$ then 
$$
H' = \trans ( \rtrans ( H', i , k ,j), i , j , k).
$$

\medskip

Now we are in the position to state the Algorithm MERGE. It is divided into
$r$ steps, where in the $i$-th step we perform reverse jeu de taquin with the
entries in the $i$-th row. Immediately before we perform the $i$-th step of MERGE to a
pair $(T,H)$, $T$ is ordered up to $(i,i)$, the first $(i-1)$ rows of $H$  are
empty and the last $r-i+1$ rows form a shifted hook tabloid. Within a row $i$,
MERGE is divided into $3$ steps. The $3$ parts of MERGE in a fixed row $i$ are numbered in 
reverse order to emphasize the connection between SPLIT~$k$ and MERGE~$k$, $k=1,2,3$. 
All the  horizontal and vertical candidates below are meant to be with respect
to the $i$-th row, i.e. their pointers are in the $i$-th row of the
current shifted hook tabloid. The marker '($*$)' in MERGE~3 is needed in the
proof that SPLIT~3 and MERGE~3 are inverse to each other in Subsection~\ref{proof}.6. 

\medskip

\fbox{ \parbox{14cm}{

{\bf MERGE~3.} 
Set $z=i$. Repeat the following until there exists no vertical candidate
strictly below the $z$-th row.

\medskip

[ Let $e$ be the vertical candidate strictly below the $z$-th row, 
which is maximal with respect to the backward paths order.

\smallskip

If $c_T(e)=(h,h)$: [ Let $h'$ be minimal such  that $T_{j,j}$ is a vertical 
unstable candidate for $h' \le j < h$.  If $e=T_{h,h}$ is unstable and the 
backward path of smallest horizontal candidate $e'$ strictly below the $h$-th row contains 
$(h+1,h+1)$ set $(T,H)=\rjs_{(h+1,h+1)}(T,H,e')$ ($*$) and
$H=\rtrans(H,i,h,h'-1)$, otherwise set $H=\rtrans(H,i,h,h')$. Furthermore $e=T_{h',h'}$. ]

\smallskip

If $e$ is a vertical candidate repeat the following:

[If $c_T(e)=(k,k)$, $k \not=z+1$, set $(T,H)=\rjs_{(k-1,k)} (T,H,e)$.
Set $(T,H) =\rjs_{MD \cup z+1} (T,H,e)$.]

until either 
\begin{enumerate}
\item $e$ is in row $z+1$,
\item there exists a vertical candidate strictly below the row of $e$,
\item $c_T(e)=(l,l)$ for an $l$ and the backward path of a horizontal candidate 
includes $(l+1,l+1)$,
\item $c_T(e)=(l,l)$ for an $l$ and $T_{l-1,l-1}$ is an unstable vertical candidate. 
\end{enumerate}

\smallskip

If $e$ is in row $z+1$ set $z=z+1$. ]

\smallskip
}}

\smallskip

Observe that in case the if-condition (If $c_T(e)=(h,h)$: [$\dots$])
at the beginning of a step of MERGE~3 is true, there exists no vertical candidate 
strictly below the $h$-th row by the maximality of $e$.

\smallskip

\fbox{ \parbox{14cm}{

{\bf MERGE~2.}
Let $i''$ be such that for $i+1 \le h \le i''$ there exists 
a vertical candidate entry in the $h$-th row and all other candidates are horizontal.

\medskip

We distinguish between two cases. In Case~2 we consider the case when 
the smallest horizontal candidate strictly below the $i''$-th row is smaller than 
the smallest vertical candidate or $i''=i$.

\smallskip

{\it Case 1.} [ Set $a=$True. Let $H_{i,i''}=(h,k)$. We define $T_{i,k-h+i}$
to be a vertical candidate in row $i$.

Repeat for $g=i$ to $g=i''$: [ If $a=$True
let $e$ be the vertical candidate in row $g$ or $\min(g+1,i'')$
that is maximal with respect to the backward paths order. If $a=$ False let
$e$ be the vertical candidate in row $g$. Let $T=\rjt_{(g,g)}(T,e).$
If the backward path includes $(g+1,g+1)$ set $a=$False. If $e$ was the candidate in the 
$(g+1)$-st row  and 
$T_{g,k}$ is the vertical candidate in row $g$, then let
$H_{i,g}=(g+1,k+1)$. ]

For $i \le g \le i''$, let $H_{i,g}=(g,g)$. ]

\smallskip

{\it Case 2.} [ Repeat for $g=i$ to $g=i''-1$: [ Let $e$ be the entry in cell $H_{i,g}$ of $T$ and
set $T=\rjt_{MD} (T,e)$.  ]
Let $i \le i' \le i''$, be such that in the previous step no entry 
moved to $(i',i')$. Set $H=\rtrans(H,i,i'',i')$.
For every $g$, $i \le g < i'$, let $H_{i,g}=(g,g)$. ]
}}

\smallskip

\fbox{ \parbox{14cm}{

{\bf MERGE~1.} Repeat for $j=i$ to $j=\lambda_i + i -1$:
[ Let $e$ be the minimal horizontal candidate with respect to the 
backward paths order in $T$. Set $(T,H)=\rjs_{(i,j)} (T,H,e)$ and erase the
entry in cell $(i,j)$ of $H$.]
}}

\medskip

Let $T$ be a shifted tabloid which is ordered up to $(i,i)$ 
and $H$ a partial shifted hook tabloid, where no entry in the $i$-th row 
is empty. Let $(T',H')$ denote the 
pair we obtain after the application of MERGE~3, MERGE~2 and MERGE~1 
to the $i$-th row of $(T,H)$. We denote 
$T'=\smerge(T,H,i)$ and note that $H'$ is equal to $H$ with the entries in the 
$i$-th row deleted.
We invite the reader to apply MERGE~3, MERGE~2 and MERGE~1 to the first rows of the 
output pairs in the examples in Section~\ref{examples} in order to find out that 
these applications result in the input pairs.

We are 
finally in the position to formulate the Algorithm MERGE. The input is 
a shifted standard tableaux $T$ of shape $\lambda$ and a shifted 
hook tabloid $H$ of shape $\lambda$.

\medskip

\fbox{ \parbox{14cm}{
{\bf MERGE.} Set $T'=T$ and $H'=H$. Repeat for $i=1$ to $i=r$: [  
$T'=\smerge(T',H',i)$ and delete the $i$-th row of $H'$.]
}}

\medskip

The output tabloid $T'$ of the algorithm is denoted by $\smerge(T,H)$.
Note that the output tabloid $H'$ is  empty.

\medskip

In the remaining section we aim to show the following theorem. 

\smallskip

\begin{theo} 
\label{key}
Let $\lambda=(\lambda_1,\lambda_2,\dots, \lambda_r)$ be a partition with distinct 
components
and $1 \le i \le r$.
\begin{enumerate}
\item Let $T$ be a shifted tabloid of shape $\lambda$ which is ordered up to $(i+1,i+1)$ and 
$H$ a partial shifted hook tabloid of shape $\lambda$, where the cells of the $i$-th row are empty.
Let $T'=\ssplit(T,H,i)$ and $H'=\hsplit(T,H,i)$. Then 
$T=\smerge(T',H',i)$.

\item Let $T$ be a shifted tabloid of shape $\lambda$ which is ordered up to $(i,i)$ 
and $H$ a partial shifted hook tabloid of shape $\lambda$, where no cell in the 
$i$-th row is empty. Let $T'=\smerge(T,H,i)$ and $H'$ denote the partial shifted hook 
tabloid $H$ with the entries in the $i$-th row deleted.
Then $T=\ssplit(T',H',i)$ and $H=\hsplit(T',H',i)$.
\end{enumerate}
\end{theo}

This theorem, once it is proved, shows that SPLIT is the desired algorithm, 
i.e. a bijection from $T_{\lambda}$ onto $S_{\lambda} \times H_{\lambda}$, for  
MERGE is its inverse.

In the following we fix a strict partition
$\lambda=(\lambda_1,\lambda_2,\dots, \lambda_r)$ and 
a row $i$, $1 \le i \le r$.

\section{The Proof of Theorem~\ref{key}}
\label{proof}

\medskip

\begin{center} \sc 6.1. Variantes of SPLIT~1 and SPLIT~2  \end{center}

\medskip

We modify SPLIT~1 and SPLIT~2 such that after this modification MERGE~$k$ 
is the inverse of SPLIT~$k$, $k=1,2,3$. The original version
of SPLIT and the modification of SPLIT, where we replace SPLIT~$k$, 
$k=1,2$, by their modifications described below, are equivalent, since we only change the 
order of some commuting steps and swap the beginning of SPLIT~3 to the end of SPLIT~2. 

\medskip

We start with a modification of SPLIT~2. Consider the following Algorithm 
POST-SPLIT~2, which we apply to a pair to which we have applied 
SPLIT~2.

\medskip

\fbox{ \parbox{14cm}{
{\bf POST-SPLIT~2.} We continue with Case~$k$, $k=1,2$, if we were in Case~$k$ in SPLIT~2.

\medskip

{\it Case 1.} [ If $T_{i',i'}$ is a vertical candidate, let $i''$, $i' \le i''$, be 
maximal such that for $i' \le j < i''$, $T_{j,j}$ is unstable and 
$T_{j+1,j+1}$ is a horizontal candidate and
set $H=\trans(H,i,i',i'')$. Otherwise omit this step and set $i''=i'$. ]

\smallskip

{\it Case 2.} [ If $T_{i'+1,i'+1}$ is a horizontal candidate, 
let $i''$, $i'+1 \le i''$, be maximal such that for $i'+1 \le j \le i''$,
$T_{j,j}$ is unstable, $T_{j+1,j+1}$ is a horizontal candidate and 
set $H=\trans(H,i,i',i'')$. Furthermore set $(T,H)=\js_{(i''+1,k)}(T,H,T_{i''+1,i''+1})$, 
where $(i''+1,k)$ was the last cell in the forward path of $T_{i''+1,i''+1}$ in the 
$(i''+1)$-st row. 
Otherwise omit this step and set 
$i''=i'$. ]

\smallskip

In both cases: $U=T$. If we were in Case~2 and $i'' \not= i'$ set
$U=\jt_{i''+2}(U,T_{i''+1,k})$. Repeat for $g=i''$ down to $g=i'+1$: [ Let $e$ be the vertical 
candidate in row $g$ of $T$ and set $U=\jt_{g+1}(U,e)$.  ]  Reject $U$.

Repeat for $g=i''$ down to $g=i'+1$: [ Let $e$ be the vertical candidate in row 
$g$ of $T$ and set $(T,H)=\js_{(g,l)} (T,H,e)$, where $(g,l)$ was the last cell in the 
forward path of $e$ in the $g$-th row in the shifted tabloid $U$ in the previous step.
{ \small (Note that the forward path of $e$ ends in $(g,l)$.)} ]
}}

\medskip

Again the following  holds for $(T,H)$ after the application of POST-SPLIT~2: 
We were in Case~2 of POST-SPLIT~2 if and only if 
either the smallest horizontal candidate strictly below the $i''$-th row
is smaller than every vertical candidate or no vertical candidate exists.

\smallskip 

Furthermore observe that the application of SPLIT~2 and SPLIT~3 is 
equivalent to the application of SPLIT~2, POST-SPLIT~2 and 
SPLIT~3, for POST-SPLIT~2 is the beginning of SPLIT~3 roughly speaking.

\medskip

{\bf Touching from above.} 
Let $(T,H)$ denote a pair to which we apply SPLIT~2, POST-SPLIT~2 and SPLIT~3.
Assume that we are in Case~2 in a certain step of SPLIT~3.
Then we say {\it $e$ touches $e'$ from above} and for 
$h < j \le h'$, $T_{j-1,j-1}$ touches $T_{j,j}$ from above. 
If $e$ is a horizontal candidate before the application of Case~2 we say that $(h-1,h-1)$ 
is the {\it place of change for $e$} and $(j,j)$ is the place of change 
for $T_{j,j}$, $h \le j < h'$.  If $e$ is already a vertical candidate before
the application of Case~2 we 
say that $(j,j)$ is the place of change for $T_{j,j}$, $h \le j \le h'$. 
Furthermore: If we are in Case~1 of SPLIT~2 let $T_{g,g}$, $i \le g \le i''$, be the 
exceptional entries, where $T$ denotes the shifted tabloid at the beginning of SPLIT~2.
If we are in Case~2 of SPLIT~2 let $T_{g,g}$, $i \le g \le i''$ and $g \not= i', i'+1$,
be the exceptional entries, $T$ again being the shifted tabloid at the beginning of SPLIT~2.
If $i' \not= i''$ and we are in Case~2 of 
SPLIT~2 then the vertical candidate in row 
$i'+1$ after the application of SPLIT~2 is defined to be an entry which touches from above 
with place of change $(i'+1,i'+1)$.

\smallskip

We replace SPLIT~1 by the following variante. Again the input is a shifted tabloid $T$, 
which is ordered 
up to $(i+1,i+1)$ and a partial shifted hook tabloid $H$, where the $i$-th row is empty.

\smallskip

\fbox{ \parbox{14cm}{
{\bf Variante of SPLIT~1.} 
Repeat for $j=\lambda_i + i - 1$ down to $j=i$: We set 
$(T,H)=\js_D(T,H,T_{i,j})$, where the set $D$ is defined as follows.
If $T_{i,j}$ is exceptional or touched from above let $D=MD$. If
$T_{i,j}$ is neither exceptional nor touched from above but 
touches from above, let $D=\{\rho\}$, where $\rho$ is its place of 
change. In all other cases $D=\emptyset$.]
}}

\medskip

The change of SPLIT~1 forces us to replace SPLIT~2 by the following slight 
modification. This is because for a pair $(T,H)$ which falls into Case~2 
in the original version of
SPLIT~2, the entry $T_{i',i'}$ either moves to
a stable position or to its place of change if we apply the Variante of SPLIT~1. 

\medskip

\fbox{ \parbox{14cm}{

{\bf Variante of SPLIT~2.} 
Choose $i'$, $i \le i' \le r$, maximal such that $H_{i,k}=(k,k)$ and 
$T_{k-1,k-1}$ is unstable (i.e. $T_{k-1,k-1} > T_{k-1,k}$) for $i <  k \le
i'$.  If $i'=i$ stop.
 
Set $U=T$. Repeat for $g=i'$ down to $g=i$: Set $U=\jt_{g+1}(U,T_{g,g})$.

\medskip

We distinguish between two cases. 
Let $h$ be minimal, $i \le h \le i'$, such that $T_{h,h}$
is not in the $(h+1)$-st row of $U$. We continue with Case 2  below if 
$h$ does not exist or $T_{i',i'}$ is stable and  the backward path with respect to the 
$i'$-th row  of a horizontal candidate in row $i'+1$ 
or below ends weakly to the left of the cell of $T_{i'-1,i'-1}$ in the 
$i'$-st row in $U$. 
In all other cases we continue with Case~1.

\medskip

Reject the tabloid $U$ we constructed so far in this step.

\medskip

{\it Case 1.} 
Repeat for $g=i'$ down to $g=h+1$:
[ Let $k$ be such that either $(g,k)$ or $(g+1,k)$ is the endcell of the
forward path of $T_{g,g}$ in the procedure for constructing $U$.
Set $T=\jt_{(g,k)} (T, T_{g,g})$ and $H_{i,g-1}=(g,k)$. ]

Let $(h,k)=\cjt(T,T_{h,h})$ and set $T=\jt(T,T_{h,h})$.
If $h - k  \le i - i'$ let $H_{i,i'}=(i,i - h + k)$ otherwise let $H_{i,i'}=(i'+h-k,i')$.

Repeat for $g=h-1$ down to $g=i$: [ Let $H_{i,g}=\cjt_{g+1}(T,T_{g,g})$ and 
$T=\jt_{g+1}(T,T_{g,g})$.]

\smallskip

{\it Case 2.} If $T_{i',i'}$ is unstable, set $i'=i'+1$.
Repeat for $g=i'-1$ down to $g=i$: 
[ Let $H_{i,g}=\cjt_{g+1} (T,T_{g,g})$ and $T=\jt_{g+1} (T, T_{g,g})$. ]

Now we apply POST-SPLIT~2.
}}

\medskip

Check that if we replace SPLIT~1 and SPLIT~2 in SPLIT by the variantes above, 
this variante of SPLIT is equivalent to the original version. In the 
following SPLIT denotes this variante and SPLIT~$k$, $k=1,2,3$, its parts. 
Furthermore
SPLIT' and SPLIT'~$k$, $k=1,2,3$, denote the versions that were valid before this
paragraph.
Note that after the application of SPLIT~2 the following holds: If we were in Case~1 
and $T_{i''+1,i''+1}$ is a horizontal candidate then $T_{i'',i''}$ is a stable vertical 
candidate, if we were in Case~2 and $T_{i''+1,i''+1}$ and $T_{i''+2,i''+2}$ are horizontal 
candidates then $T_{i''+1,i''+1}$ is stable.

\medskip

In Section~\ref{alg}, before the description of SPLIT~3, we define a set $C$, which we need
as an input for SPLIT~3. Note that we can omit the candidates in $C$ that
neither touch nor are touched from above in the application of SPLIT~3 and consequently these entries are
omitted in the following.

\medskip

\begin{center} \sc 6.2. The main lemmas I     \end{center}

\medskip

The lemmas in this subsection are needed in order to prove that SPLIT~1 and MERGE~1 
are inverse to each other. They have analogs that are used for 
showing that SPLIT~3 and MERGE~3 are each other's respective inverses as we see in 
Subsection~6.4.

\medskip

We introduce some notation concerning the relative position of 
an entry and a path. Let $e$ and $e'$ be two 
entries in a shifted tabloid $T$ in the $i$-th row
and $c_T(e)= (i,j)$, $c_T(e')=(i,k)$. If $j<k$, 
resp. $j \le k$,  we say that $e'$ is right, resp. weakly right, of $e$.
If $P$ is a path in the shifted tabloid $T$ then $e'$ is said to be 
right, respectively weakly right, of $P$, if $P$ includes an entry $e$ such that $e'$ is 
right, respectively weakly right,  of $e$. Similar definitions apply to left, weakly left, above, weakly above, 
below and  weakly below.  If the formulations of our lemmas and corollaries 
include phrases in brackets [$\dots$], the assertions are true with and without these phrases.

\medskip

\begin{figure}
\setlength{\unitlength}{1cm}
\begin{center}
\scalebox{0.40}{\includegraphics{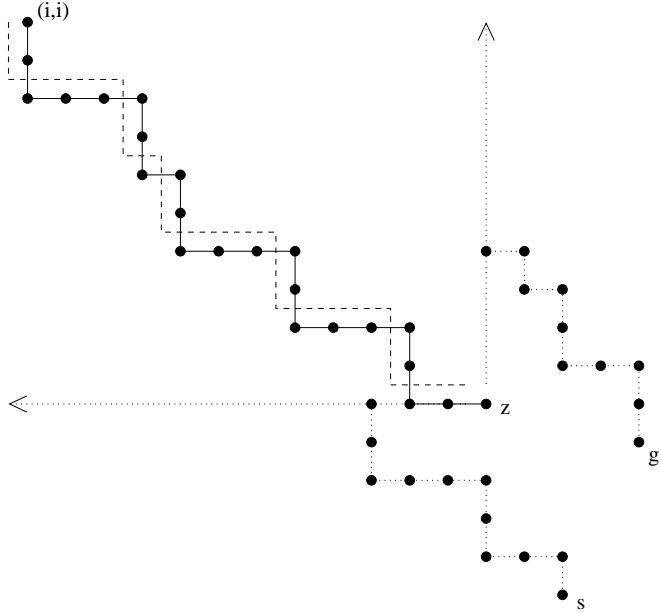}}
\end{center}
\caption{The full line indicates the backward path of entry $z$, the circles indicate the 
centers of the traversed cells. The dashed line indicates the border: The entries greater than 
$z$ with respect to the backward path order are located north-east of this border and 
the entries smaller than $z$ are located south-west of this border. Moreover observe that 
$g$ is greater than $z$ for its backward path (dotted line) enters the column of $z$ weakly 
above of $z$ and $s$ is smaller than $z$ for its backward path enters the row of $z$
weakly left of $z$.}
\label{l1}   
\end{figure}

In the following lemma we characterize the largest and the smallest entry with 
respect to the backward paths order of a given set of entries. See Figure~\ref{l1}.

\smallskip

\begin{lem}
\label{0}
Let $Z$ be a set of entries in a shifted tabloid $T$. 
\begin{enumerate}
\item Let $s \in Z$ and $P_s$ denote the backward path of $s$. Then 
$s$ is the smallest entry in $Z$ with respect to the backward paths order, if 
and only if every $z \in Z$ is either [weakly right and] above of $P_s$ or the 
backward path of $z$ enters the column of $s$ weakly above of $s$.

\item Let $g \in Z$ and $P_g$ denote the backward path of $g$. 
Then $g$ is the greatest entry in $Z$ with respect to the 
backward paths order, if and only if every $z \in Z$ is either left [and weakly below] of 
$P_g$ or the backward path
of $z$ enters the row of $g$ weakly left of $g$.
\end{enumerate}
\end{lem}

\medskip

\medskip

We extend the definition of `entry $e$ touches entry $e'$ from above': Let $T$ be 
a shifted tabloid, $e$, $e'$ two entries in $T$ and $D$ a set of cells 
in the associated shifted Ferrers diagram. Let $P$ be the forward path of $e$ with respect to $D$
and $Q$ the backward path of $e'$ in $T$. If there exists a $k$
such that $(k-1,k-1), (k-1,k) \in P$ and $(k,k) \in Q$, we say that
{\it  $e$ touches $e'$ from above with respect to $D$}. If we want to refer to the 
restricted definition for touching from above which was valid before this paragraph 
we add the phrase `in the application of SPLIT'. Note that there exist pairs
$(T,H)$ in which an entry $e$ touches another entry $e'$ from above in the
sense we just defined, but there is no step in the application of SPLIT~3 to
$(T,H)$  where $e$ touches $e'$ from above in the original sense.

\smallskip

\begin{figure}
\setlength{\unitlength}{1cm}
\begin{center}
\scalebox{0.40}{\includegraphics{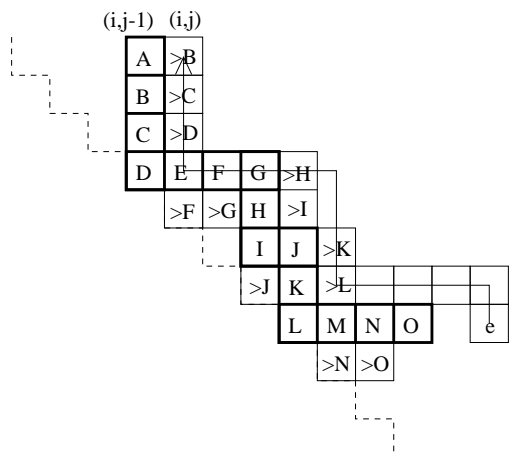}}
\end{center}
\caption{$P_{i,j-1}=\{A,B,C,\dots,M,N,O\}$. Observe that the entry to the right of `$A$' is greater than the entry in the 
cell labeled with `$B$' in $T$ (therefore labelled with `$> B$') for $P_{i,j-1}$ starts with a south step. The corresponding 
assertion holds for the other labels because of the run of $P_{i,j-1}$. Furthermore a possible backward path $Q$
of $e$ is drawn. Consider for example the cell of $Q$ labelled with `$E$'. 
Then, by reverse jeu de taquin, 
the next cell of $Q$ is forced to be the cell labelled with `$> D$' and not the cell labelled with `$D$'.}
\label{l2}   
\end{figure}

The following lemma is fundamental. It is illustrated in Figure~\ref{l2}.

\medskip

\begin{figure}
\setlength{\unitlength}{1cm}
\begin{center}
\scalebox{0.40}{\includegraphics{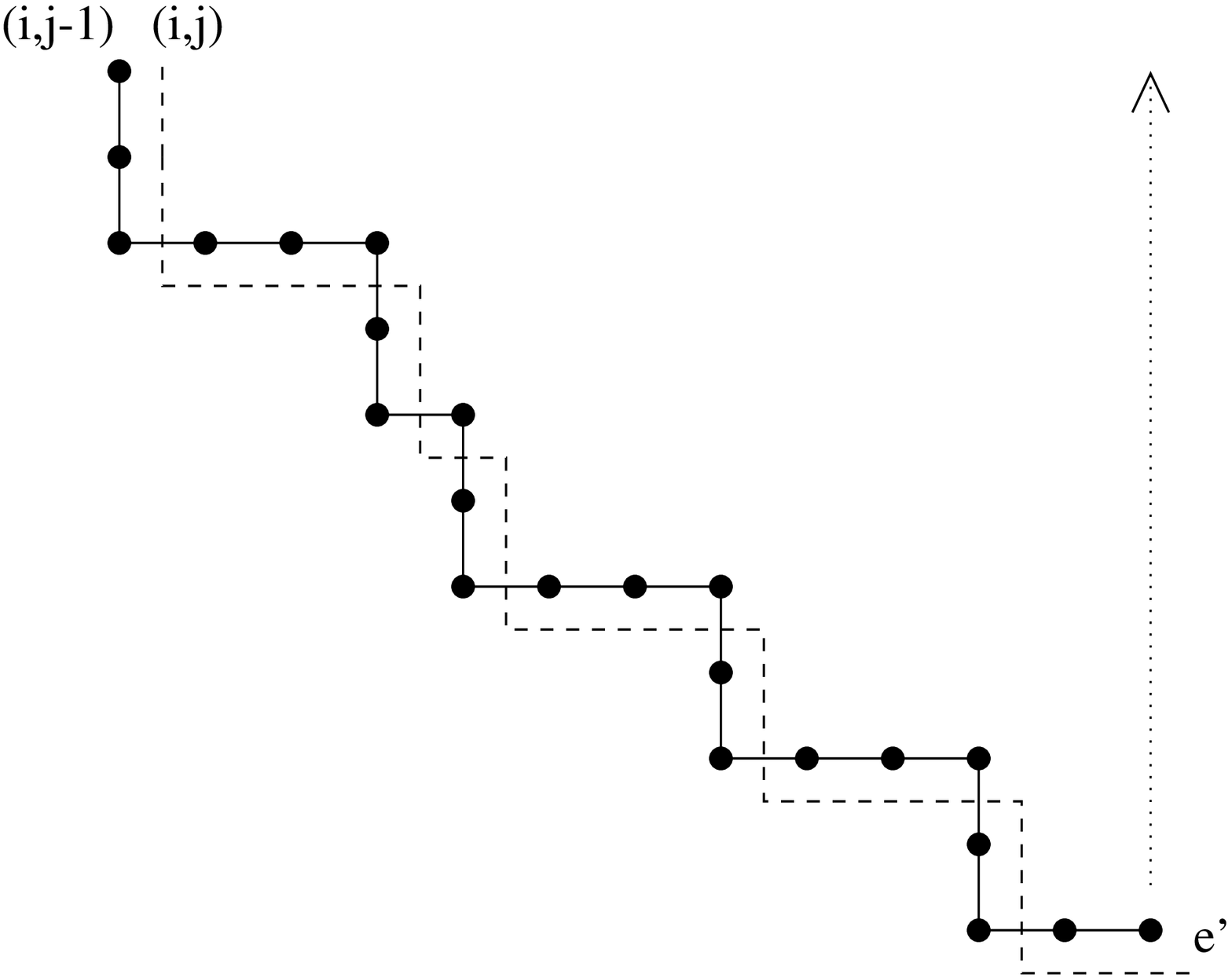}}
\end{center}
\caption{The full line is the forward path of $T_{i,j-1}$ in $T$, the circles indicate the center of the
cells the forward path traverses. The dashed line is the border, possible cells for $e$ are located north-east of this border.}
\label{el2}   
\end{figure}

\begin{lem}
\label{1}
Let $(i,j-1)$ and $(i,j)$ be two cells in a shifted tabloid $T$, which 
is ordered up to $(i,j)$. Let $P_{i,j-1}$ denote the forward path of 
$T_{i,j-1}$ in $T$ with respect to the a set $D$ and $e'$ the entry of the endcell of $P_{i,j-1}$ in $T$. 
Let $e$ be an entry in $T$ whose backward path ($i$ being the fixed row) contains $(i,j)$. 
If $T_{i,j-1}$ does not touch $e$ from above with respect to $D$, 
then $e$ is right [and weakly above] of $P_{i,j-1}$   or the 
backward path of $e$  
enters the column of  $e'$ weakly above of $e'$. (See Figure~\ref{el2}.)
\end{lem}

{\it Proof.} 
If $P_{i,j-1}$ consists solely of south steps there 
is nothing to prove. Thus we assume that there is at least one east step in 
$P_{i,j-1}$.

We show the following: Let $e''$ be below of $P_{i,j-1}$.
If $T_{i,j-1}$ does not touch $e''$ from above 
with respect to $D$ 
then every cell in the backward path $P$ of $e''$ ($i$ being the fixed row) 
in a column greater than $j-1$ is below of $P_{i,j-1}$.

If the statment were false there would exist integers $g$,$l$ such that 
$(g,l), (g,l+1) \in P_{i,j-1}$ and $(g,l+1), (g+1,l+1) \in P$. 
Since $T_{i,j-1}$ does not touch $e''$ from above with respect to $D$ the cell $(g,l)$ is
not a cell in the main diagonal and thus $(g+1,l)$ is a cell in the shifted Ferrers 
diagram. 
From $(g,l), (g,l+1) \in P_{i,j-1}$ 
it follows that $T_{g,l+1} <  T_{g+1,l}$, which is a contradiction 
to $T_{g,l+1} > T_{g+1,l}$, which follows from  $(g,l+1), (g+1,l+1) \in P$.
(Later we will often refer to this important argument as 'the argument in the proof of Lemma~\ref{1}'.)

If the statement in the lemma were false, $e$ is below 
of $P_{i,j-1}$ or the backward path of $e$ 
enters the column of $e'$
below  of $e'$. By the assertion we have just proved every entry 
in the backward path $P$ of $e$ which is in a column greater than $j-1$ 
and weakly left of the column of $e'$ 
must be below of $P_{i,j-1}$.
But, since $(i,j)$  is contained in the backward path of $e$  and $(i,j)$ 
is  not below of 
$P_{i,j-1}$, this is a contradiction if $P_{i,j-1}$ does not consist solely 
of south steps.
\qed

\smallskip

Note that if we are in the situation of Lemma~\ref{1} and $e$ enters the column of $e'$ weakly above of $e'$ then 
$e$ enters the column of $e'$ weakly above uppermost cell in $P_{i,j-1}$ and 
the column of $e'$  by the proof of the lemma.

\smallskip

\smallskip

\begin{rem}
\label{5}
Mostly we apply the lemma in the following situation.
Let $(i,j-1)$ and  $(i,j)$  be two cells in a shifted tabloid $T$ which is ordered up to 
the predecessor of $(i,j)$ in the total order. Furthermore let $P_{i,j}$ denote 
the forward path of $T_{i,j}$ in $T$ with respect to a set $D$ and let $P_{i,j-1}$ denote
the forward path of $T_{i,j-1}$ in $T'$ with respect to a set $D'$, where $T'$  denotes the shifted tabloid we 
obtain after performing jeu de taquin in $T$ with $T_{i,j}$ and with respect to $D$. Let $e'$ denote the entry in the last cell 
of $P_{i,j-1}$ in $T'$. Observe that $P_{i,j}$ 
coincides with the backward path of $T_{i,j}$ in $T'$ ($i$ being the fixed row) 
except for the part of the $i$-th row left of $(i,j)$. 
Thus, if $T_{i,j-1}$ does not touch $T_{i,j}$ from above with respect to $D'$ in $T'$, 
then $T_{i,j}$ is right [and weakly above] of 
$P_{i,j-1}$ in $T'$ or the backward path of $T_{i,j}$ in $T'$ enters the column of $e'$ weakly above of 
$e'$ by Lemma~\ref{1}.
\end{rem}


\medskip

Observe that Figure~\ref{l1} and Figure~\ref{el2} coincide after a shift of the (dashed) domain border
in Figure~\ref{l1} by the vector $(1,-1)$. This leads us to the following.

\smallskip

\begin{cor}
\label{3}
Let $(i,j-1)$, $(i,j)$, $T$ and $T'$ be as in Remark~\ref{5}.
Furthermore let 
$Z'$ be a set of cells in the last
$r-i+1$ rows of $T'$,  such that $T_{i,j}$ is the smallest 
entry under the entries in the cells of $Z'$ with respect to the backward paths order in $T'$. 
Let $T''$ denote the 
shifted tabloid we obtain after performing jeu de taquin with $T_{i,j-1}$ in $T'$ 
and with respect to $D'$
and let $j'$ denote the column of $T_{i,j-1}$ in $T''$.
Let $Z''$ denote the set of cells we obtain from $Z'$ by 
first replacing every cell $(h,k) \in Z'$ with $k \le j'$ by  $(h-1,k-1)$ and 
then deleting the replacing cells in the $(i-1)$-st row.
If $T_{i,j-1}$ touches no entry in a cell in  $Z'$ from above with respect to $D'$ in  
$T'$ then 
$T_{i,j-1}$ is smaller than every entry in a cell in $Z''$ with respect to 
the backward paths order in $T''$.
\end{cor}

{\it Proof.} By Lemma~\ref{0} (1), Lemma~\ref{1} and Remark~\ref{5}  it is clear that either 
the entry in 
the cell of $Z''$ that came from $T_{i,j}$ is greater than $T_{i,j-1}$ with respect to 
the backward paths order in $T''$ or the cell that came from $T_{i,j}$ was
deleted in the course of constructing $Z''$.
(This is due to the following observation: Let $\rho$ be a cell which is right and weakly above of a path $P$ then 
$\rho-(1,1)$ is weakly right and above of $P$ or $P$ has no cell in the row 
of $\rho-(1,1)$.)

Now suppose that $e$ is an entry in $T'$ weakly below the $i$-th row that it is greater than 
$T_{i,j}$ with respect to the backward paths order in $T'$. 
Furthermore suppose that $T_{i,j-1}$ does not touch $e$ from above with
respect to $D'$ in $T'$. Thus, by the
relative position of $e$ and $T_{i,j}$ by Lemma~\ref{0}, and  since 
the backward path ($i$ being the fixed row) of $T_{i,j}$ in $T'$ includes $(i,j)$, the backward path of $e$ 
in $T'$ includes $(i,j)$ as well. (Observe that the backward path of $e$ coincides
with the backward path of $T_{i,j}$ after they intersect.)
By Lemma~\ref{1} this implies that $e$ is either 
right [and weakly above] of $P_{i,j-1}$ in $T'$ or the backward path of 
$e$ enters the column of the end $e'$ of $P_{i,j-1}$ in $T'$ weakly 
above of $e'$. From this the assertion follows by Lemma~\ref{0} (1).
\qed

\medskip

In Section~\ref{alg} after SPLIT'~1 we claim that for an entry $e$ whose forward 
path terminates
in a cell $(k,k)$ in the course of SPLIT'~1, we have $T_{k,k}=e$ and $H_{i,k}=(k,k)$
after the application of SPLIT'~1.  This is also true for 
SPLIT~1 and proved in the subsequent lemma.

\begin{lem}
\label{6}
Let $T$ be a shifted tabloid which is ordered up to $(i+1,i+1)$ and $H$ a partial shifted 
hook tabloid where the $i$-th row is empty. Let $(T^j,H^j)$ denote the pair we obtain in the course of applying 
SPLIT~1 to $(T,H)$ in the $i$-th row after performing jeu de taquin with $T_{i,j}$. Assume that the 
forward path of $T_{i,j}$ in $T^{j+1}$ contains cells on the main diagonal and 
let $(k,k)$ be the first of these in the run of the forward path. Let $h<j$ be such that the 
forward path $P$ of $T_{i,h}$ in $T^{h+1}$ ends weakly right of 
the $k$-th column. Furthermore suppose that the forward path of $T_{i,l}$ in $T^{l+1}$  
does not contain a cell on the main diagonal for 
$h < l < j$.  Then $T_{i,h}$ touches $T_{i,j}$ from above in $T^{h+1}$, to be more 
accurate $(k-1,k-1), (k-1,k) \in P$ and $(k,k)$ is the cell of $T_{i,j}$ in $T^{h+1}$.
\end{lem}

{\it Proof.} 
We show that the backward path of $T_{i,j}$ in $T^{h+1}$ does not 
contain a cell $(g,g)$ on the main diagonal with $i < g < k$. Furthermore this 
backward path contains $(i,h+1)$ and $(k,k)$. The assertion is then a consequence of 
Lemma~\ref{1}: Let $T$ in Lemma~\ref{1} equal $T^{h+1}$, $j$ equal $h+1$
and let $D=\{(g,k)\}$, where $(g,k)$ is the first cell in the $k$-th column
in the forward path $P$ of $T_{i,h}$ in $T^{h+1}$. By Lemma~\ref{1} $T_{i,j}$ is
either right of $P$ or the backward path of $T_{i,j}$ enters column $k$
weakly above of $(g,k)$ or $T_{i,h}$ touches $T_{i,j}$ from above in $T^{h+1}$. Since
$g < k$ and $T_{i,j}$ is either strictly left of the $k$-th column in
$T^{h+1}$ or in $(k,k)$ the first two options are impossible and thus
$T_{i,h}$ touches $T_{i,j}$ from above. 
The only possible place for touching
is $(k-1,k-1)$, $(k-1,k)$, $(k,k)$ for the backward path of $T_{i,j}$ in
$T^{h+1}$ does not contain a cell on the main diagonal strictly between the
$i$-th and the $k$-th column.

Observe that the backward path ($i$ being the fixed row) of 
$T_{i,j}$ in $T^j$ contains no 
cell $(g,g)$ with $i < g < k$ and it contains $(i,j)$ and therefore 
$(i,h+1)$, for it coincides with the forward path of $T_{i,j}$ in $T^{j+1}$
except for the part strictly left of $(i,j)$. If $j-1=h$ there is nothing to prove. Otherwise the forward path 
of $T_{i,j-1}$ in $T^j$ contains no cell on the main diagonal and 
thus the backward path of $T_{i,j}$ in $T^{j-1}$ does not contain a
cell on the main diagonal strictly between column $i$ and $k$ and it contains 
$(i,j-1)$, for this backward path coincides with the backward path of 
$T_{i,j}$ in $T^j$ until it meets the backward path of 
$T_{i,j-1}$ in $T^{j-1}$ and after that it coincides with the backward path of 
$T_{i,j-1}$ in $T^{j-1}$. The claim at the beginning of the proof now follows by induction with respect to $j-h$.
\qed   

\medskip

Now we are ready to show the assertion before the formulation of Lemma~\ref{6}
(We use the notation of Lemma~\ref{6}.): Assume that
the forward path of $T_{i,k}$ in $T^{k+1}$ contains a cell $\rho$ on the main diagonal and 
let $h < k$ be such that the forward path of $T_{i,h}$ in $T^{h+1}$ ends weakly 
right of the column of $\rho$. Let $D$ be the appropriate set of cells in SPLIT~ 1 such that 
$T^h=\jt_D(T^{h+1},T_{i,h})$, then it suffices to show that the 
forward path of $T_{i,h}$ in $T^{h+1}$ with respect to
$D$ ends strictly left of the column of $\rho$. We show the assertion by induction 
with respect to $k-h$.
To this end let $j$, $j > h$, be minimal such that the 
forward
path of $T_{i,j}$ in $T^{j+1}$ with respect to the appropriate set in SPLIT~1 contains a cell 
$\rho'$ on the main diagonal. By 
Lemma~\ref{6} the forward path of $T_{i,h}$ in $T^{h+1}$ with respect to $D$ 
ends in a column strictly left of 
the column of $\rho'$. If $k=j$ the assertion follows immediately.  Otherwise the column of $\rho'$ is smaller than the column of $\rho$ by 
induction since $k-j < k - h$ and the assertion follows too.

\smallskip

Actually the corollary also implies the following inversion 
of the statement above: If $H_{i,k}=(k,k)$ after the application of SPLIT~1, 
then $T_{k,k}$ is an entry with which we have performed jeu de taquin in the course of 
applying SPLIT~1 and $(k,k)$ is the end of its forward path.

\smallskip 

In fact this is a consequence of the following observation which can be shown
by Lemma~\ref{1} and Lemma~\ref{6}: Assume that the backward path of  a horizontal candidate in 
$(T^j,H^j)$ contains $(k,k)$. Then there exists an $l$, $j \le l \le \lambda_i + i - 1$,
such that the forward path of $T_{i,l}$ in $T^{l+1}$ with respect to the appropriate $D$ in SPLIT~1
contains $(k,k)$ and
$T_{i,l}$ is the smallest horizontal candidate in $(T^j,H^j)$, whose 
backward path contains $(k,k)$.
\medskip

\begin{center} \sc 6.3. SPLIT~1 and MERGE~1 are each other's respective inverse \end{center}

\medskip

We define property AFTER\_SPLIT~1 for a pair $(T,H)$ and  will observe that 
the application of SPLIT~1 yields a pair with that property. 
A pair $(T,H)$ of a shifted tabloid  $T$ and a partial shifted hook tabloid $H$ has 
property AFTER\_SPLIT~1 if 

\begin{enumerate}
\item the last $r-i+1$ rows of $H$ form a shifted hook tabloid,

\item there exists no vertical candidate (with respect to the $i$-th row) in $(T,H)$,

\item the subtabloid of $T$ consisting of the last $r-i+1$ rows is standard, except for 
$T_{k,k} \in C$ there might hold $T_{k,k} \ge T_{k,k+1}$, where $C$ is the set of entries 
$T_{k,k}$ with $H_{i,k}=(k,k)$. 

\item for $T_{k,k} \in C$  such that $T_{k,k}$ either proves to be exceptional 
or touched from above in the application of SPLIT 
and every diagonal cell  $(j,j) \not= (k,k)$ in the backward path of $T_{k,k}$, there exists 
a horizontal candidate $e$ in $T$ whose 
backward path includes $(j,j)$ and $e <_{T} T_{k,k}$,   and 

\item $T_{k,k} \in C$ is unstable and neither exceptional nor touched from 
above in the application of SPLIT, then $T_{k,k}$  touches another entry 
from above in the 
application of SPLIT and is in its place of change.

\end{enumerate}

\medskip

{\bf Claim~1.}
\begin{enumerate}
\item Let $T$ be a shifted tabloid which is ordered up to $(i+1,i+1)$ and 
$H$ a partial shifted hook tabloid such that the $i$-th row is empty. Let $(T',H')$ denote the 
pair we obtain after the application of SPLIT~1. Then $(T',H')$ has property AFTER\_SPLIT~1
and the application of MERGE~1 to $(T',H')$ yields $(T,H)$.
\item Let $(T',H')$ be a pair with property AFTER\_SPLIT~1 and $(T,H)$ denote the 
pair we obtain after the application of MERGE~1. Then $T$ is ordered up to $(i+1,i+1)$ and 
$H$ is a partial shifted hook tabloid such that the $i$-th row is empty and the application of 
SPLIT~1 to $(T,H)$ yields $(T',H')$. 
\end{enumerate}

\smallskip

{\it Proof.}

{\it re 1.} 
Let $(T^j,H^j)$ denote the pair we obtain after the performance of
jeu de taquin with $T_{i,j}$ and the corresponding shift 
in the partial shifted hook tabloid in SPLIT~1 and 
set $(T^{\lambda_{i}+i}, H^{\lambda_{i}+i})=(T,H)$. 
By construction the current cell of $T_{i,j+1}$ in $T^{j+1}$ is in 
$\{H^{j+1}_{i,j+1},H^{j+1}_{i,j+2},\dots,H^{j+1}_{i,\lambda_i+i-1} \}$.
Furthermore by induction we assume that $T_{i,j+1}$ 
is the smallest horizontal candidate in a cell in 
$\{H^{j+1}_{i,j+1},H^{j+1}_{i,j+2},\dots,H^{j+1}_{i,\lambda_i+i-1}\}$ with respect to 
the backward paths order in $T^{j+1}$.
Let $D$ be the set such that $T^j$ is obtained 
from $T^{j+1}$ in SPLIT~1 by performing jeu de taquin with $T_{i,j}$ and with respect to 
$D$. By construction  no horizontal candidate in a cell in 
$\{H^{j+1}_{i,j+1},H^{j+1}_{i,j+2},\dots,H^{j+1}_{i,\lambda_i+i-1}\}$ 
is touched from above by $T_{i,j}$ with respect to $D$ in $T^{j+1}$:
Let $k$ be minimal such that the backward path in $T^{j+1}$ of a horizontal candidate with respect to 
$H^{j+1}$ contains $(k,k)$. Then the forward path of $T_{i,j}$
with respect to $D$ ends strictly left of the column of $k$ by Lemma~\ref{6}. 
(See also the last paragraph in the previous subsection.)
Thus, by Corollary~\ref{3}, $T_{i,j}$ is smaller or equal than
every horizontal candidate in  a cell in $\{H^{j}_{i,j},H^{j}_{i,j+1},
\dots,H^j_{i,\lambda_i+i-1}\}$
with respect to the backward paths order in $T^j$. Therefore
if we apply the $j-i+1$-st step of MERGE~1 to 
the pair $(T^j,H^j)$ 
we reobtain $(T^{j+1},H^{j+1})$.

\smallskip

{\it re 2.}
We define the following algorithm PRE-MERGE~1. The input 
is a pair $(T',H')$ with the property AFTER\_SPLIT~1.

\smallskip

\fbox{ \parbox{14cm}{
{\bf PRE-MERGE~1.} Repeat the following: [ Let $(j,j)$ be such that $H'_{i,j}\not=(j,j)$
but the backward path of a horizontal candidate in $T'$
contains $(j,j)$. If such a $j$ does not exist stop. Otherwise let $e$ be the
smallest horizontal candidate with that property 
and set $(T',H')=\rjs_{(j,j)}(T',H',e)$.]
}}

\smallskip

Observe that the application of PRE-MERGE~1 and MERGE~1 to $(T',H')$ yields the same result 
as the application of MERGE~1 to $(T',H')$ for the smallest horizontal candidate 
whose backward path contains $(j,j)$ is the smallest horizontal candidate with
respect to the backward paths order weakly right of the $j$-th column. Let $(T'',H'')$ denote 
the pair we obtain after the application of PRE-MERGE~1 to $(T',H')$. Then 
$(T'',H'')$ has the properties 1 -- 3 from AFTER\_SPLIT  
and for every cell $(k,k)$ on the main diagonal such that
the backward path of a horizontal candidate contains $(k,k)$,
$T''_{k,k}$ is itself a horizontal candidate.

We denote the property of $(T'',H'')$  by AFTER\_SPLIT'~1.
(If we apply SPLIT'~1 to a pair 
$(T,H)$, $T$ ordered up to $(i+1,i+1)$ and the $i$-th row of $H$ is empty, then 
the output pair has property AFTER\_SPLIT'~1.)

\smallskip

We show the following: Let $(T,H)$ denote the pair we obtain 
after the application of MERGE~1 to $(T'',H'')$. If we apply 
SPLIT'~1 to $(T,H)$ we reobtain $(T'',H'')$.

\smallskip

Let $({T''}^j,{H''}^j)$ denote the pair we obtain in the course of applying MERGE~1
to $(T'',H'')$ after performing reverse 
jeu de taquin with respect to $\{(i,j-1)\}$ and the 
corresponding reshift in the partial shifted hook tabloid, and 
set $({T''}^i,{H''}^i)=(T'',H'')$. 
By induction we assume that the backward path of every horizontal candidate in a cell in
$\{{H''}^j_{i,j}, {H''}^j_{i,j+1},\dots,{H''}^j_{i,\lambda_i+i-1}\}$ in 
${T''}^j$ contains $(i,j)$ and it contains no cell on the main diagonal except for
possibly $c_{{T''}^{j}}$ and $(i,j)$. Let $e$ be the smallest horizontal candidate in a cell in
$\{{H''}^{j}_{i,j},{H''}^{j}_{i,j+1},\dots,{H''}^j_{i,\lambda_i+i-1}\}$ 
with respect to the backward paths order in ${T''}^j$.
Furthermore let $P$ be its  backward paths in ${T''}^j$ with respect to $(i,j)$.
If $j'$ is the column of $e$ in  ${T''}^j$, the subtabloid
of ${T''}^j$ consisting of the first $j'$ 
columns is ordered up to $(i,j)$ by property AFTER\_SPLIT'~1 and the
minimality of $e$.
Thus $e={T''}^{j+1}_{i,j}$ and performing jeu de 
taquin with $e$ in ${T''}^{j+1}$ and with respect to
the cells on the main diagonal together with the corresponding 
shift in ${H''}^{j+1}$ results in the pair $({T''}^j,{H''}^j)$. 
By the minimality of $e$, Lemma~\ref{0} (1) and by the construction of the reshift in 
${H''}^j$ the entries in the cells in 
$\{{H''}^{j+1}_{i,j+1},{H''}^{j+1}_{i,j+2},\dots,{H''}^{j+1}_{i,\lambda_i+i-1}\}$
are right [and weakly above] of $P$ or their backward paths in ${T''}^{j+1}$ enters column $j'$ 
weakly above of $c_{{T''}^{j}}(e)$. Thus, by the argument from 
Lemma~\ref{1}, the backward paths of the entries in these cells
in ${T''}^{j+1}$ contain $(i,j+1)$. Likewise it is easy to see that no
backward path of these entries contains a cell on the main diagonal
except for possibly its origin and terminus. 

\smallskip

We define another Algorithm POST-SPLIT~1 which we apply to 
a pair $(T'',H'')$ with property AFTER\_SPLIT'~1.

\smallskip

\fbox{ \parbox{14cm}{
{\bf POST-SPLIT~1.} 
Repeat the following:
[ Let $T''_{j,j}$ be an unstable horizontal candidate which is neither exceptional 
nor touched from above and not already in its place of change (if
$T''_{j,j}$ has a place of change). If  
$T''_{j,j}$ does not exist stop. Otherwise: If $T''_{j,j}$ has 
a place of change  $\rho$ let $D=\{ \rho\}$ otherwise let $D = \emptyset$.
Set $(T'',H'')=\js_D(T'',H'',T''_{j,j})$.] }}

\smallskip

Now the assertion follows since PRE-MERGE~1 and POST-SPLIT~1 are 
inverse to each other (in order to show that use Lemma~\ref{6}) and since the application of SPLIT'~1 
and POST-SPLIT~1 is equivalent to the application of SPLIT~1.
\qed

\medskip

In view of Claim~1 it remains to show the following:

\begin{enumerate}
\item Let $(T,H)$ be a pair with the property AFTER\_SPLIT~1 and
let $(T',H')$ denote the pair we obtain after the application of SPLIT~2 and
SPLIT~3 to $(T,H)$. We have to show that if we apply MERGE~3 and MERGE~2 to $(T',H')$ we reobtain $(T,H)$.

\item Let $(T,H)$ be such that $T$ is ordered up to $(i,i)$ and 
the subtabloid of $H$ consisting of the last $r-i+1$ rows is a shifted hook
tabloid. Let $(T',H')$ denote the pair 
we obtain after the application of MERGE~3 and MERGE~2. We have to 
show that $(T',H')$ has property AFTER\_SPLIT~1 and that the application 
of SPLIT~2 and SPLIT~3 to $(T',H')$ yields $(T,H)$.
\end{enumerate}

\medskip

\begin{center}  \sc 6.4. The main lemmas II \end{center}

\medskip

The following lemma is similar to Lemma~\ref{1}.

\smallskip

\begin{lem}
\label{2}
Let $(i,m)$ and $(i-1,n)$, $m \le n$, be two cells in a shifted tabloid $T$
such that the subtabloid of $T$ consisting of the last $r-i+2$ rows and without the entries in $(i,m)$
and $(i-1,n)$ has increasing rows and columns. 
Let $P_{i-1,n}$ denote the forward path of 
$T_{i-1,n}$ in $T'$ with respect to a set, where $T'$ denotes the shifted tabloid we obtain 
after performing jeu de taquin in $T$ with $T_{i,m}$ and with respect to $D$. Furthermore let $e$ 
denote the entry in the last cell of $P_{i-1,n}$ in $T'$. Then $T_{i,m}$ is [weakly left and]
below of $P'_{i-1,n}=P_{i-1,n} \cup \{(i-1,m),(i-1,m+1),\dots, (i-1,n-1)\}$ in 
$T'$ or the backward path of $T_{i,m}$ with respect to $(i,m)$ in $T'$ enters 
the row of $e$ weakly left of $e$. 
\end{lem}

{\it Proof.} Similar to the proof of Lemma~\ref{1}.  Observe that if the backward path 
of $T_{i,m}$ in $T'$ enters the row of $e$ weakly left of $e$ then it enters the row of 
$e$ weakly left of the leftmost cell of $P_{i-1,n}$ in the row of $e$. \qed

\smallskip


\begin{cor}
\label{4}
Let $(i-1,n)$, $(i,m)$, $T$, $T'$ be as in Lemma~\ref{2}. Furthermore let 
$Z'$ be a set of cells weakly below the $i$-th row, which includes the cell of $T_{i,m}$ in $T'$,  
such that $T_{i,m}$ is the greatest 
entry in a cell of $Z'$ in $T'$ with respect to the backward paths order.
Let $T''$ denote the 
shifted tabloid we obtain after performing jeu de taquin with $T_{i-1,n}$ in $T'$ and 
with respect to a set 
and let $i'$ denote the row of $T_{i-1,n}$ in $T''$.
Let $Z''$ denote the set of cells we obtain from $Z'$ by replacing every cell $(h,k)$ in $Z'$ with $h \le i'$ by  $(h-1,k-1)$.
Then $T_{i-1,n}$ is greater than every entry in a cell of $Z''$ in $T''$ with respect to 
the backward paths order.
\end{cor}

{\it Proof.} Similar to the proof of Corollary~\ref{3}. \qed

\medskip

\begin{center} \sc 6.5. SPLIT~2 and MERGE~2 are each other's respective inverse \end{center}

\medskip

We define property AFTER\_SPLIT~2 for a pair $(T,H)$ and will observe that 
a pair to which we have applied SPLIT~1 and SPLIT~2 has that property.
A pair $(T,H)$ of a shifted tabloid $T$ and a partial shifted hook tabloid  $H$ 
has property AFTER\_SPLIT~2  if 

\begin{enumerate}

\item the last $r-i+1$ rows of $H$ form a shifted hook tabloid,

\item there exists an $i''$ such that for every $i+1 \le g \le i''$ there is a vertical 
candidate in row $g$ and no other vertical candidates exist. 
For $i+1 \le g < g+1 \le i''$ the vertical candidate in row 
$g$ is greater than the vertical candidate in row $g+1$ with 
respect to the backward paths order ($i$ being the fixed row). If we perform 
jeu de taquin with the vertical  candidates in row $i'', i''-1, \dots ,i+1$ in
that order,  then every
candidate is either stable or its forward path starts with a step 
into the next row. 
Moreover: $T_{i'',i''}$ is a stable vertical candidate or 
$T_{i''+,1,i''+1}$ is a stable horizontal candidate or $T_{i''+1,i''+1}$ is 
not a horizontal candidate or $T_{i''+2,i''+2}$ is not a horizontal candidate,

\item  the subtabolid of $T$ consisting of the last $r-i+1$ rows 
is standard, except for $T_{p,q} \in C' \cup C$ there might hold 
$T_{p,q} > \min (T_{p+1,q},T_{p,q+1})$, where $C'$ and $C$ 
are the two sets we define before the application of 
SPLIT~3 to $(T,H)$ (see Section~\ref{alg} before the description of SPLIT~3), 

\item  $e$ is a horizontal candidate which is touched from 
above in the application of SPLIT, then 
for every cell $(j,j)\not= c_T(e)$ with $i''<j$
in the backward path of $e$, there exists a horizontal 
candidate $e'$ with $e' <_T e$ and the backward path of $e'$ contains $(j,j)$,  and 

\item $e$ is an unstable horizontal candidate which is not touched from above in the 
application of SPLIT then 
$e$ touches from above in the application of SPLIT and is in its place of change. 

\end{enumerate}

\smallskip

A remark on Property~(2) of AFTER\_SPLIT~2:
Note that the vertical candidate in row $g$ is greater than the vertical
candidate in row $g+1$, $i+1 \le g < g+1 \le i''$ in a pair to which we have 
applied SPLIT~2 by Lemma~\ref{2}.

\medskip

{\bf Claim 2.} 
\begin{enumerate}
\item Let $(T,H)$ be a pair with property AFTER\_SPLIT~1 and 
let $(T',H')$ denote the pair we obtain after the application of 
SPLIT~2. Then $(T',H')$ has property AFTER\_SPLIT~2 and the 
application of MERGE~2 to $(T',H')$ results in $(T,H)$. 
\item Let $(T',H')$ be a pair with property AFTER\_SPLIT~2 
and $(T,H)$ denote the pair we obtain after the application
of MERGE~2. Then $(T,H)$ has property AFTER\_SPLIT~1 and the 
application of SPLIT~2 to $(T,H)$ results in $(T',H')$.
\end{enumerate}
  
{\it Proof.} Left to the reader. 
Observe that we are in Case~2 of SPLIT~2, respectively Case~2 of
MERGE~2, if and only if either the smallest horizontal candidate with respect to
the backward paths order in $(T',H')$ in a row strictly below the lowest vertical candidate
is smaller than the smallest vertical candidate or there exists no vertical candidate.  The rest follows from the
fact that jeu de taquin and reverse jeu de taquin are inverse to each other
and from the fact that (with one exception) the endcells of the forward
paths are not distorted by a shift in the accompanying shifted hook
tabloid. The exception is the exceptional vertical candidate in the $h$-th row in
Case~1 of SPLIT~2, whose endcell is stored in the $i$-th row. There we need
Corollary~\ref{4}. \qed

\medskip

\begin{center} \sc 6.6. Irreducible pairs \end{center}      

\medskip

{\bf Touching from below.}
We define an analog to `$e$ touches $e'$ from above' for the 
algorithm  MERGE.
Suppose we are at the beginning of a step of MERGE~3 and the if-condition (If $c_T(e)=(h,h)$
: [$\dots$]) is true. In the following commands in the description of MERGE~3
we distinguish between two cases.
If we are in the first case we say that $e'$ touches 
$e$ from below and for $h' < j \le h$, $T_{j,j}$ touches 
$T_{j-1,j-1}$ from below. Furthermore for $h' \le j \le h$, let 
$(j,j)$ be the place of change for $T_{j,j}$ .
In the second case we say that  for $h'  <  j \le h$, $T_{j,j}$ touches $T_{j-1,j-1}$
from below and $(j,j)$ is the place of change for $T_{j,j}$.

We define the exceptional entries with respect to MERGE~2. Suppose we have applied MERGE~2, 
Case~1. Define $T_{g,g}$ to be exceptional for $i \le g \le i''$. Now suppose we 
have applied MERGE~2, Case~2. For $i \le g \le i''$ and $g \not=i',i'+1$ we define $T_{g,g}$
to be exceptional and $T_{i'+1,i'+1}$ is said to be touched from below by the
smallest candidate strictly below the $(i'+1)$-st row
if $i' \not= i''$.

\medskip

A {\it step of SPLIT~3} includes the choice of a protagonist $e \in C'$ 
and the performance of either Case~1, Case~2 or 
Case~3. A {\it step of MERGE~3} is defined analogously. The maximal vertical 
candidate $e$ at the beginning of a step of MERGE~3 is said to be the first protagonist
of the step and the candidate $e$ after the performance of the if-condition 
(If $c_{T}(e)=(h,h): [\dots]$) is said to be the second protagonist. Clearly the second
protagonist can be equal to the first protagonist in a step of MERGE~3.

\medskip

We first show that SPLIT~3 and MERGE~3 are inverse to each other for a certain 
class of pairs $(T,H)$.

\medskip

\begin{figure}
\setlength{\unitlength}{1cm}
\begin{center}
\scalebox{0.85}{\includegraphics{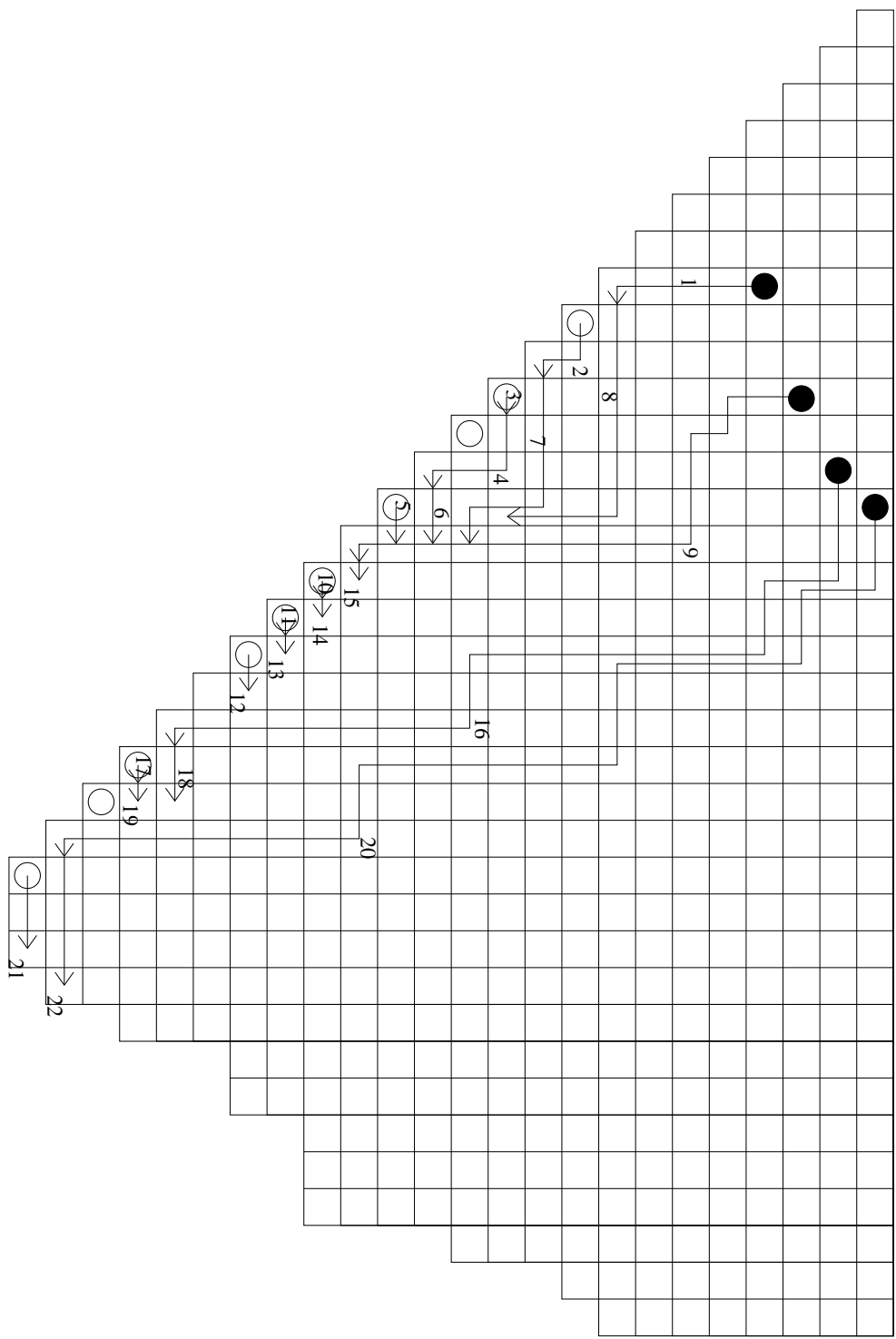}}
\end{center}
\caption{The application of SPLIT~3 to a irreducible pair. A full circle
  indicates a candidate in $C'$, an empty circle indicates a
  candidate in $C$ at the beginning of the application. The lines are the forward
  paths of jeu de taquin, the numbers near them indicate
  their order in the application of SPLIT~3. If the tip of an arrow touchs
  the border of a cell, this indicates that the entry in the cell is noticed
  to be unstable at that point in the application of the algorithm.}
\label{irred}   
\end{figure}

{\bf Irreducible pair with respect to SPLIT.} Let $(T,H)$ be a 
pair with property AFTER\_SPLIT~2. 
Suppose that in the course of applying SPLIT~3 to 
$(T,H)$ every protagonist which touches an entry from above was either 
in $C'$  at the beginning of  SPLIT~3 or is itself touched by 
an entry from above in the application of SPLIT. Observe that this is equivalent to the fact that we are
never in Case~3 in SPLIT~3 and that $C'\not=\emptyset$ at the beginning of
every step of SPLIT~3.
Then $(T,H)$ is said to be irreducible with 
respect to SPLIT. In Figure~\ref{irred} the application of SPLIT~3 to an
irreducible pair is illustrated.

\smallskip

\smallskip

{\bf Irreducible pair with respect to MERGE.} Let $T$ be a shifted 
tabloid which is ordered up to $(i,i)$ and $H$ a partial 
shifted hook tabloid, such that the subtabloid consisting
of the last $r-i+1$ rows is a shifted hook tabloid. If every
protagonist $e$ (except for possibly the last) in MERGE~3 either touches an entry from below
or terminates in row $z+1$ in a certain step of MERGE~3, then $(T,H)$ is said to
be irreducible with respect to MERGE. Note that this is equivalent
to the following: Suppose the second protagonist in a step of MERGE~3 is 
a horizontal candidate. Then this is either the last step in the application
of MERGE~3 to $(T,H)$ or this second
protagonist is equal to $e'$ in the if-condition at the beginning of the next step of
MERGE~3.

\medskip

{\bf Claim~3.} Let $(T,H)$ be irreducible with respect to SPLIT and let $(T',H')$ 
denote the pair we obtain after the application of SPLIT~3. Then $T'$ is
ordered up to $(i,i)$, the last $r-i+1$ rows of $H'$ form a shifted hook
tabloid and the application of MERGE~3 to 
$(T',H')$ results in $(T,H)$. 
 
\medskip

{\it Proof.} Let $C'_{start}$ denote the set $C'$ at the beginning 
of the application of SPLIT~3 to $(T,H)$. 
Suppose we are in the course of applying SPLIT~3: We 
choose in a certain step entry $e_1$ to be the protagonist and obtain the pair $(U,J)$,
in the next step we choose $e_2$ to be the protagonist and obtain $(U',J')$.
Let $z'$ be such that there exists a vertical candidate 
in row $z'$ of $(U,J)$ and in $C'_{start}$, which is not the protagonist in and before the step of $e_1$ and 
if there exists a vertical candidate in the $(z'+1)$-st row of $(U,J)$ and in $C'_{start}$, this 
entry was the protagonist in or before the step of $e_1$. 
If $z'$ does not exist, set $z'=i$.
Let $z''$ be the corresponding quantity for the step of $e_2$ and observe 
that $z''=z'-1$ iff $e_2$ is the vertical candidate of $(U,J)$ in the $z'$-th
row. 
Otherwise $z''=z'$. If we are in Case~2 in the step of $e_{k}$ let 
$h_{k}$ and $h'_{k}$ denote $h$ and $h'$ in the description of Case~2 of SPLIT~3, $k=1,2$.

\smallskip

By induction we assume the following: 

\begin{enumerate}

\item
Suppose we are in Case~1 in the step of $e_1$. Then 
$e_1$ is greater or equal than every vertical candidate strictly below the
$z'$-th row of $(U,J)$ with respect to the backward paths order and 
if $e_1$ is horizontal candidate then $e_2$ is the greatest vertical candidate strictly
below the $z''$-th row of
$(U,J)$. 

\item 
Suppose we are in Case~2 in the step of $e_1$. Then $e_2=U_{h_{1}',h_{1}'}$ is 
greater or equal than every vertical candidate strictly below the $z'$-th row of $(U,J)$ and 
if $e_2$ is a horizontal candidate in $(U,J)$ then $U_{h_{1}'-1,h_{1}'-1}$ is the greatest 
vertical candidate strictly below the $z'$-th row  of $(U,J)$.

\end{enumerate}

\smallskip

We aim to show that the application of one step of MERGE~3 to $(U',J')$ with
$z=z''$ either yields 

\begin{enumerate} 
\item $(U,J)$ or 
\item the pair we obtain before the performance of
$\trans$ in the step of $e_{1}$ (in which we are in Case~2 of SPLIT~3 and $e_{1}$ is a
horizontal candidate at the beginning of the step of $e_{1}$) or 
\item a pair which results in $(U,J)$ after the application
of a second step of MERGE~3 until the marker $(*)$ (see the description of MERGE~3).
\end{enumerate}

Furthermore we show the 
following (induction step): 

\begin{enumerate}
\item 
Suppose we are in Case~1 in the step of $e_2$. Then 
$e_2$ is greater or equal than every vertical candidate in $(U',J')$ strictly
below the $z''$-th row
and if $e_2$ is a horizontal candidate then $e_2=U_{k,k}$ for a $k$ and $U'_{k-1,k-1}=U_{k-1,k-1}$
is the greatest vertical candidate of $(U',J')$ strictly below the $z''$-th
row. 

\item
Suppose we are in Case~2 in the step of $e_2$ . 
In this case $U'_{h_{2}',h_{2}'}$ is greater or equal than 
every vertical candidate of $(U',J')$ strictly below the $z''$-th row 
of $(U',J')$ and if $U'_{h_{2}',h_{2}'}$ is a horizontal 
candidate then $U'_{h_{2}'-1,h_{2}'-1}$ is the greatest vertical candidate of
$(U',J')$ strictly below the $z''$-th row. 
\end{enumerate}

\smallskip

We distinguish between the two cases according to whether or not we are in Case~1 in 
the step of $e_1$.

\medskip

{\it Case~1.} There are two cases. Either

\begin{enumerate}
\item $e_1, e_2 \in C'_{start}$, 
$e_2$ is in the $z'$-th row of $U$ and $e_1$ was either the vertical candidate
of $C'_{start}$
in the $(z'+1)$-st row or the unique horizontal candidate of $C'_{start}$ in
$(T,H)$ or 
\item $c_U(e_2)=(k,k)$ for a $k$, $e_2$ is unstable in $U$ and 
the backward path of $e_1$ in $U$ contains $(k+1,k+1)$.
\end{enumerate}

\smallskip

{\it re 1.} Clearly $e_2$ is the greatest
vertical candidate in $(U,J)$ weakly below the $z'$-th row  for $e_1$
is greater or equal than every vertical candidate in $(U,J)$ strictly below
the  $z'$-th row 
by induction and $e_2$ is greater than $e_1$ by property AFTER\_SPLIT~2.
If $e_2$ changes row in the step of $e_2$, the first cell in row $z'+1$ in its forward path is
weakly right of the first cell in row $z'+1$ in $(U,J)$ in the backward path of $e_{1}$
by property AFTER\_SPLIT~2.
Thus $e_2$ is also the greatest vertical 
candidate in $(U',J')$ strictly below the $z''$-th row by Corollary~\ref{4}, 
if we are not in Case~2 in the step of $e_2$ and
in this case it is obvious that $e_2$ is the first and the second protagonist in the 
application of a step of MERGE~3 to $(U',J')$ with $z=z''$ (check that the if-condition
at the beginning of the step of MERGE~3 is either not fulfilled or its application leaves
$(U',J')$ unchanged in this step).
Otherwise $U'_{h_{2}',h_{2}'}$ is the greatest vertical candidate in $(U',J')$ 
strictly below the $z''$-th row  for firstly $e_2=U'_{h_{2}-1,h_{2}-1}$ is the greatest vertical candidate  
strictly below the $z''$-th row and strictly above the $h_{2}$-th row
of $(U',J')$ by Corollary~\ref{4}, secondly there exists no vertical candidate
strictly below the $h'_{2}$-th row and thirdly $U'_{h_{2}',h_{2}'}$ is greater than
$U'_{h,h}$ for $h_{2}-1 \le h < h_{2}'$. Thus
$U'_{h_{2}',h_{2}'}$ is
the first protagonist in the application of a step of MERGE~3 to $(U',J')$. 
Consequently $e_2$ is the second protagonist in the application 
of a step of MERGE~3 to $(U',J')$ if we are in Case~2 in the step of $e_2$.
It remains to show that the loop in the application of a step of MERGE~3 with $z=z''$ to 
$(U',J')$ terminates with 
`1. $e$ is in row $z+1$': Observe that the forward path of $e_2$ in $U$ does not contain 
a cell $(k,k)$ such that there exists a vertical candidate strictly below the
$k$-th row for such 
a vertical candidate would be greater than $e_2$ in $U$ by the arguments from Lemma~\ref{1}, 
which is a contradiction to the maximality of $e_2$ weakly below the $z'$-th
row in $(U,J)$.
Thus the loop does not stop with `2.'. 
Furthermore the 
forward path of $e_2$ in $(U,J)$ with respect to the appropriate set does not
contain a cell $(k,k)$ 
(except for possibly the endcell of the forward path of $e_{2}$ with respect
of the emptyset) with $h_{2}-1 \not= k$ such 
that the 
backward path in $U$ of a horizontal candidate contains $(k+1,k+1)$, for otherwise 
$e_2$ touches the minimal horizontal candidate whose backward paths contains $(k+1,k+1)$ from above in the application of
SPLIT by Lemma~\ref{6}, which implies $J_{i,k+1}=(k+1,k+1)$ and this is a
contradiction to the choice of $h_{2}$. Therefore the loop does 
not stop with `3.'. Finally, 
if $k$ is the row of $e_2$ in $U'$, then the subtabloid of $U$ consisting 
of the cells strictly below of the $z'$-th row and weakly above of the $k$-th
row is standard and consequently 
the loop does not stop with `4.'.
Therefore the application of one step of MERGE~3 to $(U',J')$ with $z=z''$ yields $(U,J)$ in this case.

\smallskip

{\it re 2.} There are two 
cases: either $e_1$ is a vertical candidate in $(U,J)$ or $e_1$ is a horizontal 
candidate in $(U,J)$.

Let $e_1$ be a vertical candidate, then $e_1$ is the greatest vertical candidate 
strictly below the $z'$-th row of $(U,J)$ by the induction hypothesis. If we are not in Case~2 in the step of $e_2$ then 
$e_2$ is the greatest vertical candidate strictly below the $z''$-th row 
of $(U',J')$ by 
Corollary~\ref{4} (since $e_2$ is unstable in $(U,J)$) and thus 
the first and the second protagonist in the application of a step of  
MERGE~3 to $(U',J')$ with $z=z''$.  If we are
in Case~2 in the step of $e_2$ then again $U'_{h_{2}',h_{2}'}$ is the greatest vertical candidate
strictly below the $z''$-th row of $(U',J')$ 
and with this the first protagonist in the application of a step of MERGE~ 3 to 
$(U,J)$. Thus $e_2$ is the second protagonist in the application of a step of MERGE~3 to 
$(U',J')$ with $z=z''$ in this case. It remains to show that the loop at the
end of the  application of a step of MERGE~3 to 
$(U',J')$ stops with `2.' and at the cell of $e_2$ in $U$. This is left to the reader, for it is similar to the 
previous case.

Now let $e_1$ be a horizontal candidate, then $e_1$ is greater than every vertical 
candidate strictly below the $z'$-th row of $(U,J)$ and $e_2$ is the greatest vertical 
candidate strictly below  the $z'$-th row of $(U,J)$ by the induction
hypothesis. Thus either $e_2$ 
(if we are in Case~1 in the step of $e_{2}$) or 
$U'_{h_{2}',h_{2}'}$ (if we are in Case~2 in the step of $e_{2}$) 
is the greatest vertical candidate strictly below the
$z''$-th row of $(U',J')$ 
and therefore the application of a step of MERGE~3 to $(U',J')$ with $z=z''$ yields $(U,J)$.
Note that the loop in MERGE~3 terminates with `3.' and at the cell of $e_2$ in
$U$, $e_{1}$ being the horizontal candidate in '3.'.

\medskip

{\it Case~2.} Clearly $e_2=U_{h_{1}',h_{1}'}$
and $e_2$ is greater or equal than every vertical candidate strictly below the
$z'=z''$-th row  of $(U,J)$ with respect to the backward paths order by the induction hypothesis. Thus either 
$e_2$ (if we are in Case~1 in the step of $e_{2}$) or $U'_{h_{2}',h_{2}'}$ (if
we are in Case~2 in the step of $e_{2}$) is greater or equal than every vertical candidate 
strictly below the $z''$-th row in $(U',J')$.

Suppose $e_2$ is a vertical candidate in $(U,J)$. Then either $e_2$ or 
$U'_{h_{2}',h_{2}'}$ is the first protagonist  
in the application of a step of MERGE~3 with $z=z''$ to
$(U',J')$ and therefore $e_2$ is the second protagonist. Thus the application
of a step of MERGE~3 
yields $(U,J)$. Observe that the loop in MERGE~3 terminates with `4.' and at
the cell of $e_2$ in $U$.

Suppose $e_2$ is a horizontal candidate in $(U,J)$. 
First we suppose that we are in Case~1 in the step of $e_{2}$.
Then $U'_{h_{1}'-1,h_{1}'-1}$ is the greatest vertical candidate strictly
below the $z''$-th row of $(U',J')$ and consequently the first protagonist 
in the application of a step of MERGE~3 to $(U,J)$. Thus $e'=e_{2}$ in the if-condition
at the beginning of this step of MERGE~3 and $e_{1}$
is the second (horizontal) protagonist. We obtain the pair which we had
immediately before the performance of $\trans$ in the step of $e_{1}$. 
Now suppose we are in Case~2 in the step of $e_{2}$. Clearly  $U'_{h_{2}'-1,h_{2}'-1}$ is the greatest vertical 
candidate strictly below the $z''$-th row of $(U',J')$ and consequently the 
first protagonist in the application of a step of MERGE~3 to $(U',J')$ with $z=z''$. The second 
protagonist is the horizontal candidate $e_2$ and since
$U'_{h_{1}'-1,h_{1}'-1}$ is the first protagonist in the next step of MERGE~3
we have  $e'=e_{2}$ in the if-condition at the beginning of  this next step. Thus the application of another step of MERGE~3 
until the marker $(*)$  yields $(U,J)$.  \qed

\medskip

{\bf Claim~4.} Let $(T,H)$ be irreducible with respect to MERGE and $(T',H')$ 
denote the pair we obtain after the application of MERGE~3 to $(T,H)$. Then $(T',H')$ 
has property AFTER\_SPLIT~2 and the application of SPLIT~3 to $(T',H')$ yields $(T,H)$.

\smallskip

{\it Proof.} In the course of applying MERGE~3 to $(T,H)$ we construct two sets
$C'$ and $C$. We start with $C=C'=\emptyset$. If the if-condition 
(If $c_T(e)=(h,h)$: [$\dots$]) at the beginning of a step of MERGE~3 is true
and we have just applied the appropriate commands, 
we set 
$C=C \cup \{T_{h'+1,h'+1},T_{h'+2,h'+2}, \dots, T_{h+1,h+1}\}$
and 
$C'=C' \setminus \{T_{h'+1,h'+1}, T_{h'+2,h'+2}, \dots, T_{h+1,h+1} \}$ if we 
are in the first case in these commands, and we set 
$C=C \cup \{T_{h'+1,h'+1}, T_{h'+2,h'+2}, \dots, T_{h,h} \}$ and 
$C'=C' \setminus \{ T_{h'+1,h'+1}, T_{h'+2,h'+2}, \dots, T_{h,h} \}$
if we are in the second case in these commands. After the if-condition and
its commands
we  let $C'= C' \cup \{ e \}$, whether or not the if-condition was fulfilled. 

\smallskip

Suppose that after a certain step in the application of MERGE~3 to $(T,H)$ 
we obtain the pair $(U',J')$, the sets $C'$ and $C$ and a row $z$. 
First we show that the second protagonist $e$ of this step is also the  
protagonist in the application of one step of SPLIT~3 to $(U',J')$ with the sets $C'$ and
$C$. This is equivalent to the fact that $e$ is the candidate in $C'$ whose row
is maximal. 

\smallskip

If the second protagonist is a horizontal candidate in a
certain step of MERGE~3, then it is the horizontal candidate $e'$ in the 
if-condition at the beginning of the next step (if we are not in the final step) of MERGE~3 and thus leaves $C'$ at that point, 
for $(T,H)$ is irreducible with respect to MERGE~3.
Consequently $C'$ includes at most one horizontal candidate and if this is the
case then 
all other candidates in $C'$ are strictly above the row of the unique horizontal candidate. 
(In order to see the second assertion in the previous sentence note that whenever $c_{T}(e)=(h,h)$ at the beginning of 
a step of MERGE~3 then there exists no vertical candidate strictly below the
$h$-th row by the maximality of $e$.)
Therefore it suffices to show 
the assertion for second protagonists $e$ that are vertical candidates.

\smallskip

In a step of MERGE~3 the second protagonist either terminates in row $z+1$ or in 
a cell $(k,k)$ on the main diagonal, if the protagonist is a vertical candidate. In the latter case by the maximality of the 
second protagonist
(observe that after the if-condition in a step of MERGE~3 the second protagonist is 
greater or equal than every vertical candidate strictly below the $z$-th row;
this is because the backward path of the first protagonist always contains the
cell of the second protagonist at the beginning of a step of MERGE~3)
every 
backward path of a vertical candidate strictly below the $k$-th row in the
tabloid we obtain after this step contains $(k+1,k+1)$ by the argument in the
proof of Lemma~\ref{1} and  
therefore these vertical candidates are greater than $e$ with respect to the backward paths order. However, 
the vertical second protagonist is 
still greater than every vertical candidate strictly below the $z$-th row and
strictly above the $k$-th row, 
since this part of the shifted hook tabloid remains unchanged in this step of
MERGE~3. Consequently in the course of MERGE~3 the second vertical
protagonists (and with this the elements of $C'$) are
arranged in  increasing rows (as long as the first protagonist is equal to the
second protagonist) until from time to time former second
vertical protagonists leave $C'$ after the performance of the commands in the
if-condition. 
We conclude that the second protagonist in a step of MERGE~3 either 
starts (by its maximality and the arguments in this paragraph) and ends 
(by the stopping criteria `4. $c_T(e)=(l,l)$ and $T_{l-1,l-1}$ is an unstable
vertical candidate' in MERGE~3 and the arguments in this paragraph)
below the lowest vertical candidate in $C'$ or is equal to the lowest
vertical candidate in $C'$, whose backward path in this step ends below the
second lowest vertical candidate in $C'$ (if it exists) by the fourth stopping criteria.
Consequently the 
row of the vertical second protagonist in a step of MERGE~3 is always below the other vertical 
candidates in $C'$ and by the observation about the horizontal candidates in
$C'$ in the previous paragraph
this proves our claim in the second paragraph of the proof.

\medskip

Let $(U,J)$ denote the pair before the step of MERGE~3, in which $e$ is the
second protagonist and let $C'_{before}$ and $C_{before}$ denote the accompanying sets. 
We have to show that if we apply one or two steps of SPLIT~3 to $(U',J')$ 
with the sets $C'$ and $C$, we reobtain $(U,J)$ together with $C'_{before}$ and $C_{before}$ 
or a pair of a shifted tabloid and a shifted hook tabloid together with the
two sets which  
we have obtained in the application of MERGE~3 to $(T,H)$ before the pair $(U,J)$.  
Proving this claim is now a matter of considering 
every circumstance in the step of $e$ in MERGE~3.
(Case distinction:
First suppose that we are in the  first case in the if-condition at the
beginning of the step of MERGE~3 (this leaves
no choice open for the rest of the step) and second
combine the second case in the if-condition and the case that the if-condition
is false with each of the four stopping 
criterions of the loop in MERGE~3.)
\qed

\smallskip 

Thus we have shown that for irreducible pairs SPLIT~3 and MERGE~3
are inverse to each other. Now we 
show the assertion for reducible pairs by induction with respect 
to the number of protagonists not contained in $C'$ at the beginning of the
application of SPLIT~3 which touch from above but are not touched from above in the application of SPLIT, 
respectively with respect to the number of second protagonists unequal to the
last second protagonist which 
neither touch an entry from below nor terminate in row $z+1$ in the
application of MERGE~3.

\medskip

\begin{center} \sc 6.7. Reducible pairs    \end{center}

\medskip

\begin{figure}
\setlength{\unitlength}{1cm}
\begin{center}
\scalebox{0.85}{\includegraphics{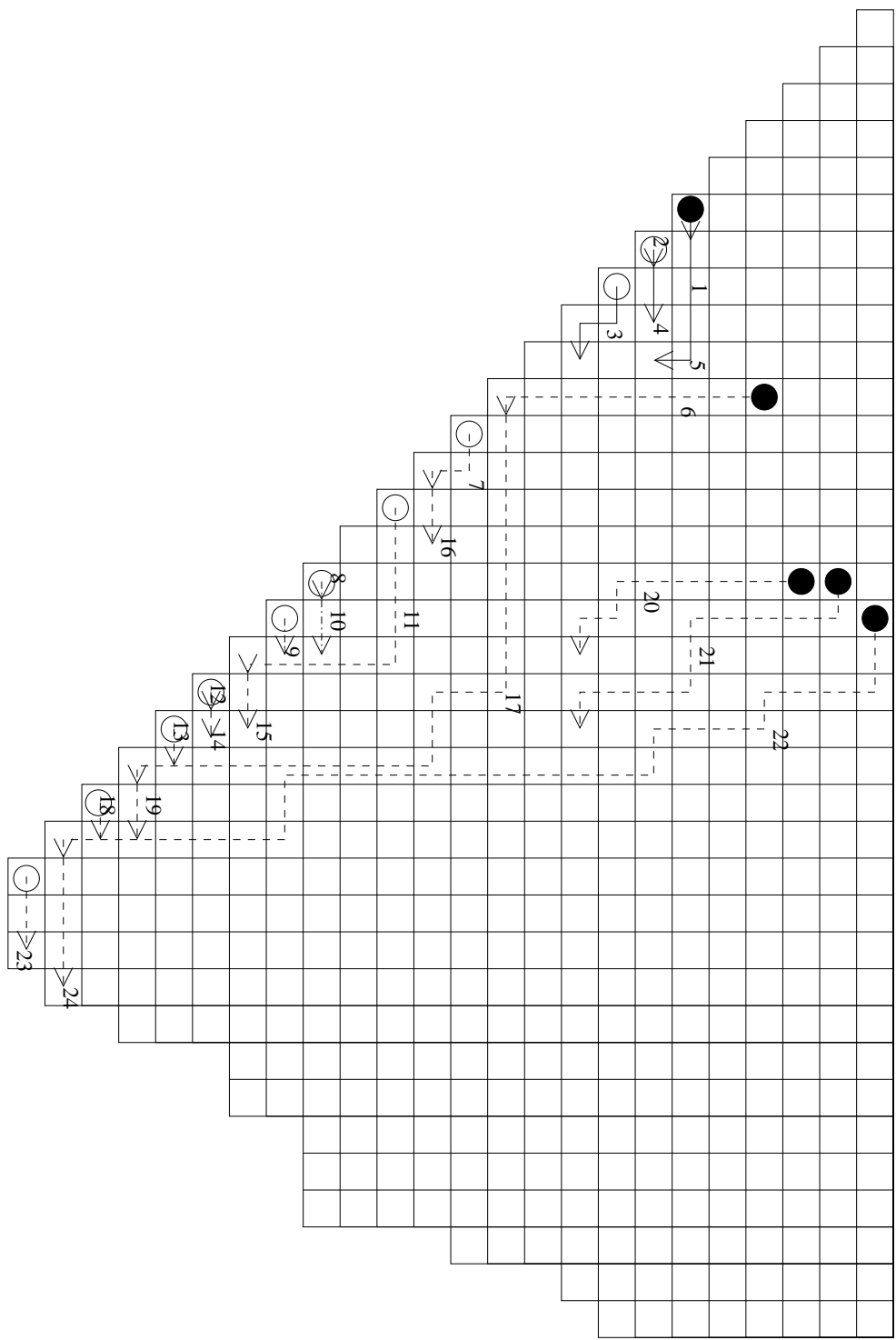}}
\end{center}
\caption{The application of SPLIT~3 to a reducible pair. (See caption of
  Figure~\ref{irred} for an explanation.) Observe that for distinct irreducible 'parts'
  different line styles are used.}
\label{red}   
\end{figure}

Observe that the pair in Example~\ref{best} in Section~\ref{examples} is reducible. In
our proof of Claim~5 we sometimes refer to that example. Moreover in
Figure~\ref{red} the application of SPLIT~3 to a reducible pair is illustrated.

\smallskip

\smallskip

For the proof of Claim~5 we need the following observation: Suppose we are in Case~3
of SPLIT~3 and let $(U,J)$ denote the current pair before the application 
of Case~3. Then the protagonist $e$ in the description of Case~3 of SPLIT~3 is a vertical candidate.
In order to show that suppose that $e$ was the first horizontal candidate in
the course of SPLIT~3 which contradicts this assertion.
Then we are either in the first step of SPLIT~3 or $e$ joined $C'$ in the 
previous step of SPLIT~3 in which we were in Case~2 and $e$ equals
$T_{h',h'}$ in the description of Case~2. 
We only consider the latter case here, for the first case is similar.
By the relative
position of $e$ and $e'$ in the current step of SPLIT~3 and the argument from
Lemma~\ref{1} the backward path of $e'$ in $(U,J)$ 
contains $(h'+1,h'+1)$. By Lemma~\ref{6} $e$ touches the smallest horizontal 
candidate $e''$ in $(U,J)$ whose backward path in $U$ contains $(h'+1,h'+1)$ from
above in the application of SPLIT~3 and this is a contradiction, since $e''$
is weakly above of the row of $e'$ by its minimality and by the fact that $e'$ is
on the main diagonal.

\medskip

{\bf Claim~5.} Let $(T,H)$ be a pair with property AFTER\_SPLIT~2
and $(T',H')$ the pair we obtain after the application of SPLIT~3. 
Then $T'$ is ordered up to $(i,i)$,
the subtabloid of $H'$ consisting of the last $r-i+1$ rows is a shifted 
hook tabloid
and the application of MERGE~3 to $(T',H')$ yields $(T,H)$.

\smallskip

{\it Proof.} Let $C'_{start}$ and $C_{start}$ denote $C'$ and $C$ at the beginning 
of the application of SPLIT~3 to $(T,H)$. We  
assume that every entry in $C_{start}$ is either touched from above 
or touches an entry from above in the application of SPLIT~3. By Claim~3 
we furthermore assume that $(T,H)$ is reducible with respect to SPLIT.

\smallskip

Let $e$ be the candidate in $C_{start}$ with maximal row $k$, which is not touched from above in the
application of SPLIT ($e=49$ in Example~\ref{best}). 
Furthermore, if it exists,  let $e' \in C_{start}$ be the entry with minimal row 
strictly below the $k$-th row such that $e'$ is touched from above in a
certain step in the application of SPLIT~3 by a protagonist which is at the beginning of the
step in a row strictly above of the $k$-th row ($e'=63$ in the example).  
Let 
$C_1$ denote the entries in $C_{start}$ which are strictly above the $k$-th row and
weakly 
below the row of $e'$ and $C_2$ denote the entries 
in $C_{start}$ which are strictly below the $k$-th row and strictly above the row of 
$e'$ in $T$.

\smallskip

Let SPLIT~3.1 denote the application of SPLIT~3 with $C'=C'_{start}$
and $C=C_1$ to $(T,H)$ and let SPLIT~3.2. denote the application 
of SPLIT~3 with $C'=\{e\}$ and $C=C_2$ to $(T,H)$. 
Observe that the application of SPLIT~3 with 
$C'=C'_{start}$ and $C=C_{start}$ to $(T,H)$ is equivalent to 
the simultaneous application of SPLIT~3.1 and SPLIT~3.2, 
where after every step of SPLIT~3.$i$, $i=1,2$, we can decide 
to make either the next step 
of SPLIT~3.1 or the next step of SPLIT~3.2, but whenever a protagonist in 
SPLIT~3.1 moves weakly below the $k$-th row the application of SPLIT~3.2 has already 
terminated. This is because $k$ is maximal and there exists a step of SPLIT~3 in which we are either in Case~3
with $e$ being the $e'$ in the description of Case~3 of SPLIT~3 or $C'=\emptyset$ at the
beginning of the step and $e$ is the candidate with minimal row in $C$ at that
point.

\smallskip 

Observe that after the termination of SPLIT~3.2 and before a protagonist of
SPLIT~3.1 moves weakly below the $k$-th row, $e$ is the greatest vertical 
candidate under the vertical candidates that arose in the application of
SPLIT~3.2 by the proof of Claim~3, for $e$ is the final protagonist in
SPLIT~3.2 for $k$ is maximal. 
The first protagonist of SPLIT~3.1 that moves weakly below of the $k$-th row 
is a vertical candidate (since whenever we are in 
Case~3 in SPLIT~3 the protagonist is a vertical candidate; this is proved at
the beginning of this subsection) and greater than $e$ with respect to the
backward paths order by Corollary~\ref{6}.  Moreover every subsequent protagonist of 
SPLIT~3.1, which moves weakly below the $k$-th row, is greater than the previous.
Thus the vertical candidates from SPLIT~3.1 dominate
the vertical candidates from SPLIT~3.2 once they come weakly below  $k$-th row.

Let $V_1$ denote the vertical pointers that came from 
SPLIT~3.1 and let $V_2$ denote the vertical pointers that 
came from SPLIT~3.2 after the application of SPLIT. By the induction 
hypothesis we reobtain 
$(T,H)$ if we apply MERGE~3 to the pointers in $V_1$ (this procedure is denoted by MERGE~3.1) and, 
seperately, to the pointers in $V_2$ (this procedure is denoted by MERGE~3.2) such that when we start
applying it to the pointers originated in $V_2$ the lowest vertical pointer that 
came from $V_1$ is strictly above the $k$-th row.  
If we apply MERGE~3 not seperately to the vertical pointers in 
$V_{1}$ and $V_{2}$, but to all vertical pointers in $(T',H')$ at
once, then the first steps of this application are equal to the first steps of MERGE~3.1
until the vertical pointers of MERGE~3.1 have moved strictly above the $k$-th
row, for the candidates originated in $V_{1}$ dominate the candidates
originated in
$V_{2}$ with respect to the backward paths order until they are strictly above
the $k$-th row as we saw in the previous paragraph.
\qed

\medskip

{\bf Claim~6.} Let $(T,H)$ be such that $T$ is ordered up to $(i,i)$, 
the subtabloid of $H$ consisting of the last $r-i+1$ rows is a shifted hook tabloid
and let $(T',H')$ 
denote the pair we obtain after the application of MERGE~3. 
Then $(T',H')$ has property AFTER\_SPLIT~2 and the application of SPLIT~3 to 
$(T',H')$ yields $(T,H)$.

\smallskip

{\it Proof.} 
By Claim~4 we assume that $(T,H)$ is reducible with respect to MERGE~3.
Let $e$ be the protagonist in the application of MERGE~3 
to $(T,H)$ which neither touches an entry from 
below nor terminates in row $z+1$, whose row $k$ of the place of change is
maximal.
We bipartite 
the vertical pointers of $(T,H)$ into two sets $V_1$ and 
$V_2$. The set $V_2$ includes the vertical pointer in $(T,H)$ that results in 
$e$ and if $p_1 \in V_2$ and the vertical pointer $p_2$ results in 
a protagonist that touches the protagonist that corresponds to $p_1$ from below then $p_2 \in V_2$.
All other vertical pointers are in $V_1$. 
Observe that by the algorithm  
a candidate associated with $V_{2}$ can be a protagonist only 
after the candidates associated with  $V_{1}$ have moved strictly above 
the $k$-th row or changed to horizontal candidates: By the proof of Claim~3
and 4 a
candidate $e'$ associated with $V_{2}$
and unequal to $e$ can only be a protagonist if $e$ was already a protagonist and
has moved to a row above of $e'$. Therefore we only have to show that $e$ is 
a protagonist only after the pointers associated with  $V_{1}$ have moved
strictly above the $k$-th row or changes to horizontal candidates. But this is obvious since if $e$ is a
protagonist in a step of MERGE~3
then every vertical candidate weakly below the $k$-th row belongs to
$V_{2}$. (After $e$ is protagonist for the first time and until it changes
into a horizontal candidate, $e$ is greater than every vertical candidate
strictly above of the current row of $e$ and strictly below the $z$-th row.
Moreover if there exist vertical candidates in this phase, which are strictly below the
current row of $e$ then they are greater than $e$ with respect to the backward
paths order and their places of change are  strictly below 
the current row of $e$. See the proof of Claim~3 and 4.)

We denote the application of MERGE~3 to $(T,H)$ restricted to the vertical pointers in 
$V_1$ by MERGE~3.1 and the application of MERGE~3 to $(T,H)$ restricted to the 
vertical pointers in $V_2$ by MERGE~3.2. Then the application of MERGE~3 is equivalent to 
the simultaneous application of MERGE~3.1 and MERGE~3.2, where MERGE~3.2 
starts after the vertical pointers that came from $V_1$ are strictly above of
the $k$-th row or have changed to horizontal candidates. 

\smallskip

By induction hypothesis we reobtain $(T,H)$ from $(T',H')$ if we apply SPLIT~3 simultaneously 
to the protagonists of MERGE~3.1 together with  the minimal horizontal candidates
$e'$ from the  if-condition  
(we denote this restricted application of SPLIT~3 by SPLIT~3.1) and 
to the protagonists of MERGE~3.2 together with  these horizontal candidates
(we denote this restricted application of SPLIT~3 by SPLIT~3.2) 
such that when a candidate in SPLIT~3.1  moves weakly below the
$k$-th row the application of SPLIT~3.2 has already finished. 
Note that $e$ is not touched from above in the application
of SPLIT~3 to $(T,H)$ (if we apply it to all candidates at once), 
for $e$ does not touch from below in the application of MERGE~3 to $(T,H)$.
Consequently if
we apply SPLIT~3 unrestricted to $(T,H)$, then this application starts 
with the application of SPLIT~3.1, until we are either in Case~3 of SPLIT~3 
with the $e'$ in the description of Case~3 in SPLIT~3 being equal to $e$ or SPLIT~3.1
has terminated. In the first case  
we perform SPLIT~3.2 and after the termination of SPLIT~3.2  we continue with the
application of SPLIT~3.1. In the
second case the set $C'$ is empty immediately after the termination of SPLIT~3.1 and $e$ is the candidate in $C$ with
minimal row. Thus we terminate with the application of SPLIT~3.2 in this case. 
\qed

\end{document}